\documentclass[reqno,12pt]{amsart}

\usepackage{geometry,amsmath,amsfonts,amssymb,amsthm,amscd,mathrsfs,graphicx,enumerate,multicol,wrapfig,subfigure,hyperref,stmaryrd, accents,emptypage}
\usepackage[all]{xy}
\usepackage[nodisplayskipstretch]{setspace}
\usepackage{setspace}
\usepackage[T1]{fontenc}

\newcommand{\ra}{\rightarrow}

\newcommand{\lra}{\longrightarrow}
\newcommand{\8}{\infty}
\newcommand{\p}{\prime}
\newcommand{\pt}{\partial}

\newcommand{\al}{\alpha}
\newcommand{\Om}{\Omega}
\newcommand{\om}{\omega}
\newcommand{\gam}{\gamma}
\newcommand{\Gam}{\Gamma}

\newcommand{\vp}{\varphi}
\newcommand{\lam}{\lambda}

\newcommand{\q}{\theta}

\newcommand{\be}{\beta}
\newcommand{\dt}{\delta}

\newcommand{\Zbb}{\mathbb{Z}}

\newcommand{\Cbb}{\mathbb{C}}

\theoremstyle{plain} 

\newtheorem{THM}{Theorem}[section]
\newtheorem{DEF}[THM]{Definition}
\newtheorem{EX}[THM]{Example}

\newtheorem{PROP}[THM]{Proposition}
\newtheorem{LEM}[THM]{Lemma}

\newtheorem*{CONJ1}{Conjecture}

\newcommand{\bt}{\bullet}

\newcommand{\Ac}{\mathcal{A}}

\newcommand{\Uc}{\mathcal{U}}

\newcommand{\Oc}{\mathcal{O}}

\newcommand{\Hom}{\mathrm{Hom}}

\newcommand{\Dc}{\mathcal{D}}

\newcommand{\Fc}{\mathcal{F}}

\newcommand{\Abb}{\mathbb{A}}

\newcommand{\Xfr}{\mathfrak{X}}
\newcommand{\Ec}{\mathcal{E}}

\newcommand{\Ufr}{\mathfrak{U}}

\newcommand{\Mfr}{\mathfrak{M}}

\newcommand{\Scl}{\mathcal{S}}
\newcommand{\Xc}{\mathcal{X}}

\newcommand{\Xscr}{\mathscr{X}}

\newcommand{\Uscr}{\mathscr{U}}

\showoutput
\usepackage{xcolor}
\definecolor{airforceblue}{rgb}{0.36, 0.54, 0.66}
\definecolor{burgundy}{rgb}{0.5, 0.0, 0.13}
\definecolor{majorelleblue}{rgb}{0.38, 0.31, 0.86}
\definecolor{darkblue}{rgb}{0.0, 0.0, 0.55}

\hypersetup{urlcolor=burgundy, linkcolor=burgundy,
citecolor=darkblue,
colorlinks=true}

\newcommand{\RNum}[1]{\uppercase\expandafter{\romannumeral #1\relax}}

\usepackage{chngcntr}
\counterwithin{section}{part} 

\counterwithout{equation}{section}

\title[Analytic and Alg. Defs SRS]
{Analytic and Algebraic Deformations of Super Riemann Surfaces
\\}
\author{\small Kowshik Bettadapura}
\date{}

\begin{document}

\begin{abstract} 
By analytic deformations of complex structures, we mean perturbations of the Dolbeault operator. By algebraic deformations of complex structures, we mean deformations of holomorphic glueing data. For complex manifolds there is, infinitesimally, a correspondence between these two types deformations. 
In this article we argue that an analogous correspondence holds between the analytic and algebraic deformations of a super Riemann surface.
\end{abstract}

\maketitle
\thispagestyle{empty}

\newpage
\thispagestyle{empty}

\setcounter{tocdepth}{1}
\tableofcontents

\onehalfspacing

\newpage
\section*{Introduction}

\subsection*{Super Riemann Surfaces}
Supersymmetry in physics is a proposed symmetry between bosons and fermions. From their quantum spin statistics, parameters directing the mechanics of bosons are commutative while those for fermions are anti-commutative. Mathematically, these parameters describe manifolds and supermanifolds respectively. We refer to \cite{LEI,BER,YMAN,QFAS,VARAD} for further historic commentary on the development and evolution of supersymmetry and supermanifolds, from a mathematical viewpoint. In superstring theory, Friedan in \cite{FRIEDAN} identified the super Riemann surface as the wordsheet to be swept out by the superstrings in superspacetimes. Once clearly defined, further work on the mathematical properties of super Riemann surfaces and their moduli were undertaken, e.g., in \cite{RABCRANE, ROSLY, ROTHMOD, FALQMOD}; in addition to some extensions of super Riemann surfaces to Riemann surfaces with higher degrees of supersymmetry in \cite{COHN,FALQ}. More recently, Donagi and Witten in \cite{DW1, DW2} have revived interest in this topic with some landmark results on the geometry of their moduli. Along the way, a certain pairing was derived between `analytic' deformations and supermoduli in \cite{DW2}. The author's attempt to understand this pairing forms the basis for the work undertaken in this article.
\subsection*{Donagi and Witten's Pairing}
Motivation for this article stems from a question posed by the author in the previous article \cite{BETTSRS} on the geometry of deformations of super Riemann surfaces. In \cite{BETTSRS}, super Riemann surfaces and their deformations are described as supermanifolds equipped with certain atlases, following similar descriptions in \cite{FRIEDAN, RABCRANE}. This description is referred to here as the \emph{algebraic viewpoint}. Donagi and Witten in \cite{DW2}, en route to their presentation of an invariant pairing formula over $\Mfr_g$, consider deformations of super Riemann surfaces as superconformal perturbations of the Dolbeault operator. This coincides with similar descriptions employed by D'Hoker and Phong in their work on superstring scattering amplitudes \cite{HOKERPHONG1, HOKDEFCPLX}. This description is referred to as the \emph{analytic viewpoint}. In an effort to better understand the pairing formula derived by Donagi and Witten, it was proposed in \cite{BETTSRS} that this formula extends appropriately to the algebraic deformations. We deduce such an extension in this article by establishing a correspondence between the algebraic and analytic viewpoints on deformations of super Riemann surfaces `to second order'. We follow work by Kodaira in \cite{KS} on a similar question concerning infinitesimal deformations of complex manifolds.
\subsection*{The Obstruction Class}
Upon establishing a correspondence between the analytic and algebraic viewpoints on super Riemann surfaces and deformations, we can obtain some insight on the splitting problem and the role of obstruction classes. Following early work in supermanifold theory by Berezin in \cite{BER} and Batchelor in \cite{BAT}, Green in \cite{GREEN} gave a general classification of supermanifolds via non-abelian sheaf cohomology. Here supermanifolds are characterised as being \emph{split} or \emph{non-split} and the cohomology groups housing obstruction classes to splitting, referred to as `obstruction spaces', were explicitly related to sheaves of non-abelian groups. In \cite{BETTPHD, BETTSRS} the total space of deformations (algebraic) of super Riemann surfaces were considered as a complex supermanifolds in their own right. As such, one can ask whether these total spaces will split and what their obstruction spaces might look like. From the analytic viewpoint however, it is unclear as to how this classification of complex supermanifolds might translate since, here, supermanifolds are smooth whence Batchelor's theorem in \cite{BAT} implies they split. Therefore, under our correspondence between algebraic and analytic deformations of super Riemann surfaces we can ask: \emph{what happens to the obstruction classes to splitting the total space of algebraic deformations?} We find that the \emph{primary} obstruction class will, under the Dolbeault isomorphism, form part of the complex structure for an analytic deformation.

\subsection*{Applications and Outlook}
With the correspondence between analytic and algebraic deformations of super Riemann surfaces we can: (1) adapt and address the splitting problem for analytic deformations; and (2) translate Donagi and Witten's pairing between analytic deformation parameters and supermoduli to the algebraic case. We conclude the article with one further comment regarding the `higher obstruction theory' of super Riemann surface deformations. In \cite{BETTSRS}, it was observed that if the \emph{primary} obstruction to splitting a super Riemann surface deformation vanishes, then the deformation will be split to second order. That is, there will not exist an exotic, third order obstruction class. Subsequently, it was conjectured: \emph{if the primary obstruction to splitting an algebraic deformation of a super Riemann surface vanishes, then the total space of the deformation must itself be split}. We refer to this as the `No-Exotic-Obstructions' conjecture. With Donagi and Witten's pairing adapted to the algebraic case, in addition to their results on the non-projectedness of the supermoduli space in \cite{DW1,DW2}, we speculate as to what a proof of the this conjecture might look like.

\section*{Outline}

\noindent
This article is divided into four parts. We begin in Part \ref{kfnkbrvuybv94oi3n3po} with some preliminary theory on deformations of complex manifolds, supermanifolds and super Riemann surfaces. Key results presented here include Theorem \ref{rfjrnvonoieep} on the relation between the Kodaira-Spencer class and perturbations of the Dolbeault operator; Theorem \ref{34fiufguhfojp3} on obstruction classes of supermanifolds and Theorem \ref{h8937893h0fj309j03} on thickenings over Riemann surfaces; Lemma \ref{rfnoifnoi3fip3mpm3} and \ref{rfnoifnoi3feem3} on the structure of superconformal maps and vector fields parametrised by $\Abb^{0|n}_\Cbb$; and finally Lemma \ref{445f4g4g4g4} and Theorem \ref{fjbvbviueovinep} on the structure of super Riemann surfaces. In Part \ref{rcjoehc894h984h48} we look at the analytic and algebraic deformations of a super Riemann surface. Definitions of each kind are presented. See Definition \ref{rf4fg784hf98h0f3} for the algebraic case and Definition \ref{rf4fg70f324fh894hg9843445} for the analytic case. The respective deformation parameters up to equivalence are described in Theorem \ref{fhvurh84j093j3j33333} in the analytic case; and in Theorem \ref{rfhf9h3f8309fj3fjoooo3} in the algebraic case. We conclude with the statement of the main result of this article on the correspondence between these deformations in Theorem \ref{rfg748gf7f983h04i4i4i43hf03} and give the outline of its proof. Part \ref{lfrpovviovuoh} is devoted to the proof Theorem \ref{rfg748gf7f983h04i4i4i43hf03}. We draw on a number of the key results presented earlier in the article and aim therefore to give a self-contained proof. 
We conclude this article with Part \ref{ppjpvburbviyvbiuvr} where some applications of the correspondence in Theorem \ref{rfg748gf7f983h04i4i4i43hf03} are given. In Theorem \ref{jf0jg43r3jgp3345topgj4gj4} we consider how the the splitting problem might adapt to analytic deformations; and Theorem \ref{jf04jg0jg4jgpo4jopgj4gj4} concerns the translation of Donagi and Witten's invariant pairing formula to the algebraic case. The article concludes with some comments on the No-Exotic-Obstructions conjecture.

\newpage
\section*{List of Notation}

\noindent
In order of appearance:
\begin{align*}
\begin{tabular}{l l l}
Symbol & & Description\\
\hline
$X$& &Complex manifold;
\\
$X^\8$ & & Smooth manifold underlying $X$;
\\
$\overline \pt$ & & Dolbeault operator on $X^\8$;
\\
$\Xscr$& & Total space of a deformation of $X$;
\\
$\Uscr = (\Uscr_\al)_{\al\in I}$& & Open covering of $\Xscr$;
\\
$(\Uscr, F)$ & & Atlas for $\Xscr$;
\\
$(U, f)$ & & Atlas for $X$;
\\
$\Ac^{p, q}(X^\8, TX^\8)$ & & Smooth, $TX^\8$-valued $(p, q)$-forms;
\\
$(\Xscr^\8, \overline \pt_t)$ & & $\Xscr$ as a smooth manifold with complex structure;
\\
$\Xfr = (X, \Oc_\Xfr)$& & Supermanifold over $X$;
\\
$T^*_{X, -}$& &Odd cotangent bundle;
\\
$S(X, T^*_{X, -})$& &Split model;
\\
$(X, T^*_{X, -})$& & Model;
\\
$(\Ufr, \rho)$& & Atlas for $\Xfr$;
\\
$\om_{(\Ufr, \rho)}$& &Obstruction to splitting $(\Ufr, \rho)$;
\\
$\Abb^{m|n}_\Cbb$ & & Complex affine superspace;
\\
$\Cbb^{m|n}$& & Complex Euclidean superspace;
\\
$\Cbb^{1|1}_{(x|\q)}$& &Complex Euclidean space with coordinate system $(x|\q)$;
\\
$D_{(x|\q)}$& &Superconformal generator on $\Cbb^{1|1}_{(x|\q)}$;
\\
$(\xi_1, \ldots, \xi_n)$ & & Parameters on $\Abb^{0|n}_\Cbb$;
\\
$\Dc$& &Superconformal structure;
\\
$\Scl = (\Xfr, \Dc)$ & & Super Rieman surface;
\\
$C$, $C^\8$ & & Riemann surface; underlying smooth manifold
\\
$T_{C, -}^*$ & & Spin structure on $C$;
\\
$\Xc$& &Deformation of $S(C, T^*_{C, -})$;
\\
$(\Uc, \vartheta)$ & & Algebraic deformation (atlas for $\Xc$);
\\
$\Theta(T^*_{\Xc, -})$& &Odd Kodaira-Spencer class of $\Xc$;
\\
$\widetilde\pt$& & Dolbeault operator on supermanifolds;
\\
$\Xc^\8, S(C, T^*_{C, -})^\8$& &Smooth supermanifolds underlying $\Xc, S(C, T^*_{C, -})$;
\\ 
$(\Xc^\8, \widetilde\pt_\xi)$& &Analytic family or deformation;
\\
$\Xc^{\mathrm{an.}}$& & Complex supermanifold associated to $(\Xc^\8, \widetilde\pt_\xi)$;
\\
$\Mfr_g$& &Supermoduli space of genus $g$ curves.
\end{tabular}
\end{align*}

\newpage
\part{Preliminaries}
\label{kfnkbrvuybv94oi3n3po}

\section{Deformations of Complex Structures}

\noindent
We review here some relevant material on the deformations of complex manifolds from the text by Kodaira in \cite{KS}. The main result (Theorem \ref{rfjrnvonoieep}) ought to viewed as a guiding precedent for the main result derived in this article.

\subsection{Preliminary Descriptions}

\subsubsection{The Analytic Viewpoint}
Let $X$ be a complex manifold. Note that it can also be viewed as a smooth manifold $X^\8$ with a choice of complex structure $\overline \pt$.\footnote{Here $\overline \pt$ is the Dolbeault operator. A smooth function $f$ is holomorphic if and only if $\overline \pt f = 0$. As we can decide on the holomorphicity of smooth functions in this way, we can therefore identify $\overline \pt$ itself with a complex structure on $X^\8$. The operator $\overline \pt$ can also be related to an almost complex structure $J$ on $X^\8$. Integrability of $J$ to a complex structure amounts to flatness of $\overline \pt$, i.e., $\overline \pt^2=0$. We refer to \cite[p. 108]{HUYB} for further details on this correspondence.} A deformation of $X$ over a base $B$ with distinguished point $b_0\in B$  is a complex manifold $\Xscr$ equipped with a projection (or, surjective submersion) $\pi : \Xscr\ra B$ whose fibers $\pi^{-1}(b)\subset \Xscr$ are complex submanifolds; and $\pi^{-1}(b_0)\cong X$. The complex structure of the fibers vary with parameters on the parameter space $B$. The smooth structure is taken to be constant which means $\Xscr$ and $X\times B$ are diffeomorphic. Thus $\Xscr \ra B$ captures a variation of the complex structure $\overline \pt$ on a fixed, smooth manifold $X^\8$. This viewpoint with $X$ a smooth manifold with complex structure $(X^\8, \overline \pt)$ will be referred of as the `analytic viewpoint'. 

\subsubsection{The Algebraic Viewpoint}
The deformation $\Xscr$ can also be described without any reference to the underlying, smooth manifold $X^\8$ and complex structure $\overline \pt$. This description was the first proposal by Kodaira and Spencer, as explained in \cite[p. 182]{KS}, of a `variation of complex structure'. It references the glueing data of the complex manifold $X$ directly. Let $\Uscr = (\Uscr_\al)$ be a covering of $\Xscr$ with transition data $(F_{\al\be} : \Uscr_{\al\be} \stackrel{\cong}{\ra}\Uscr_{\be\al})$ identifying intersections.\footnote{Here the open sets $(\Uscr_\al)$ cover $\Xscr$ with $\Uscr_\al\cong \Cbb^{\dim X}$; and on choosing distinguished subspaces $\Uscr_{\al\be}\subset \Uscr_\al$ for all pairs $\al, \be$ we recover $\Xscr$ as the equivalence class $\Xscr \cong (\bigsqcup_\al \Uscr_\al)/(F_{\al\be})$, where $(F_{\al\be} : \Uscr_{\al\be} \stackrel{\cong}{\ra}\Uscr_{\be\al})$ are biholomorphisms between the distinguished subspaces.}
The tuple $((\Uscr_\al), (\Uscr_{\al\be}), (F_{\al\be}))$, denoted by $(\Uscr, F)$ in shorthand, is referred to as an \emph{atlas} for $\Xscr$. 
We take $B = \Cbb^n$ with $b_0 = {\bf 0}$ and global coordinate system $t$.
As with $\Xscr$, let $(U, f) = ((U_\al), (U_{\al\be}), (f_{\al\be}))$ be an atlas for $X$ defining its complex structure. With $\pi: \Xscr\ra (\Cbb^m, {\bf 0})$ we require:
\begin{enumerate}[$\bt$]
	\item for each $\al$, $\pi^{-1}({\bf 0})\cap \Uscr_\al\cong U_\al$;
	\item for each $\al, \be$, $\pi^{-1}({\bf 0})\cap \Uscr_{\al\be}\cong U_{\al\be}$;
	\item a commutative diagram for each $\al,\be$:
	\[
	\xymatrix{
	\ar[d] U_{\al\be} \ar[r]^{f_{\al\be}} & U_{\be\al}\ar[d]
	\\
	\Uscr_{\al\be} \ar[r]^{F_{\al\be}}& \Uscr_{\be\al}
	}
	\]
	where the vertical arrows are the inclusions $U_{\al\be}\cong \pi^{-1}({\bf 0})\cap \Uscr_{\al\be}\subset \Uscr_{\al\be}$.
\end{enumerate}
The above bullet-points justify: on $\Xscr$ one has local coordinates $x_t$ on $\Uscr_\al$ and $y_t$ on $\Uscr_\be$. At $t = {\bf 0}$ we obtain local coordinates on $X$.  For a neighbourhood $V$ of $b_0 = {\bf 0}$, the transition data for $\Xscr$ on $\pi^{-1}(V) \cap \Uscr_\al\cap \Uscr_\be$ is given by,\footnote{more explicitly, $F_{\al\be}|_{\pi^{-1}(V)} : \pi^{-1}(V)\cap \Uscr_{\al\be} \stackrel{\cong}{\ra} \pi^{-1}(V)\cap \Uscr_{\be\al}$}  
\begin{align}
y_t
&=
F_{\al\be}(x_t, t)
\label{fjvbkbvirbviubveubo}
\\
&=
f_{\al\be}(x)
+
\sum^\8_{|I|>0} f_{\al\be}^I(x)~t_I 
\label{vh7gf794gf98hf083h093jf0}
\end{align}
where $I$ is a multi-index and $|I|$ its length. For $t = (t_1, \ldots, t_m)$ and $I = (i_1, \ldots, i_k)$ we write $t_I = t_{i_1,\ldots, i_k} = t_{i_1}\cdots t_{i_k}$. The viewpoint of deformations $\pi:\Xscr\ra (B, {\bf 0})$ as variations of the holomorphic glueing data of the central fiber $X = \pi^{-1}({\bf 0})$ in \eqref{fjvbkbvirbviubveubo} will be referred to as the `algebraic viewpoint'. 
\\\\
{\bf Remark.} It is a difficult problem as to whether the expression in \eqref{vh7gf794gf98hf083h093jf0} can be chosen on all intersections so as to be convergent. Treating \eqref{vh7gf794gf98hf083h093jf0} formally however, one can bypass this problem and obtain a preliminary classification.

\subsection{Vector Dolbeault Forms and the Kodaira-Spencer Class}

\subsubsection{Analytic Viewpoint (Vector Dolbeault Forms)}
Let $X = (X^\8, \overline \pt)$ be a complex manifold and denote by $\Ac^{p, q}(X^\8)$ the space smooth $(p, q)$-forms on $X^\8$.\footnote{That is, the smooth sections of the vector bundle of $(p, q)$-forms.} Let $\Ac^{p, q}(X^\8,TX^\8)$ denote the space of smooth, vector $(p, q)$-forms (i.e., $(p, q)$-forms valued in the tangent bundle $TX^\8$). The total space $\Xscr$ of a deformation $\Xscr\ra (\Cbb^m, {\bf 0})$ of $X$ is itself a complex manifold and so can be represented by a smooth manifold with complex structure $(\Xscr^\8, \overline \pt_t)$. The pair $(\Xscr^\8, \overline \pt_t)$ will be referred to as an \emph{analytic deformation}. Recall now that $\Xscr^\8$ and $X^\8\times \Cbb^m$ are diffeomorphic so that $T\Xscr^\8 \cong TX^\8\times T\Cbb^m$. Kodaira in \cite[p. 255]{KS} proves:
\begin{center}
\begin{minipage}[t]{0.7\linewidth}
\begin{PROP}\label{r84hg98h48g409g3} 
Let $\pi: (\Xscr^\8, \overline \pt_t)\ra (\Cbb^m, {\bf 0})$ be an analytic deformation of $X = (X^\8, \overline \pt)$. Then over a sufficiently small open neighbourhood $V\subset (\Cbb^m, {\bf 0})$ there exists a vector $(0, 1)$-form $\vp(t)\in \Ac^{0, 1}(X^\8, TX^\8)$, holomorphic in $t$, such that 
\begin{enumerate}[(i)]
	\item \[\overline \pt_t |_{\pi^{-1}(V)}= \overline \pt - \vp(t);\]
	\item  $\vp({\bf 0}) = 0$ and;
	\item $\vp(t)$ satisfies the Mauer-Cartan equation
	\[
	\overline \pt \vp(t) - \frac{1}{2}\left[\vp(t), \vp(t)\right] = 0.
	\]
\end{enumerate}
\qed
\end{PROP}
\vspace{\dp0}
\end{minipage}
\end{center}
A function $g$ on $\Xscr^\8$ is then holomorphic on the fiber over $t_0\in V\subset (\Cbb^m, {\bf 0})$ if and only if  $\overline \pt_tg|_{t=t_0} = 0$. We see in Proposition \ref{r84hg98h48g409g3} that this vector $(0, 1)$-form $\vp(t)$ essentially characterises an analytic deformation $\pi : (\Xscr^\8, \overline \pt_t)\ra(\Cbb^m, {\bf 0})$. An important consequence of the Mauer-Cartan equation satisfied by $\vp(t)$ in Proposition \ref{r84hg98h48g409g3}\emph{(iii)} is the following:
\begin{center}
\begin{minipage}[t]{0.7\linewidth}
\begin{PROP}\label{ruguyeiueihf89h389} 
The vector $(0, 1)$-form $\vp(t)$ in Proposition \ref{r84hg98h48g409g3} satisfies:
\[
\overline \pt \left(\left.\frac{\pt \vp}{\pt t}\right|_{t=0}\right) = 0.
\]
\qed
\end{PROP}
\vspace{\dp0}
\end{minipage}
\end{center}
Hence, associated to any analytic deformation $(\Xscr,^\8, \overline\pt_t)\ra (\Cbb^m, {\bf 0})$ is a representative of a vector Dolbeault form $\pt\vp/\pt t|_{t=0}$. That is, $\pt\vp/\pt t|_{t=0}$ defines a class $\om_{(\Xscr^\8,\overline \pt)}$ in the Dolbeault cohomology group $H^{0, 1}_{\overline \pt}(X^\8, TX^\8)$. We will refer to $\om_{(\Xscr^\8,\overline \pt_t)}$ as the \emph{Dolbeault form associated to $(\Xscr^\8, \overline \pt_t)$}. 

\subsubsection{Algebraic Viewpoint (Kodaira-Spencer Class)}
Let $\Xscr\ra(\Cbb^m, {\bf 0})$ be a deformation of a complex manifold $X$ with $(\Uscr, F)$ an atlas for $\Xscr$. Denote coordinates on $\Uscr_\al$ resp., $\Uscr_\be$ by $x_t$ resp., $y_t$. On the intersections $\Uscr_\al\cap\Uscr_\be$ we have the relation $y_t = F_{\al\be}(x_t, t)$ from \eqref{fjvbkbvirbviubveubo}. A choice of atlas $(\Uscr, F)$ for $\Xscr$ will be referred to as an \emph{algebraic deformation}. Since $F = (F_{\al\be})_{\al, \be\in I}$ are the transition functions of an atlas, they satisfy the cocycle condition on triple intersections $\Uscr_\al\cap\Uscr_\be\cap \Uscr_\gam$:
\begin{align}
F_{\al\gam}(x_t, t) = F_{\be\gam}\left(F_{\al\be}(x_t, t), t\right).
\label{fkvnrb49h9fio3jfp3}
\end{align}
Expanding $F$ in the (formal) infinite sum \eqref{vh7gf794gf98hf083h093jf0} over each intersection and imposing \eqref{fkvnrb49h9fio3jfp3} on triple intersections, Kodaira in \cite[p. 198]{KS} was led to the following result serving to classify deformations of complex manifolds infinitesimally:
\begin{center}
\begin{minipage}[t]{0.7\linewidth}
\begin{PROP}\label{uniunioiojio} 
Let $(\Uscr, F)$ be an algebraic deformation of a complex manifold $X$ over $(\Cbb^m, {\bf 0})$; and let $T_X$ denote the sheaf of holomorphic vector fields on $X$. Then there exists a class $\om_{(\Uscr, F)}\in \mbox{\emph{\v H}}^1(\Uscr, T_X)$ represented on intersections $\Uscr_\al\cap \Uscr_\be$ by the $1$-cocycle:
\[
(\om_{(\Uscr, F)})_{\al\be}
=
\left.\frac{\pt F_{\al\be}}{\pt t}\right|_{t=0}\otimes\frac{\pt}{\pt y}
\]
where $y$ is the complex coordinate on the open set $U_\be = \Uscr_\be\cap\pi^{-1}({\bf 0})\subset X$.\qed
\end{PROP}
\vspace{\dp0}
\end{minipage}
\end{center} 
To eliminate dependence on the choice of atlas in Proposition \ref{uniunioiojio} we can limit over the common refinement of coverings to get isomorphisms\footnote{see e.g., \cite[Lemma 5.2, p. 210]{KS}} $\mbox{\v H}^1(\Uscr, T_X)\cong \mbox{\v H}^1(X, T_X)\cong H^1(X, T_X)$. The image of $\om_{(\Uscr, F)}$ in $H^1(X, T_X)$ under these isomorphisms is referred to as the \emph{Kodaira-Spencer class} of the deformation $\Xscr$.

\subsection{Analytic and Algebraic Deformations}

\subsubsection{The Vector Dolbeault Theorem}
In Proposition \ref{ruguyeiueihf89h389} and \ref{uniunioiojio} respectively we obtained cohomology classes classifying the infinitesimal deformations of a complex manifold $X$. The classification in the former used Dolbeault cohomology while the latter used sheaf cohomology. Generally, there exists a natural isomorphism between these two cohomologies, reviewed in a number of classical texts on complex manifold theory (e.g., \cite[p. 110]{HUYB}). The statement of the isomorphism in considerable generality is as follows: let $E$ be a holomorphic vector bundle on a complex manifold $X$ and denote by $\Ec$ its sheaf of holomorphic sections. We can view $E$ as a smooth vector bundle $E^\8$ on $X^\8$ with a flat connection $\overline \pt_E$ (see \cite[p. 110, Theorem 2.6.26]{HUYB}). Let $\Ac^{p, q}(X^\8, E^\8)$ be the space of smooth sections of $(p, q)$-forms on $X^\8$ valued in $E^\8$. Flatness of $\overline \pt_E$ ensures that $(\Ac^{p, q}(X^\8, E^\8), \overline\pt_E)$ will be a differential complex with cohomology the $E^\8$-valued, Dolbeault cohomology group $H^{p, q}_{\overline \pt}\left(X^\8, E^\8\right)$. The Dolbeault isomorphism is then the isomorphism between the following cohomology groups:
\begin{align}
Dol : H^q(X, \Om^p_X\otimes \Ec) 
\stackrel{\cong}{\lra}
H^{p, q}_{\overline \pt}\left(X^\8, E^\8\right)
\label{rhf74g9h3f830fj03333333}
\end{align}
where $\Om^p_X$ is the sheaf of holomorphic $p$-forms on $X$.

\subsubsection{Correspondence of Infinitesimal Deformations}
We wish to apply \eqref{rhf74g9h3f830fj03333333} to the setting $p = 0, q = 1$ and $E^\8 = TX^\8$. In this case see that we have an isomorphism of cohomology groups,
\begin{align}
Dol: 
H^1(X, T_X)
\stackrel{\cong}{\lra}
H^{0, 1}_{\overline\pt}(X^\8, TX^\8).
\label{rhf74g9h3f830fj0333333333}
\end{align}
The range of $Dol$ in \eqref{rhf74g9h3f830fj0333333333} contains the Dolbeault forms of deformations of $X$ by Proposition \ref{ruguyeiueihf89h389}; while the domain of $Dol$ in \eqref{rhf74g9h3f830fj0333333333} contains the Kodaira-Spencer classes of deformations of $X$ by Proposition \ref{uniunioiojio}.\footnote{Not every such class need to come from an infinitesimal deformation however.}
In \cite[p. 256]{KS} we find a proof of the following:
\begin{center}
\begin{minipage}[t]{0.7\linewidth}
\begin{THM}\label{rfjrnvonoieep} 
Fix a complex manifold $X$ and base $(\Cbb^m, {\bf 0})$. Then to any analytic deformation $(\Xscr^\8, \overline \pt_t)$ of $X$ over $(\Cbb^m, {\bf 0})$  there exists an algebraic deformation $(\Uscr, F)$ of $X$ over $(\Cbb^m, {\bf 0})$ such that 
\[
Dol~\om_{(\Xscr^\8, \overline \pt_t)} = \om_{(\Uscr, F)}
\]
and vice-versa.
\qed
\end{THM}
\vspace{\dp0}
\end{minipage}
\end{center} 
Since the classes $\om_{(\Xscr^\8, \overline \pt_t)}$ and $\om_{(\Uscr, F)}$ classify the infinitesimal deformations of $X$ over $(\Cbb^m, {\bf 0})$, we refer to the correspondence in Theorem \ref{rfjrnvonoieep} as a correspondence of infinitesimal deformations.

\section{Complex Supermanifolds}

\subsection{Preliminaries}\label{knebvou4bvoioi3nvo3}
Supermanifolds are a very mild kind of non-commutative manifold. They are spaces whose function algebra is not strictly commutative, but rather weakly commutative or, `supercommutative'. Much like the locally-ringed-space-definition of smooth or complex manifolds, a smooth or complex supermanifold is defined as a locally ringed space $\Xfr = (X, \Oc_\Xfr)$, with structure sheaf $\Oc_\Xfr$ locally isomorphic to an exterior algebra---the prototypical supercommutative algebra. To clarify the role of smoothness, an equivalent definition of a supermanifold is as a locally ringed space $(X, \Oc_\Xfr)$ where $X$ is smooth or complex; and the structure sheaf $\Oc_\Xfr$ is locally isomorphic to the exterior algebra of a smooth or holomorphic vector bundle $\wedge^\bt T^*_{X -}$. In this definition, we say $\Xfr$ is \emph{modelled on the pair $(X, T^*_{X, -})$}. Here $X$ and $T_{X, -}^*$ are referred to as the \emph{reduced space of $\Xfr$}  and the \emph{odd cotangent bundle on $X$} respectively. The prototypical supermanifold is the locally ringed space $S(X, T^*_{X, -})\stackrel{\Delta}{=} (X, \wedge^\bt T^*_{X, -})$. It is referred to as the \emph{split model}.
\begin{center}
\begin{minipage}[t]{0.7\linewidth}
\begin{DEF}\label{tvirri4f4f4f44orjpjrpov} 
\emph{A supermanifold $\Xfr$ modelled on $(X, T^*_{X -})$ is said to be \emph{split} if it is isomorphic to the split model $S(X, T^*_{X, -})$. Otherwise $\Xfr$ is non-split. A \emph{splitting} of $\Xfr$ is a choice of isomorphism (if one exists) between $\Xfr$ and the split model, $S(X, T^*_{X, -})$.}
\end{DEF}
\vspace{\dp0}
\end{minipage}
\end{center}
There are generally obstructions for a supermanifold to be split. However, as illustrated in \cite{BER, DW1, BETTHIGHOBS}, these obstructions might be `exotic'---i.e., falsely claim there is no splitting while in fact there might be one. Hence, a choice of atlas on $\Xfr$ does not necessarily establish a unique `complex supermanifold structure'. This problem is circumvented in `degree two' however and the obstruction class there does indeed define an invariant of the supermanifold structure on $\Xfr$, as observed by Berezin in \cite[Theorem 4.6.2., p. 158]{BER}; and employed by Donagi and WItten in \cite{DW1} to deduce their non-projectability result for the supermoduli space of curves. 
We collect these remarks on the obstructions to splitting in what follows.\footnote{\label{ftnoet3646f6r4}In Theorem \ref{34fiufguhfojp3} we use the notation: $T^*_{X, (\pm)^j} = T^*_X$ or $T^*_{X, -}$ if $j$ is even or odd respectively.} 
\begin{center}
\begin{minipage}[t]{0.7\linewidth}
\begin{THM}\label{34fiufguhfojp3} 
Let $\Xfr$ be a complex supermanifold modelled on $(X, T^*_{X, -})$. Then any atlas $(\Ufr, \rho)$ for $\Xfr$ defines an obstruction class $\om_{(\Ufr,\rho)}$ in the cohomology group $H^1(X, \mathcal Hom_{\Oc_X}(T^*_{X, (\pm)^j}, \wedge^jT^*_{X, (\pm)^j}))$ for some $j \in\{ 2, \ldots, \mathrm{rank}~T^*_{X, -}\}$.
If $j = 2$ and $\om_{(\Ufr,\rho)}\neq0$, then $\Xfr$ is non-split.\qed
\end{THM}
\vspace{\dp0}
\end{minipage}
\end{center}
As remarked in Theorem \ref{34fiufguhfojp3}, any atlas $(\Ufr, \rho)$ for a supermanifold $\Xfr$ will define an obstruction class $\om_{(\Ufr,\rho)}\in H^1(X, \mathcal Hom_{\Oc_X}(T^*_{X, (\pm)^j}, \wedge^jT^*_{X, (\pm)^j}))$ for $j = 2, \ldots, \mathrm{rank}~T^*_{X, -}$. This $j$ is referred to as the \emph{splitting type} of the atlas $(\Ufr, \rho)$. 
\begin{center}
\begin{minipage}[t]{0.7\linewidth}
\begin{DEF}\label{tvirri4f4f4f44orjpjrpov} 
\emph{Any atlas for a supermanifold with splitting type $j = 2$ will be referred to as \emph{primary}. Similarly, any obstruction to splitting a supermanifold atlas of splitting type $j = 2$ will be referred to as a \emph{primary obstruction}.}
\end{DEF}
\vspace{\dp0}
\end{minipage}
\end{center}
\subsection{Obstructed Thickenings}
Fix a model $(X, T^*_{X,-})$. The cohomology groups $H^1(X, \mathcal Hom_{\Oc_X}(T^*_{X, (\pm)^j}, \wedge^jT^*_{X, (\pm)^j}))$ are referred to as \emph{obstruction spaces}; and elements therein are referred to as \emph{obstruction classes}. Generally, if we are given some obstruction class $\om$, there need \emph{not} exist a supermanifold with an atlas whose obstruction is this given $\om$. Such obstruction classes $\om$ were referred to as \emph{obstructed thickenings} in \cite{BETTOBSTHICK} and were said to `fail to integrate to a supermanifold' in \cite{BETTAQ}. In analogy with deformations of complex manifolds, this failure for obstruction classes $\om$ to integrate is measured by cohomology groups in degree two, identified originally by Eastwood and Lebrun in \cite{EASTBRU} and studied further by the author in \cite{BETTOBSTHICK, BETTAQ}. On a Riemann surface however, all such cohomology groups will vanish. Now in \cite{BETTOBSTHICK}, thickenings were categorised into three types: (1) supermanifolds; (2) pseudo-supermanifolds; and (3) obstructed thickenings. A pseudo-supermanifold is a thickening which, in a sense, becomes obstructed `eventually'. As there cannot exist any obstructed thickenings over a Riemann surface $X$, it follows that thickenings over $X$ must coincide with supermanifolds over $X$. Hence we have:
\begin{center}
\begin{minipage}[t]{0.7\linewidth}
\begin{THM}\label{h8937893h0fj309j03} 
Let $(X, T^*_{X, -})$ be a model with $X$ a Riemann surface. Then, given any obstruction class $\om$ there will exist a supermanifold $\Xfr$ modelled on $(X, T^*_{X, -})$ with atlas $(\Ufr, \rho)$ such that $\om_{(\Ufr, \rho)} = \om$. \qed
\end{THM}
\vspace{\dp0}
\end{minipage}
\end{center}
{\bf Remark.} One could view Theorem \ref{h8937893h0fj309j03} as the analogue of a classical result in the deformation theory of Riemann surfaces: \emph{any class in the cohomology group $H^1(X, T_X)$, for $X$ a Riemann surface, will be the Kodaira-Spencer class of some deformation $\Xc\ra (\Cbb^m, {\bf 0})$ of $X$ and some $m$.}\footnote{e.g., $m = \dim H^1(X, T_X)$, in which case one gets versal deformations, being deformations whose Kodaira-Spencer map is an isomorphism.}

\section{Superconformal Maps and Vector Fields}

\subsection{Superconformal Structures on $\Cbb^{1|1}$}

\subsubsection{Affine and Euclidean Superspace}
Following Manin in \cite{YMAN}, over a field $k$ one defines affine superspace $\Abb^{m|n}_k$ as the spectrum $\mathrm{Spec}~k[x_1, \ldots, x_m|\q_1, \ldots, \q_n]$, where $x_i$ are even and in the center of the algebra; and $\q_i$ are odd and so anti-commute amongst each other. In analogy then with constructions in algebraic geometry, complex Euclidean space $\Cbb^{m|n}$ are the complex $(0|1)$-points of the corresponding affine space and can be written $\Cbb^{m|n} = \Hom(\mathrm{Spec}~\Cbb[\xi], \Abb^{m|n}_\Cbb)$. As a supermanifolds, 
\begin{align*}
\Abb^{m|n}_\Cbb
=
S(\Abb^m_\Cbb, \oplus^n\Oc_{\Abb^m_\Cbb})
&&
\mbox{and}
&&
\Cbb^{m|n} 
=
S(\Cbb^m, \oplus^n\Oc_{\Cbb^m})
\end{align*}
where $\Oc_{\Abb^m_\Cbb}$ and $\Oc_{\Cbb^m}$ are the structures sheaves on $\Abb^m_\Cbb$ and $\Cbb^m$ respectively. Complex supermanifolds of dimension $(m|n)$ are locally isomorphic to $\Cbb^{m|n}$; while superspace schemes over a field $k$ are locally isomorphic to $\Abb^{m|n}_k$.

\subsubsection{Superconformal Generator}
In dimension $(1|1)$ we can form the derivation
\begin{align}
D_{(x|\q)}
\stackrel{\Delta}{=}
\frac{\pt}{\pt \q} + \q\frac{\pt}{\pt x}. 
\label{rhf748gf7hf98f3hf03j9}
\end{align}
The derivations $\Cbb[x|\q]\ra\Cbb[x|\q]$ inherit the $\Zbb_2$-grading on $\Cbb[x|\q]$ and, in accordance with this grading, the derivation $D_{(x|\q)}$ above is odd. It has the important property that $D_{(x|\q)}^2 = \pt/\pt x$. Hence $\{D_{(x|\q)}, D^2_{(x|\q)}\}$ forms a basis for the derivations on $\Abb^{1|1}_\Cbb$; and for the vector fields on $\Cbb^{1|1}$. The vector field $D_{(x|\q)}$ in \eqref{rhf748gf7hf98f3hf03j9} is referred to as a \emph{superconformal generator}. It defines a `superconformal structure' $\Cbb^{1|1}$.

\subsection{Superconformal Maps}
We consider now two coordinate systems $(x|\q)$ and $(y|\eta)$ on $\Cbb^{1|1}$. Denote by $\Cbb^{1|1}_{(x|\q)}$ resp. $\Cbb^{1|1}_{(y|\eta)}$ the superspace $\Cbb^{1|1}$ with a choice of coordinate system $(x|\q)$ resp. $(y|\eta)$. The general form of a map between superspaces $\Phi: \Cbb_{(x|\q)}^{1|1}\ra\Cbb_{(y|\eta)}^{1|1}$ parametrised by $\Abb_\Cbb^{0|n} =\mathrm{Spec}~\Cbb[\xi_1, \ldots, \xi_n]$ is:
\begin{align}
&&y &= \vp^+(x|\q, \xi)
\notag
\\
&& ~&=
f_0(x) + \sum_{|I|>0~\scriptsize{even}} \xi_I~f^I(x) + \sum_{|J|>0~\scriptsize{odd}}\q\xi_J~f^J(x)
\label{fjbckvhevibuebvoei}
\\
&&
~&\stackrel{\Delta}{=} f_+(x|\xi) + \q~f_-(x|\xi);
\label{rfh78gf7hf893f093jf39}
\\
\text{and}&&
\eta &=
\vp^-(x|\q,\xi)
\notag
\\
&& 
~&=
\q~\zeta(x)
+
\sum_{|J|>0~\scriptsize{odd}}\xi_J~\psi^J(x) + 
\sum_{|I|>0~\scriptsize{even}} \q\xi_I~\zeta^I(x)
\label{fjbkdbvjkbjkvbkeeee}
\\
&&
~&\stackrel{\Delta}{=}
\q~ \zeta(x|\xi) + \psi(x|\xi).
\label{fnibvbenvoieep}
\end{align}
Evidently $\Phi$ is a certain map $\Abb_\Cbb^{0|n}\times \Cbb^{1|1}_{(x|\q)}\ra \Cbb^{1|1}_{(y|\eta)}$. 
\begin{center}
\begin{minipage}[t]{0.7\linewidth}
\begin{DEF}\label{ruhf94hrfjfijepwww} 
\emph{The map $\Phi: \Cbb^{1|1}_{(x|\q)}\ra\Cbb^{1|1}_{(y|\eta)}$ over $\Abb^{0|n}_\Cbb$ (i.e., parametrised by $\Abb_\Cbb^{0|n}$) is \emph{superconformal} if and only if it preserves the superconformal structures on $\Cbb^{1|1}_{(x|\q)}$ and $\Cbb^{1|1}_{(y|\eta)}$ respectively. That is, if and only if there exists some nowhere vanishing function $h$ on $\Cbb^{1|1}_{(y|\eta)}$ such that $\Phi_*D_{(x|\q)} = (\Phi^*h)D_{(y|\eta)}$.}
\end{DEF}
\vspace{\dp0}
\end{minipage}
\end{center} 
To emphasise that $\Phi$ is a map over $\Abb^{0|n}_\Cbb$ we will write $\Phi_{/\Abb_\Cbb^{0|n}}$. The following lemma which appears in \cite{RABCRANE} gives conditions under which $\Phi$ will be superconformal.
~\\
\begin{center}
\begin{minipage}[t]{0.7\linewidth}
\begin{LEM}\label{rfnoifnoi3fip3mpm3} 
Let $\Phi_{/\Abb^{0|n}_\Cbb}: \Cbb_{(x|\q)}^{1|1}\ra\Cbb_{(y|\eta)}^{1|1}$ be a map given by \eqref{rfh78gf7hf893f093jf39} and \eqref{fnibvbenvoieep}. Then $\Phi$ will be superconformal if and only if the following equations are satisfied:
\begin{align*}
\zeta^2
&=
\frac{\pt f_+}{\pt x}
+
\psi\frac{\pt \psi}{\pt x}
\\
f_-
&=
\zeta\psi
\end{align*}
\qed
\end{LEM}
\vspace{\dp0}
\end{minipage}
\end{center}
With $\Phi_{/\Abb^{0|n}_\Cbb}$ superconformal, the factor $h$ in Definition \ref{ruhf94hrfjfijepwww} is given on the image of $\Phi$ by:
\begin{align}
\Phi^*h
=
h\circ \Phi(x|\q, \xi)
=
\zeta(x) + \q \frac{\pt \psi}{\pt x}.
\label{fjbvkbvuirbvuoireiiejfpoffeuh}
\end{align}
We can therefore see that $h$ will be invertible whenever $\zeta$ is invertible, as will be the case if $\Phi$ is itself invertible.

\begin{EX}\label{fjbkbvhjbrbvenvoe}
To second order (i.e., up to terms of order $\xi^3$ and higher), Lemma \ref{rfnoifnoi3fip3mpm3} asserts that the coefficients of a superconformal map $\Phi_{/\Abb^{0|n}_\Cbb}: \Cbb_{(x|\q)}^{1|1}\ra\Cbb_{(y|\eta)}^{1|1}$ (using \eqref{fjbckvhevibuebvoei} and \eqref{fjbkdbvjkbjkvbkeeee} now) and satisfy the following relations:
\begin{align}
\zeta(x)^2 
&=
\frac{\pt f_0}{\pt x};
\label{ghbbjksklee2}
\\
f^i(x) &= \zeta(x)\psi^i(x);
\label{rf4f9hf083h09f3j0}
\\
\frac{1}{2}\frac{\pt f^{ij}}{\pt x}
&=
\zeta(x)\zeta^{ij}(x) - \frac{1}{2}\left(\frac{\pt \psi^i}{\pt x}\psi^j(x) - \psi^i(x)\frac{\pt \psi^j}{\pt x}\right)
\notag
\\
&= 
\zeta(x)\zeta^{ij}(x)
-
\frac{1}{2} Wr.(\psi^i, \psi^j)
\label{rvuih984hf80309fj3fjpp}
\end{align}
where in the last equation above it suffices to consider $i<j$; and where $Wr.(\psi^i, \psi^j)$ in \eqref{rvuih984hf80309fj3fjpp} is the Wronskian.
\end{EX}

\subsection{Superconformal Vector Fields}
A superconformal vector field is a vector field $\nu$ which preserves the superconformal structure. And so, for $\Cbb^{1|1}_{(x|\q)}$ equipped with its superconformal structure $D_{(x|\q)}$, a vector field $w$ on $\Cbb^{1|1}_{(x|\q)}$ is superconformal if and only if $[w, D_{(x|\q)}] = v(x|\q)D_{(x|\q)}$ for some  function $v$. From this definition, the following result can be obtained by explicit calculation. It also appears in a number of articles in the literature, e.g., in \cite{FRIEDAN, WITTRS}.
\begin{center}
\begin{minipage}[t]{0.7\linewidth}
\begin{LEM}\label{rfnoifnoi3fip3mpeeem3} 
A basis for the superconformal vector fields on $\Cbb^{1|1}_{(x|\q)}$ over $\Cbb$ is given by
\begin{align*}
&&
W^-_{(x|\q)}
&= 
w^-(x)\left(\frac{\pt}{\pt \q} - \q\frac{\pt}{\pt x}\right);
\\
\text{and}&&
W^+_{(x|\q)}
&=
w^+(x)\frac{\pt}{\pt x} + \frac{1}{2}\frac{\pt w^+}{\pt x}\q~\frac{\pt}{\pt \q},
\end{align*}
where $w^+(x)$ and $w^-(x)$ are smooth or holomorphic functions on $\Cbb$ (i.e., independent of $\q$).
\qed
\end{LEM}
\vspace{\dp0}
\end{minipage}
\end{center}
In Definition \ref{ruhf94hrfjfijepwww} we formulated the notion of a superconformal map \emph{over} a parameter space $\Abb^{0|n}_\Cbb$. Recall that it is a certain map $\Abb^{0|n}_\Cbb\times \Cbb^{1|1}\ra\Cbb^{1|1}$. A similar notion can be formulated for superconformal vector fields. 
 \begin{center}
\begin{minipage}[t]{0.7\linewidth}
\begin{DEF}\label{34t577g879h8h} 
\emph{A vector field $W$ on the product $\Cbb_{(x|\q)}^{1|1}\times \Abb^{0|n}_\Cbb$ is \emph{superconformal} if and only if there exists some function $f(x|\q, \xi)$ such that
$[W, D_{(x|\q)}] = f(x|\q, \xi)~D_{(x|\q)}$.}
\end{DEF}
\vspace{\dp0}
\end{minipage}
\end{center}
Unlike with superconformal maps, we do not necessarily require $f$ in Definition \ref{34t577g879h8h} to be invertible, e.g., it could vanish modulo the ideal $(\xi_1, \ldots, \xi_n)$. As in Lemma \ref{rfnoifnoi3fip3mpeeem3} we have the following.
\begin{center}
\begin{minipage}[t]{0.7\linewidth}
\begin{LEM}\label{rfnoifnoi3feem3} 
Any superconformal vector field $W_{(x|\q, \xi)}$ on $\Cbb_{(x|\q)}^{1|1}\times \Abb^{0|n}_\Cbb$ can be written 
\[
W_{(x|\q, \xi)} = \sum_{|I|\geq0} \xi_I W_{(x|\q)}^I
\]
where $I$ ranges over multi-indices; and $W_{(x|\q)}^I$ are superconformal vector fields on $\Cbb^{1|1}_{(x|\q)}$.
\end{LEM}
\vspace{\dp0}
\end{minipage}
\end{center}
\emph{Proof of Lemma \ref{rfnoifnoi3feem3}.}
Since $D_{(x|\q)}(\xi_I) = 0$ for any multi-index $I$ it follows that any vector field of the form $\sum_{|I|\geq 0}\xi_IW^I$, with $W^I$ superconformal, will itself be superconformal. Conversely, suppose $W_{(x|\q, \xi)}$ is a vector field on $\Cbb^{1|1}_{(x|\q)}\times \Abb^{0|n}_\Cbb$. A basis for the tangent space is given by $\{\pt/\pt x, \pt/\pt \q, \pt/\pt\xi_i\}_{i= 1, \ldots n}$ and so we can write:
\[
W_{(x|\q, \xi)} 
= 
\sum_{|I|>0}
\xi_I \left(W^I_{(x|\q)} + V^{I,i}(x|\q)\frac{\pt}{\pt\xi_i}\right).
\]
If $W_{(x|\q, \xi)}$ is superconformal, then linearity of the Lie bracket along with $D_{(x|\q)}(\xi_I) = 0$ shows that $W^I_{(x|\q)}$ must be superconformal. regarding the vector field $V^{I, i}\pt/\pt \xi_i$ we have, up to a sign factor:
\begin{align}
\left[V^{I,i}(x|\q)\frac{\pt}{\pt\xi_i}, D_{(x|\q)}\right] = D_{(x|\q)}(V^{I,i}(x|\q))\frac{\pt}{\pt \xi_i}
\label{rg64874h89gh483j903}
&&
\mbox{(since $[\pt/\pt\xi_i, D_{(x|\q)}] =0$).}
\end{align}
As we are assuming $W_{(x|\q, \xi)}$ is superconformal, the term \eqref{rg64874h89gh483j903} must vanish so therefore $D_{(x|\q)}(V^{I, i}) = 0$. This means:
\[
\frac{\pt V^{I, i}}{\pt \q} + \q\frac{\pt V^{I, i}}{\pt x} = 0.
\]
The above equation can be satisfied if and only if $V^{I, i}= 0$ since $\{1, \q\}$ are independent.\footnote{To clarify further, note that $\pt V^{I, i}/\pt \q$ will not contain any terms proportional to $\q$ while the term $\q\pt V^{I, i}/\pt x$ is explicitly proportional to $\q$.} The lemma now follows.
\qed

\section{Super Riemann Surfaces}

\subsection{Definition}
Recall that on $\Cbb_{(x|\q)}^{1|1}$, the superconformal structure $D_{(x|\q)}$ defines a basis for $T\Cbb^{1|1}_{(x|\q)}$, given by $\{D_{(x|\q)}, D^2_{(x|\q)}\}$. On a $(1|1)$-dimensional supermanifold then, a globally-defined basis for its tangent bundle can then be identified with a superconformal structure. More precisely, following the presentation by Donagi and Witten in \cite{DW1}:
\begin{center}
\begin{minipage}[t]{0.7\linewidth}
\begin{DEF}\label{4847f98f03jf93} 
\emph{Let $\Xfr$ be a $(1|1)$-dimensional, complex supermanifold. A subsheaf $\Dc\subset T_\Xfr$ is said to be a \emph{superconformal structure} if there exists a covering $\Ufr$ of $\Xfr$ such that for any open set $\Ufr_\al\in \Ufr$ there exists coordinates $(x|\q)$ in which $\Dc(\Ufr_\al)$ is generated by $D_{(x|\q)}$.} 
\end{DEF}
\vspace{\dp0}
\end{minipage}
\end{center}
\begin{center}
\begin{minipage}[t]{0.7\linewidth}
\begin{DEF}\label{tvirri4f4f44orjpjrpov} 
\emph{A \emph{super Riemann surface} is a $(1|1)$-dimensional supermanifold $\Xfr$ over a Riemann surface equipped with a choice of superconformal structure $\Dc$.}
\end{DEF}
\vspace{\dp0}
\end{minipage}
\end{center}

\subsection{The Sheaf of Superconformal Vector Fields}
We present now some further remarks made and results obtained by Donagi and Witten in \cite{DW1}. Let $\Scl = (\Xfr, \Dc)$ be a super Riemann surface. Since $\Dc\subset T_\Xfr$ is a subsheaf, we have a short exact sequence:
\begin{align}
0 \lra \Dc\lra T_\Xfr \lra T_\Xfr/\Dc\lra0.
\label{fcbdbjkdbhvjfdhve}
\end{align}
In Lemma \ref{rfnoifnoi3fip3mpeeem3} we have a basis for the superconformal vector fields on $\Cbb^{1|1}_{(x|\q)}$. Since $\Scl$ is covered by superconformal domains $(\Cbb^{1|1}_{(x|\q)})$, it makes sense to form a sheaf of superconformal vector fields on $C$. Now observe that we can write, modulo the superconformal generator:
\begin{align}
W^-_{(x|\q)} 
&=
w^-(x)\left(\frac{\pt}{\pt \q} - \q\frac{\pt}{\pt x}\right)
\notag
\\
&=
2w^-(x)\frac{\pt}{\pt \q} - w^-(x)D_{(x|\q)} 
\notag
\\
&\equiv 2w^-(x)\frac{\pt}{\pt \q}\mod D_{(x|\q)};
\label{rfh7rttrgf78hf893h0}
\\
W^+_{(x|\q)}
&=
w^+(x)\frac{\pt}{\pt x} + \frac{1}{2}\frac{\pt w^+}{\pt x}\q~\frac{\pt}{\pt \q}
\notag
\\
&=
w^+(x)\frac{\pt}{\pt x} + \frac{1}{2}\frac{\pt w^+}{\pt x}\q D_{(x|\q)}&&\mbox{(since $\q^2 = 0$)}
\notag
\\
&\equiv 
w^+(x)\frac{\pt}{\pt x} \mod D_{(x|\q)}.
\label{rfh7rttrgf78hf893h}
\end{align}
Evidently we can identify the quotient $T_\Xfr/\Dc$ with the sheaf of superconformal vector fields on $\Scl$ (as sheaves of modules over $\Cbb$).\footnote{Indeed, note that for $W^+$ an even, superconformal vector field, $fW^+$ will be even if and only if $f$ is constant---in contrast to the odd case.} From the equivalences in \eqref{rfh7rttrgf78hf893h0} and \eqref{rfh7rttrgf78hf893h}, we can readily deduce:
\begin{center}
\begin{minipage}[t]{0.7\linewidth}
\begin{LEM}\label{445f4g4g4g4} 
For any super Riemann surface $(\Xfr, \Dc)$,
\[
T_\Xfr/\Dc \cong \Dc\otimes_{\Cbb} \Dc.
\]
\end{LEM}
\vspace{\dp0}
\end{minipage}
\end{center}
\emph{Proof of Lemma \ref{445f4g4g4g4}.} Since $\Xfr$ is a supermanifold, its tangent sheaf will be globally $\Zbb_2$-graded and we can write $T_\Xfr = T_{\Xfr, \mathrm{even}}\oplus T_{\Xfr, \mathrm{odd}}$. Locally, in coordinates $(x|\q)$, the even and odd components of $T_\Xfr$ are generated by $\{\pt/\pt x, \pt/\pt \q\}$. Using \eqref{rfh7rttrgf78hf893h0} and \eqref{rfh7rttrgf78hf893h} it is clear then that the map 
$T_\Xfr/\Dc \ra (T_{\Xfr, \mathrm{even}}/\Dc)\otimes (T_{\Xfr, \mathrm{odd}}/\Dc)$ will give the desired isomorphism of $\Cbb$-modules.\qed

\subsection{The Structure Theorem}
An equivalent description of a super Riemann surface is then as a $(1|1)$-dimensional supermanifold $\Xfr$ together with an extension of its tangent sheaf as in \eqref{fcbdbjkdbhvjfdhve} and satisfying the relation in Lemma \ref{445f4g4g4g4}. Alternatively, from Definition \ref{4847f98f03jf93} we could also consider a super Riemann surface as a supermanifold admitting a `superconformal atlas'. This is a more classical viewpoint, appearing in early works on super Riemann surface theory such as \cite{FRIEDAN, RABCRANE}; and adopted more  recently by the author in \cite{BETTPHD, BETTSRS}. A superconformal atlas is an atlas $(\Ufr, \rho)$ where the transition functions $\rho$ are superconformal isomorphisms. From this description one can arrive the following, classical result which we might view as a structure-type theorem for super Riemann surfaces:
\begin{center}
\begin{minipage}[t]{0.7\linewidth}
\begin{THM}\label{fjbvbviueovinep} 
There exists a one-to-one correspondence between super Riemann surfaces $(\Xfr, \Dc)$ and spin structures on the reduced space $C$ of $\Xfr$.\qed
\end{THM}
\vspace{\dp0}
\end{minipage}
\end{center}
Any supermanifold of dimension $(1|1)$ is split. In view of Theorem \ref{fjbvbviueovinep} then, any super Riemann surface $\Scl$ can be written $\Scl = S(C, T^*_{C, -})$, where $C$ is a Riemann surface and $T_{C, -}^*$ is a spin structure on $C$, i.e., a holomorphic line bundle (or, invertible sheaf) on $C$ equipped with an isomorphism $T_{C, -}^*\otimes T_{C, -}^*\cong \om_C$, where $\om_C$ is the canonical bundle.

\part{Algebraic and Analytic Deformations}
\label{rcjoehc894h984h48}

\section{Algebraic Deformations}

\subsection{Definition}
In the remarks surrounding Theorem \ref{fjbvbviueovinep}, a super Riemann surface is a $(1|1)$-dimensional supermanifold which admits a superconformal atlas. A \emph{deformation} of a super Riemann surface, from the algebraic viewpoint, can be described analogously to the case of complex manifolds,
i.e., by deforming the gluing data. We refer to \cite{BETTSRS} for further details on this point. Presently, we give the following definition. 
\begin{center}
\begin{minipage}[t]{0.7\linewidth}
\begin{DEF}\label{fjbvbrfrfrfrfrviueovinep} 
\emph{Let $\Scl$ be a super Riemann surface with superconformal atlas $(\Ufr, \rho) = ((\Ufr_\al), (\Ufr_{\al\be}), (\rho_{\al\be}))$. A \emph{deformation} of $\Scl$ over $\Abb^{0|n}_\Cbb = \mathrm{Spec}~\Cbb[\xi_1, \ldots, \xi_n]$, denoted $\Xc\ra\Abb^{0|n}_\Cbb$, is a complex supermanifold which admits an atlas $(\Uc, \vartheta)$ where, for each $\al$, there exists a smooth, superconformal isomorphism $\Ufr_\al \stackrel{\cong}{\ra}\Uc_\al$ over $\Abb^{0|n}_\Cbb$; and, with respect to these isomorphisms, for each $\al, \be$, the transition functions $\vartheta_{\al\be} : \Uc_{\al\be}\stackrel{\cong}{\ra}\Uc_{\be\al}$ can be identified with superconformal isomorphisms $\Ufr_{\al\be}\stackrel{\cong}{\ra}\Ufr_{\be\al}$ over $\Abb^{0|n}_\Cbb$.}
\end{DEF}
\vspace{\dp0}
\end{minipage}
\end{center}
\begin{center}
\begin{minipage}[t]{0.7\linewidth}
\begin{DEF}\label{rf4fg784hf98h0f3} 
\emph{An \emph{algebraic deformation} of a super Riemann surface is a choice of atlas for the deformation $\Xc\ra\Abb^{0|n}_\Cbb$ as described in Definition \ref{fjbvbrfrfrfrfrviueovinep}.}
\end{DEF}
\vspace{\dp0}
\end{minipage}
\end{center}

\subsection{Transition Data}
From Definition \ref{fjbvbrfrfrfrfrviueovinep} and \ref{rf4fg784hf98h0f3}, an algebraic deformation of a super Riemann surface over $\Abb^{0|n}_\Cbb$ will have transitions functions of the following kind: firstly set $\Abb^{0|n}_\Cbb = \mathrm{Spec}~\Cbb[\xi_1, \dots, \xi_n]$ and let $(x|\q)$ resp., $(y|\eta)$ be coordinates on $\Ufr_\al$ resp., $\Ufr_\be$. Then local coordinates on $\Uc_\al$ resp., $\Uc_\be$ are $(x_\xi|\q_\xi)$ resp., $(y_\xi|\eta_\xi)$ related as follows:
\begin{align}
y_\xi &= \vartheta^+_{\al\be}(x_\xi|\q_\xi) 
\label{rbvyrvyuev7998hf083h0}
\\
\notag
&=
f_{\al\be}(x) + \sum_i \q\xi_i~f^i_{\al\be}(x) +\sum_{i< j} \xi_i\xi_j~g^{ij}(x) +\cdots;
\\
\eta_\xi &= 
\vartheta^-_{\al\be}(x_\xi|\q_\xi) 
\label{fgjhbhbhbwuodn}
\\
\notag
&=
\q~\zeta_{\al\be}(x)
+
\sum_i \xi_i~\psi_{\al\be}(x) 
+
\sum_{i<j}\q\xi_i\xi_j~\zeta^{ij}(x) 
+\cdots
\end{align}
Since $\vartheta = (\vartheta_{\al\be})$ is superconformal, Example \ref{fjbkbvhjbrbvenvoe} shows that we must have the following relations on intersections:
\begin{align}
\frac{\pt f_{\al\be}}{\pt x} 
&=
\zeta_{\al\be}(x)^2;
\label{fjbvrbviubvonoioie1}
\\
\psi_{\al\be}^i(x) &= \zeta_{\al\be}(x)f^i_{\al\be}(x);
\label{fjbvrbviubvonoioie12}
\\
\frac{1}{2}\frac{\pt g^{ij}_{\al\be}}{\pt x}
&=
\zeta_{\al\be}(x)\zeta_{\al\be}^{ij}(x) - \frac{1}{2}Wr.\left(\psi^i_{\al\be}, \psi^j_{\al\be}\right).
\label{fjbvrbviubvonoioie13}
\end{align}
The above relations however are still insufficient to guarantee that the atlas is superconformal. We need to ensure $(\vartheta_{\al\be})_{\al,\be\in I}$ will satisfy the cocycle conditions, being:
\begin{align*}
\vartheta_{\al\be}\circ \vartheta_{\be\al} = {\bf 1}_{\be\al}
~(\mbox{$\forall$ intersections})
&&
\mbox{and}
&&
\vartheta_{\al\gam} = \vartheta_{\be\gam}\circ\vartheta_{\al\be}
~\mbox{($\forall$ triple intersections)}
\end{align*}
where ${\bf 1}_{\be\al}$ above refers to the identity map $\Uc_{\be\al}\ra \Uc_{\be\al}$. Calculations identifying conditions ensuring $\vartheta = (\vartheta_{\al\be})$ will define a $1$-cocycle  were undertaken in \cite{BETTSRS},\footnote{see e.g., p. 10, Corollary 3.9 in the arXiv version} where the following necessary condition was found:
\begin{center}
\begin{minipage}[t]{0.7\linewidth}
\begin{PROP}\label{fjbvbrviueovinep} 
Let $(\Uc, \vartheta)$ be a superconformal atlas for a deformation $\Xc\ra\Abb_\Cbb^{0|n}$ of a super Riemann surface, i.e., an algebraic deformation of $\Scl$ over $\Abb^{0|n}_\Cbb$; and suppose $\vartheta$ is given by \eqref{fgjhbhbhbwuodn} and \eqref{fgjhbhbhbwuodn}. Then on all intersections $U_\al\cap U_\be$ we have
\[
Wr.(\psi^i_{\al\be}, \psi^j_{\al\be}) = 0
\]
for all $i, j$.
\qed
\end{PROP}
\vspace{\dp0}
\end{minipage}
\end{center}

\subsection{Global Properties of Deformations}
The total space $\Xc$ of a deformation $\Xc\ra\Abb^{0|n}_\Cbb$ of a super Riemann surface will be a complex supermanifold in its own right. Therefore, it can be studied as such. In \cite{BETTPHD, BETTSRS}, some global properties of $\Xc$ were obtained which we present in the following.
\begin{center}
\begin{minipage}[t]{0.7\linewidth}
\begin{THM}\label{rfrfefefefee} 
Let $\Xc\ra \Abb^{0|n}_\Cbb$ be a deformation of a super Riemann surface $S(C, T^*_{C, -})$. Then as a supermanifold:
\begin{enumerate}[$\bt$]
	\item $\Xc$ has dimension $(1|n+1)$;
	\item the model for $\Xc$ is $(C, T^*_{\Xc, -})$ where $T^*_{\Xc, -}$ is an extension of holomorphic bundles on $C$,
	\begin{align}
	0 \lra \oplus^n\Oc_C\lra T^*_{\Xc, -} \lra T^*_{C, -}\lra0
	\label{rfh78gf87gf9h3fj3f3}
	\end{align}
	where $\Oc_C$ is the structure sheaf of $C$;
	\item when $n = 1$, the class of $T^*_{\Xc, -}$ as the extension of bundles in \eqref{rfh78gf87gf9h3fj3f3} coincides with the primary obstruction to splitting $\Xc$, up to some constant factor.
\end{enumerate}
\qed
\end{THM}
\vspace{\dp0}
\end{minipage}
\end{center}
For any deformation $\Xc\ra\Abb^{0|n}_\Cbb$ modelled on $(C, T^*_{\Xc, -})$, the class of $T^*_{\Xc, -}$ as an extension of holomorphic bundles in \eqref{rfh78gf87gf9h3fj3f3} is taken to be the \emph{odd Kodaira-Spencer class of $\Xc$}. An interesting corollary now of Proposition \ref{fjbvbrviueovinep} is the following.
\begin{center}
\begin{minipage}[t]{0.7\linewidth}
\begin{THM}\label{rfrfef4r4f455efefee} 
The $\Xc\ra\Abb^{0|n}_\Cbb$ be a deformation of a super Riemann surface with model $(C, T^*_{\Xc, -})$. Then there exists a class $\Theta\in \mathrm{Ext}^1(T^*_{C, -}, \Oc_C)$ and a constant $c$ such that 
\[
\Theta(T^*_{\Xc, -}) = c~\Theta
\]
where $\Theta(T^*_{\Xc, -})$ is the odd Kodaira-Spencer class of $\Xc$ (i.e., the class of $T^*_{\Xc, -}$ as the extension in \eqref{rfh78gf87gf9h3fj3f3} of holomorphic bundles).
\end{THM}
\vspace{\dp0}
\end{minipage}
\end{center}
\emph{Proof of Theorem \ref{rfrfef4r4f455efefee}.}
If a tuple of functions are linearly dependent, then their Wronskian will vanish. The converse need not generally hold. If the given tuple of functions are analytic however, Peano late in the nineteenth century observed that the converse statement will in fact hold. That is, a vanishing Wronskian will imply linear dependence. Subsequently, at the turn of the century, Bochner gave a more rigorous proof. For these references and a modern day proof, see \cite{WRONDEP}. For our purposes, recall that since $\Xc\ra\Abb^{0|n}_\Cbb$ is a deformation of a super Riemann surface, it will entertain a superconformal atlas $(\Uc, \vartheta)$ where $\vartheta$ is given by \eqref{fgjhbhbhbwuodn} and \eqref{fgjhbhbhbwuodn}. With respect to this atlas we obtain a cocycle representative for the odd Kodaira-Spencer class:
\begin{align}
\Theta(T^*_{\Xc, -})_{\al\be} = \sum_i \xi_i~\psi^i_{\al\be}.
\label{rg874gf73hf98j3f0j30}
\end{align}
For each $i$, $\psi^i = (\psi^i_{\al\be})$ will be a $T_{C, -}$-valued $1$-cocycle and so will define an element $\Theta^i \in\mathrm{Ext}^1(T^*_{C, -}, \Oc_C)$. Thus we already see from \eqref{rg874gf73hf98j3f0j30} that $\Theta(T^*_{\Xc, -}) = \sum_i \Theta^i$. Now for each $i, j$, the terms $\psi^i_{\al\be}$ and $\psi^j_{\al\be}$ will be holomorphic and so, analytic functions.\footnote{Generally, for each $i$, $\psi^i = (\psi^i_{\al\be})$ is an abelian, $1$-cocycle on $\Uc$ and valued in $T_{C, -}$. On each intersections $\psi^i_{\al\be}$ is a section of $T_{C, -}(U_{\al\be})$. By local triviality we have $T_{C, -}(U_{\al\be})\cong \Oc_C(U_{\al\be})$ and so $\psi^i_{\al\be}$ can itself be identified with an analytic function $U_{\al\be}\ra\Cbb$.} From Proposition \ref{fjbvbrviueovinep} we know that their Wronskian will vanish. As such by the afore-mentioned remarks on the Wronskian of analytic functions, we see that $\psi_{\al\be}^i$ and $\psi_{\al\be}^j$ will be linearly dependent. Hence, their respective cohomology classes $\Theta^i$ and $\Theta^j$ will be proportional over $\Cbb$. As this is holds for all $i, j$, the theorem now follows.
\qed

\subsection{Superconformal Models}
Recall (e.g., in (sub-)Section \ref{knebvou4bvoioi3nvo3}) that supermanifolds were described as being modelled on the data of a space $X$ and locally free sheaf $T^*_{X, -}$. Pairs $(X, T^*_{X, -})$ are referred to as models. In Theorem \ref{rfrfefefefee} and Theorem \ref{rfrfef4r4f455efefee} we obtain constraints on the models which could, in principle, support a supermanifold with superconformal atlas. The following definition labelling this class of models  will be useful.
\begin{center}
\begin{minipage}[t]{0.7\linewidth}
\begin{DEF}\label{rfrfef55w7w7w} 
\emph{Let $C$ be a Riemann surface and $T^*_-$ a locally free sheaf on $C$. The model $(C, T^*_-)$ is said to be \emph{superconformal} if: (1) $T^*_-$ can be described as an extension of holomorphic bundles as in \eqref{rfh78gf87gf9h3fj3f3} for some spin structure $T_{C, -}^*$ on $C$; and (2) the corresponding extension class $\Theta(T^*_-)$ can be written
\[
\Theta(T^*_-) = c~\Theta
\] 
for some constant $c$ and class $\Theta\in \mathrm{Ext}^1(T^*_{C, -}, \Oc_C)$.}
\end{DEF}
\vspace{\dp0}
\end{minipage}
\end{center}
If $\Xc\ra\Abb^{0|n}_\Cbb$ is a deformation of a super Riemann surface, its total space $\Xc$ will be a complex supermanifold modelled on $(C, T^*_{\Xc, -})$. This model is referred to as the \emph{total space model} of the deformation. By Definition \ref{rfrfef55w7w7w}, it is a superconformal model.

\subsection{The Odd Kodaira-Spencer Class and Splitting}
Regarding questions of splitting the total space of a deformation $\Xc\ra \Abb^{0|n}_\Cbb$, recall from
Theorem \ref{rfrfefefefee} that the odd Kodaira-Spencer class coincides with the obstruction to splitting $\Xc$ when $n = 1$. More generally, for $n>1$, the \emph{primary} obstruction\footnote{recall Definition \ref{tvirri4f4f4f44orjpjrpov}} to splitting $\Xc$ maps onto this odd Kodaira-Spencer class.\footnote{c.f., Lemma \ref{hf8h8fh30f390jf93jf39}.} Theorem \ref{34fiufguhfojp3} will then imply:
 \begin{center}
\begin{minipage}[t]{0.7\linewidth}
\begin{THM}\label{rfrfefefefee555w7w7w} 
Let $\Xc\ra\Abb^{0|n}_\Cbb$ be a deformation of a super Riemann surface with total space model $(C, T^*_{\Xc, -})$. Denote by $\Theta(T^*_{\Xc, -})$ the class of $T_{\Xc, -}^*$ as the extension of holomorphic bundles in \eqref{rfh78gf87gf9h3fj3f3}. If $\Theta(T^*_{\Xc, -})\neq0$, then $\Xc$ cannot be split, i.e., is non-split.
\qed
\end{THM}
\vspace{\dp0}
\end{minipage}
\end{center}

\subsection{Remarks on Higher Obstructions}
\label{rhg74hg89h03j09g408}
Concerning higher obstructions, in order for any to exist on a deformation $\Xc$ of a super Riemann surface, it is necessary for the primary obstruction to vanish and therefore $T_{\Xc, -}^*\cong \oplus^n\Oc_C\oplus T^*_{C, -}$ by Theorem \ref{rfrfefefefee555w7w7w}. It is conjectured in \cite{BETTSRS}, and confirmed to second order,\footnote{i.e., for deformations $\Xc\ra\Abb^{0|2}_\Cbb$} that deformations of super Riemann surfaces do \emph{not} admit higher obstructions. More explicitly:
\begin{center}
\begin{minipage}[t]{0.7\linewidth}
~\\
{\bf Conjecture.}
\emph{If the primary obstruction to splitting the total space of a deformation $\Xc\ra\Abb^{0|n}_\Cbb$ of a super Riemann surface vanishes, then $\Xc$ is split as a supermanifold.}
\\
\vspace{\dp0}
\end{minipage}
\end{center}
We refer to the above conjecture as the `No Exotic Deformations' conjecture. Further comments on it and a proposal of a proof sketch are given in the remarks concluding this article (see Section \ref{rhf47g794hg893093j9444}).

\section{Analytic Deformations}

\subsection{Preliminaries}
Following remarks by Witten in \cite{WITTRS} and Donagi-Witten in \cite{DW2}, only the even variables in superspace will be subject to conjugation. As such, in local coordinates $(x|\q)$ in complex superspace, a smooth function will be a function of the variables $(x, \widetilde x|\q)$.\footnote{Let $\Xfr$ be a supermanifold (smooth) with local coordinates $(x, \widetilde x|\q)$. Let $|\Xfr|\subset \Xfr$ be the reduced space (also smooth). This embedding defines a surjection of structure sheaves $i^\sharp: \Oc_\Xfr\ra\Oc_{|\Xfr|}\ra0$. Under this map we have: $i^\sharp (\widetilde x) = \overline{i^\sharp (x)}$, where $i^\sharp (x)$ is a local, complex coordinate on $|\Xfr|$.}
The Dolbeault operator, expressed in these local coordinates, is then the following, even, smooth vector $(0,1)$-form
\[
\widetilde\pt
=
d\widetilde x~\frac{\pt}{\pt\widetilde x}.
\]
In Proposition \ref{r84hg98h48g409g3} we see how certain perturbations of the Dolbeault operator can be identified with deformations of the complex structure. Similarly, certain perturbations of $\widetilde\pt$ can be realised as deformations of a super Riemann surface. This is the viewpoint adopted by Donagi and Witten in \cite{DW2} in their study of supermoduli; and by D'Hoker and Phong in their computation of superstring scattering amplitudes, reviewed for instance in \cite{HOKERPHONG1}.

\subsection{The Superconformal Vector $(0, 1)$-Form}
Our description of analytic deformations of a super Riemann surfaces involves smooth supermanifolds. We begin therefore with the following classical result in supermanifold theory due to Batchelor in \cite{BAT}.
 \begin{center}
\begin{minipage}[t]{0.7\linewidth}
\begin{THM}\label{rfrfefefef3333455555ee555w7w7w} 
Any smooth supermanifold is split.
\qed
\end{THM}
\vspace{\dp0}
\end{minipage}
\end{center}
Let $\Xc\ra \Abb^{0|n}_\Cbb$ be a deformation of a super Riemann surface $S(C, T^*_{C, -})$. Recall from Theorem \ref{rfrfefefefee} that its total space $\Xc$ is modelled on $(C, T_{\Xc, -}^*)$ where $T^*_{\Xc, -}$ is the extension of holomorphic bundles on $C$ in \eqref{rfh78gf87gf9h3fj3f3}. Hence, its split model is $S(C, T^*_{\Xc, -})$; and $\Xc$ and $S(C, T^*_{\Xc, -})$ are locally isomorphic. Denote by $\Xc^\8$ and $S(C, T^*_{\Xc, -})^\8$ the respective smooth supermanifolds underlying $\Xc$ and $S(C, T^*_{\Xc, -})$. Theorem \ref{rfrfefefef3333455555ee555w7w7w} asserts that $\Xc^\8$ and $S(C, T^*_{\Xc, -})^\8$ are \emph{diffeomorphic}.\footnote{\label{g84g784h9f8h344fdvvbkk}Explicitly, $S(C, T^*_{\Xc, -})^\8\cong S(C, T^*_{C, -})^\8\times \Abb^{0|n}_\Cbb$, where $C^\8$ is the smooth Riemann surface underlying $C$. As a locally ringed space, $S(C, T^*_{\Xc, -})^\8$ is a smooth supermanifold with reduced space $C^\8$ and structure sheaf $\Oc_{C^\8}\otimes \wedge^\bt_{\Oc_{C^\8}} (\oplus^n\Oc_C\oplus T_{C, -}^*)$.}
Generally, as in the case of complex manifolds, we will view complex supermanifolds $\Xfr$ as smooth supermanifolds $\Xfr^\8$ equipped with a complex structure $\widetilde \pt$. Accordingly, we identify $S(C,T^*_{C, -}) \equiv \big(S(C,T^*_{C, -})^\8, \widetilde\pt\big)$ and $\Xc \equiv \big(\Xc^\8, \widetilde \pt_\xi\big)$. In analogy now with Proposition \ref{r84hg98h48g409g3} we have:
\begin{center}
\begin{minipage}[t]{0.7\linewidth}
\begin{THM}\label{rfeiddddddofofepfpe} 
Let $\Xc\ra\Abb^{0|n}_\Cbb = \mathrm{Spec}~\Cbb[\xi_1, \ldots, \xi_n]$ be a deformation of a super Riemann surface $S(C, T^*_{C, -})$. Write $S(C,T^*_{C, -}) = \big(S(C,T^*_{C, -})^\8, \widetilde\pt\big)$ and $\Xc = (\Xc^\8, \widetilde \pt_\xi)$; and fix a diffeomorphism $\Xc^\8\cong S(C, T^*_{\Xc, -})^\8$. Then there exists a smooth, even superconformal vector $(0, 1)$-form $\vp(\xi)$ on $S(C, T^*_{\Xc, -})^\8$ such that:
\begin{enumerate}[$\bt$]
	\item $\vp(\xi)\equiv 0$ modulo the ideal $(\xi_i)_{i = 1, \ldots, n}$ and;
	\item
$
\widetilde\pt_\xi
=
\widetilde\pt - \vp(\xi).
$
\end{enumerate}
\end{THM}
\vspace{\dp0}
\end{minipage}
\end{center}
\emph{Proof of Theorem \ref{rfeiddddddofofepfpe}.}
In what follows it will be convenient to think of $\Xc$ as a collection of fibers $\{\Xc_\xi\}_{\xi\in \Abb^{0|n}_\Cbb}$. We proceed now with the following observation. Let $(x|\q)$ and $(x_\xi|\q_\xi)$ be local coordinates on $\Xc_0 = S(C, T^*_{C, -})$ and $\Xc_\xi$. Then Definition \ref{fjbvbrfrfrfrfrviueovinep} requires the existence of a smooth, superconformal isomorphism over $\Abb^{0|n}_\Cbb$ of the respective coordinate neighbourhoods sending $(x|\q)\mapsto (x_\xi|\q_\xi)$, i.e., that there exists some smooth, superconformal isomorphism $u$ such that: 
\begin{align}
x_\xi = u^+(x,\widetilde x|\q, \xi)
&&
\mbox{and}
&&
\q_\xi
=
u^-(x, \widetilde x|\q, \xi).
\label{rfg784gf973hf8h0j39j93333}
\end{align}
Now let $D_{(x|\q)}$ and $D_{(x_\xi|\q_\xi)}$ denote local generators for the respective superconformal structures on $\Xc_0$ and $\Xc_\xi$. Note that
\begin{align}
\left[\widetilde\pt, D_{(x|\q)}\right] = 0
&&
\mbox{and}
&&
\left[\widetilde\pt_\xi, D_{(x_\xi|\q_\xi)}\right] = 0.
\label{rfh79hf983hf080f9f3f3f3}
\end{align}
With $(x_\xi|\q_\xi)$ and $(x|\q)$ related by a superconformal transformation $u$ as in \eqref{rfg784gf973hf8h0j39j93333}, it must preserve the superconformal structure. That is, we have:
\begin{align}
u_* D_{(x|\q)} = hD_{(x_\xi|\q_\xi)}
\label{rfg783gf73h98fh38f04f4f4f4}
\end{align}
 for some non-vanishing function $h = h(x, \widetilde x|\q, \xi)$. With a diffeomorphism $\Xc^\8\cong S(C, T^*_{\Xc, -})^\8\cong S(C, T^*_{C, -})^\8\times\Abb^{0|n}_\Cbb$ (c.f., footnote \eqref{g84g784h9f8h344fdvvbkk}) we write:
\begin{align}
\widetilde\pt_\xi = \widetilde\pt - \vp(\xi)
\label{fkvjrbvjkrbv4ivonoivn}
\end{align}
where $\vp(\xi)$ is a smooth, vector $(0,1)$-form on $S(C, T^*_{\Xc, -})$ such that $\vp \equiv 0$ modulo the ideal $(\xi_1, \ldots, \xi_n)$. In what follows it will be convenient to write $D_{(x_\xi|\q_\xi)} = g(x, \widetilde x|\q, \xi)D_{(x|\q)}$. By \eqref{rfg783gf73h98fh38f04f4f4f4} we have $g = h^{-1}$. To see why $\vp(\xi)$ will be superconformal we compute, using the observation in \eqref{rfh79hf983hf080f9f3f3f3}:\footnote{\label{riurgu4hg894hf803hf093j0}We use the following rules of commutation between operators and functions: if $O, O^\p$ are operators and $g$ a function, then 
\begin{align*}
[O, g] = O(g)
&&
\mbox{and}
&&
[O, gO^\p]
=
O(g)O^\p + g[O, O^\p].
\end{align*}}
\begin{align*}
0
&=
\left[\widetilde\pt_\xi, D_{(x_\xi|\q_\xi)}\right] 
\\
&=
\left[
\widetilde\pt - \vp(\xi),
gD_{(x|\q)}
\right]
\\
&
=
\left[\widetilde\pt, gD_{(x|\q)}\right]
-
\left[
\vp(\xi), gD_{(x|\q)}\right]
\\
&=
\widetilde\pt g~D_{(x||q)}
-
\vp(\xi)g~D_{(x|\q)}
- g\left[\vp(\xi), D_{(x|\q)}
\right]
&&\mbox{(c.f., footnote \eqref{riurgu4hg894hf803hf093j0}).}
\end{align*}
In using that $g = h^{-1}$ is invertible we see from the above equation that:
\begin{align}
\left[\vp(\xi), D_{(x|\q)}
\right]
=
\left(\frac{\widetilde\pt g - \vp(\xi)g}{g}\right)D_{(x|\q)}.
\label{fg8gf798430f03f3f33}
\end{align}
Hence $\vp(\xi)$ is a smooth, \emph{superconformal} vector $(0, 1)$-form on $S(C, T^*_{\Xc, -})$. It is even since, by \eqref{fkvjrbvjkrbv4ivonoivn}, it can be expressed as a difference of even, $(0, 1)$-forms.
\qed

\subsection{Analytic Families and Deformations}
In Theorem \ref{rfeiddddddofofepfpe} we began with a deformation $\Xc\ra \Abb^{0|n}_\Cbb$ of a super Riemann surface and viewed its total space as a smooth supermanifold with a complex structure. This leads to the following, preliminary notion of a deformation which we refer to simply as an `analytic family'.
\begin{center}
\begin{minipage}[t]{0.7\linewidth}
\begin{DEF}\label{rf4fg784hf98h0f323445} 
\emph{An \emph{analytic family} of a super Riemann surfaces over $\Abb^{0|n}_\Cbb$ with central fiber
 $\big(S(C, T^*_{C, -})^\8, \widetilde \pt\big)$ is a choice of complex structure $\widetilde \pt_\xi$ on the product $S(C, T^*_{C, -})^\8\times \Abb^{0|n}_\Cbb$ such that, with $\Abb^{0|n}_\Cbb = \mathrm{Spec}~\Cbb[\xi_1, \ldots, \xi_n]$, the difference $\widetilde\pt_\xi - \widetilde\pt$ is an even, smooth, superconformal vector $(0, 1)$-form which vanishes modulo the ideal $(\xi_1, \ldots, \xi_n)$.}
\end{DEF}
\vspace{\dp0}
\end{minipage}
\end{center}
The definition of an analytic deformation we will propose here preempts the correspondence we wish to eventually establish between analytic and algebraic deformations. This requires appropriately accounting for Proposition \ref{fjbvbrviueovinep}. We do so as follows. 
Let $\pt$ be the holomorphic de Rham differential. In local coordinates $(x, \widetilde x)$ we have $\pt = dx~\pt/\pt x$. With this differential we can form the Wronskian of differential forms and so, with this in mind, we propose:
\begin{center}
\begin{minipage}[t]{0.7\linewidth}
\begin{DEF}\label{rf4fg70f324fh894hg9843445} 
\emph{An \emph{analytic deformation} of a super Riemann surface $\big(S(C, T^*_{C, -})^\8, \widetilde \pt\big)$ over $\Abb^{0|n}_\Cbb$ is an analytic family $(\Xc^\8, \widetilde \pt_\xi)$ with complex structure $\widetilde\pt_\xi = \widetilde\pt - \vp(\xi)$ and vector $(0,1)$-form satisfying:
\begin{align*}
Wr.\left(\left.\frac{\pt \vp}{\pt \xi_i}\right|_{\xi=0}\right.&, \left.\left.\frac{\pt \vp}{\pt \xi_j}\right|_{\xi=0}\right)
\\
\stackrel{\Delta}{=}&~
\pt \left(\left.\frac{\pt \vp}{\pt \xi_i}\right|_{\xi=0}\right)\otimes \left.\frac{\pt \vp}{\pt \xi_j}\right|_{\xi=0}
-
\left.\frac{\pt \vp}{\pt \xi_i}\right|_{\xi=0}\otimes\pt\left(\left.\frac{\pt \vp}{\pt \xi_j}\right|_{\xi=0}\right)
\\
=&~ 0
\end{align*}
for all $i, j$.}
\end{DEF}
\vspace{\dp0}
\end{minipage}
\end{center}
In Lemma \ref{rfnoifnoi3feem3} a basis for the superconformal vector fields on $\Cbb^{1|1}\times \Abb^{0|n}_\Cbb$ was described. Adapting this description to super Riemann surfaces, we have the following explicit description for analytic deformations in local coordinates $(x, \widetilde x|\q, \xi)$ on $S(C, T^*_{C, -})^\8\times\Abb^{0|n}_\Cbb$:\footnote{keeping in mind that $\vp(\xi) = \widetilde\pt - \widetilde\pt_\xi$ has even parity.}
\begin{align}
\widetilde\pt_\xi
=
\widetilde \pt 
+
\sum_i \xi_i~\chi^i\left(\frac{\pt}{\pt \q} - \q\frac{\pt}{\pt x}\right)
+
\sum_{i<j}
\xi_i\xi_j
\left(
h^{ij}\frac{\pt}{\pt x} + \frac{1}{2}\frac{\pt h^{ij}}{\pt x}\q~\frac{\pt}{\pt \q}\right)
+\cdots
\label{tg74g94hf8hf033}
\end{align}
where $\chi^i$ and $h^{ij}$ are smooth, $(0,1)$-forms; and the ellipses `$\cdots$' comprise terms proportional to $\xi^3$ and higher. The expression \eqref{tg74g94hf8hf033} representing families of super Riemann surfaces was used by Donagi and Witten in \cite{DW2}. By Definition \ref{rf4fg70f324fh894hg9843445} now, the family $(\Xc^\8, \widetilde\pt_\xi)$ giving the vector $(0, 1)$-form in \eqref{tg74g94hf8hf033} will be a deformation if $Wr.(\chi^i,\chi^j) \stackrel{\Delta}{=} \pt \chi^i\otimes\chi^j - \chi^i\otimes\pt\chi^j = 0$ for all $i,j$.

\section{Equivalences of Deformations: Analytic}

\subsection{Preliminary Definition}
An analytic deformation $\widetilde\pt_\xi$ is referred to in \cite{HOKERPHONG1} as a \emph{gauge slice} and in \cite{DW2} as a \emph{gauge field}; with transformations and equivalences referred to as gauge transformations or gauge equivalences. We present here, using the language introduced so far in this article (i.e., terminology from Definition \ref{rf4fg784hf98h0f323445} and \ref{rf4fg70f324fh894hg9843445}), the treatment of gauge transformations as can be found in \cite{DW2}. Let $\big(S(C, T^*_{C, -})^\8, \widetilde\pt\big)$ be a super Riemann surface and $\big(S(C, T^*_{C, -})^\8\times \Abb^{0|n}_\Cbb, \widetilde\pt_\xi\big)$ an analytic family. Fix a smooth, superconformal vector field $\nu$ on $S(C, T^*_{C, -})^\8\times \Abb^{0|n}_\Cbb$.\footnote{c.f., Definition \ref{34t577g879h8h}.}  A \emph{transformation} of $\pt_\xi$ is then obtained after conjugating by the formal group element defined by $\nu$. That is, on forming the exponential $e^\nu$, a transformation of $\widetilde\pt_\xi$ is identified with the action on the complex structure:
\begin{align}
\widetilde\pt_\xi \longmapsto \widetilde\pt_\xi^\p = e^{-\nu}\widetilde\pt_\xi e^\nu.
\label{rthg78gf974hf8030f9j34r4f4f4f}
\end{align}
This leads to:
\begin{center}
\begin{minipage}[t]{0.7\linewidth}
\begin{DEF}\label{rf4fg784hf98h0f32347f6d7giu45} 
\emph{Two analytic families $\widetilde\pt_\xi$ and $\widetilde\pt_\xi^\p$ over $\Abb^{0|n}_\Cbb$ of a super Riemann surface are said to be \emph{equivalent} if and only if there exists a smooth, superconformal vector field $\nu$ relating $\widetilde\pt_\xi$ and $\widetilde\pt_\xi^\p$ as in \eqref{rthg78gf974hf8030f9j34r4f4f4f}.}
\end{DEF}
\vspace{\dp0}
\end{minipage}
\end{center}

\subsection{Equivalences (Linearised)}
By Lemma \ref{rfnoifnoi3feem3} (and as in \eqref{tg74g94hf8hf033}), in local coordinates $(x, \widetilde x|\q, \xi)$ on $S(C, T^*_{C, -})^\8\times\Abb^{0|n}_\Cbb$, any \emph{even}, superconformal vector field $\nu$ can be written:
\begin{align}
\nu
=
\sum_i \xi_i~\nu^i(x,\widetilde x)\left(\frac{\pt}{\pt \q} - \q\frac{\pt}{\pt x}\right)
+
\sum_{i<j}
\xi_i\xi_j
\left(
\nu^{ij}(x, \widetilde x)\frac{\pt}{\pt x} + \frac{1}{2}\frac{\pt \nu^{ij}}{\pt x}\q~\frac{\pt}{\pt \q}\right)+\cdots
\label{fvnibviuh89rhf3j9jeoeio}
\end{align}
where $\nu^i$ and $\nu^{ij}$ are smooth functions; and where the ellipses `$\cdots$' denote terms proportional to $\xi^3$ and higher. Donagi and Witten in \cite{DW2} obtained conditions ensuring when two analytic deformations $\widetilde\pt_\xi$ and $\widetilde\pt_\xi^\p$ will be gauge equivalent \emph{to second order}, i.e., modulo the ideal $(\xi^3)$. Their result is as follows:
\begin{center}
\begin{minipage}[t]{0.7\linewidth}
\begin{THM}\label{fhvurh84j093j3j33333} 
Let $\widetilde\pt_\xi$ and $\widetilde\pt^\p_\xi$ be two analytic families of a super Riemann surface $(S(C, T^*_{C, -})^\8, \widetilde\pt)$ over $\Abb^{0|n}_\Cbb$. Suppose they are respectively given by the following expressions:
\begin{align*}
\widetilde\pt_\xi
&= 
\widetilde\pt
+
\sum_i \xi_i~\chi^i\left(\frac{\pt}{\pt \q} - \q\frac{\pt}{\pt x}\right)
+
\sum_{i<j}
\xi_i\xi_j
\left(
h^{ij}\frac{\pt}{\pt x} + \frac{1}{2}\frac{\pt h^{ij}}{\pt x}\q~\frac{\pt}{\pt \q}\right);
\\
\widetilde\pt_\xi^\p
&= 
\widetilde\pt
+
\sum_i \xi_i~\chi^{i\p}\left(\frac{\pt}{\pt \q} - \q\frac{\pt}{\pt x}\right)
+
\sum_{i<j}
\xi_i\xi_j
\left(
h^{ij\p}\frac{\pt}{\pt x} + \frac{1}{2}\frac{\pt h^{ij\p}}{\pt x}\q~\frac{\pt}{\pt \q}\right);
\end{align*}  
where the terms proportional to $\xi^3$ and higher are omitted. The analytic families $\widetilde\pt_\xi$ and $\widetilde\pt_\xi^\p$ are linearly equivalent if and only if there exists a global, smooth, superconformal vector field $\nu$ as in \eqref{fvnibviuh89rhf3j9jeoeio} such that:
\begin{align}
\chi^{i\p}
&=
\chi^i + \widetilde\pt \nu^i;
\label{Pfnvfbvjkvnvoeefe}
\\
h^{ij\p}
&=
h^{ij} + \widetilde\pt \nu^{ij}
+
\nu^i\chi^j - \nu^j\chi^i.
\label{fnvknvkjnivnoienvo}
\end{align}
\qed
\end{THM}
\vspace{\dp0}
\end{minipage}
\end{center}
The term `linear' in the statement of Theorem \ref{fhvurh84j093j3j33333} reflects the use of the term `linearised' by Donagi and Witten in \cite{DW2}. It emphasises that the equivalence in \eqref{fnvknvkjnivnoienvo} is taken modulo terms quadratic in $(\nu^i)_{i = 1, \ldots, n}$. For our purposes in this article, it will be useful to include these terms.

\subsection{Equivalences (Second Order)}
Our objective in this article is to establish a correspondence between equivalence classes of analytic deformations with those of algebraic deformations, to second order. To that extent it will be useful to rewrite the equivalence in \eqref{fnvknvkjnivnoienvo}. We begin with the following.
\begin{center}
\begin{minipage}[t]{0.7\linewidth}
\begin{LEM}\label{rfroifmiofmo4mf4} 
Let $\nu$ and $\chi$ be odd, superconformal vector fields on $\Cbb^{1|1}$. Then
\[
\nu\chi \equiv \frac{1}{2}\left[\nu, \chi\right]
\]
where the above equivalence is taken modulo the superconformal structure.
\end{LEM}
\vspace{\dp0}
\end{minipage}
\end{center}
\emph{Proof of Lemma \ref{rfroifmiofmo4mf4}.}
On $\Cbb^{1|1}_{(x|\q)}$ with superconformal generator $D_{(x|\q)} = \pt/\pt\q + \q\pt/\pt x$ we can write the following by Lemma \ref{rfnoifnoi3fip3mpeeem3},
\begin{align}
\nu = v(x,\widetilde x)\left(\frac{\pt}{\pt \q} - \q\frac{\pt}{\pt x}\right)  
&&
\mbox{and}
&&
\chi = 
X(x, \widetilde x)
\left(\frac{\pt}{\pt \q} - \q\frac{\pt}{\pt x}\right).
\label{fbvrg7g7398fh3h03}
\end{align}
Since $\nu$ and $\chi$ are both odd, their bracket is given by the anticommutator, i.e., $[\nu,\chi] = \nu\chi + \chi\nu$. Calculating this, we get
\begin{align*}
\frac{1}{2}[\nu,\chi]
= 
- vX\frac{\pt}{\pt x} - \frac{1}{2}\frac{\pt (vX)}{\pt x}\q\frac{\pt}{\pt\q}
\equiv 
-\nu\chi\frac{\pt}{\pt x}&&\mbox{(c.f., \eqref{rfh7rttrgf78hf893h}).}
\end{align*}
And similarly,
\[
\nu\chi
=
-vX\frac{\pt}{\pt x} - v\frac{\pt X}{\pt x}\q~\frac{\pt}{\pt\q}
=
-vX\frac{\pt}{\pt x} - v\frac{\pt X}{\pt x}\q~D_{(x|\q)}
\equiv -vX\frac{\pt}{\pt x}.
\]
The lemma now follows.
\qed
\\\\
As a consequence of Lemma \ref{rfroifmiofmo4mf4}, the linear equivalence in \eqref{fnvknvkjnivnoienvo} becomes:
\begin{align*}
h^{ij\p}
=
h^{ij}
+
\widetilde\pt \nu^{ij}
+ 
\frac{1}{2}
\left(
\left[\nu^i, \chi^j\right] - \left[\nu^j, \chi^i\right]
\right).
\end{align*}
Upon including the terms quadratic in $(\nu^i)_{i=1, \ldots, n}$, we arrive at what we refer to as an \emph{equivalence to second order}:\footnote{for the corresponding bosonic analogue, see \cite[p. 23]{DW2}}
\begin{align}
h^{ij^\p}
=
h^{ij}
+
\widetilde\pt \nu^{ij}
+ 
\frac{1}{2}
\Big(
\left[\nu^i, \chi^j\right] - \left[\nu^j, \chi^i\right]
\Big)
+
\frac{1}{4}
\left(
[\nu^i, \widetilde\pt \nu^j]
- 
[\nu^j, \widetilde\pt \nu^i]
\right).
\label{fbvrbvyurbvuir98h38h3}
\end{align}
It is the above equivalence which we will consider in forming our correspondence between analytic and algebraic deformations.
\\\\
{\bf Remark.} We mention now for completeness that both \eqref{Pfnvfbvjkvnvoeefe} and \eqref{fnvknvkjnivnoienvo} can be subsumed in the requirement that $\vp(\xi)$ be a vector $(0, 1)$-form on $\Xc^\8 = (S(C, T^*_{C, -})^\8\times\Abb^{0|n}_\Cbb, \widetilde\pt_\xi)$. That is, under equivalences of deformation parameters we have:
\begin{align}
\vp(\xi) \longmapsto \vp(\xi) + \widetilde\pt_\xi\nu.
\label{fjnfrbviutbiuvneioneuive}
\end{align}
As can be checked, both \eqref{Pfnvfbvjkvnvoeefe} and \eqref{fnvknvkjnivnoienvo} will follow from \eqref{fjnfrbviutbiuvneioneuive} above. To second order, we need to account for the quadratic terms in $\nu$ in \eqref{fbvrbvyurbvuir98h38h3}. In doing so, we modify \eqref{fjnfrbviutbiuvneioneuive} to get: $\vp(\xi) \mapsto \vp(\xi) + \widetilde\pt_\xi\nu + \frac{1}{2}\nu\widetilde\pt\nu$.

\section{Equivalences of Deformations: Algebraic}

\subsection{Definition}
Let $\Xc\ra \Abb^{0|n}_\Cbb$ be a deformation of $S(C, T^*_{C, -})$. By Definition \ref{fjbvbrfrfrfrfrviueovinep} it admits a superconformal atlas $(\Uc, \vartheta) = ((\Uc_\al), (\Uc_{\al\be}),(\vartheta_{\al\be}))$ and any such atlas was referred to in Definition \ref{rf4fg784hf98h0f3} as an \emph{algebraic deformation}. In Theorem \ref{fhvurh84j093j3j33333} we described the transformation undergone by analytic families under gauge field $u = e^\nu$, leading to the notion of equivalence. Presently, we will establish the transformations undergone by algebraic deformations. Recall from \eqref{rbvyrvyuev7998hf083h0} and \eqref{fgjhbhbhbwuodn} the general form of the transition functions $\vartheta_{\al\be}$ on intersections $\Uc_\al\cap \Uc_\be$. By superconformality of $\vartheta_{\al\be}$, the coefficient functions are subject to the relations in \eqref{fjbvrbviubvonoioie1}, \eqref{fjbvrbviubvonoioie12} and \eqref{fjbvrbviubvonoioie13}. A transformation of algebraic deformations consists of a collection of superconformal isomorphisms $\lam = (\lam_\al : \Uc_\al\stackrel{\cong}{\ra}\tilde \Uc_\al)$ which acts on the given atlas $(\Uc,\vartheta)$ by
\begin{align}
\Uc_\al \longmapsto\tilde\Uc_\al \stackrel{\Delta}{=} \lam_\al(\Uc_\al)
&&
\mbox{and}
&&
\vartheta_{\al\be} \longmapsto \tilde \vartheta_{\al\be} \stackrel{\Delta}{=} \lam^{-1}_\be \circ \vartheta_{\al\be}\circ \lam_\al.
\label{rhfthg8g9j40gj433}
\end{align}
Since compositions of superconformal maps will be superconformal, the new (transformed) atlas $(\tilde\Uc, \tilde\vartheta)$ in \eqref{rhfthg8g9j40gj433} will be superconformal and hence be an algebraic deformation. Now in Definition \ref{ruhf94hrfjfijepwww} we defined the notion of a superconformal map over $\Abb^{0|n}_\Cbb = \mathrm{Spec}~\Cbb[\xi_1, \ldots, \xi_n]$. A superconformal isomorphism is an invertible, superconformal map. A superconformal \emph{automorphism} over $\Abb^{0|n}_\Cbb$ is a superconformal isomorphism which coincides with the identity modulo the ideal $(\xi_1, \ldots,\xi_n)$.
\begin{center}
\begin{minipage}[t]{0.7\linewidth}
\begin{DEF}\label{787grgirh89hf8033} 
\emph{Two algebraic deformations $(\Uc, \vartheta)$ and $(\tilde\Uc, \tilde\vartheta)$ are said to be \emph{equivalent} if there exists some collection of superconformal \emph{automorphisms} $\lam = (\lam_\al)$ relating the deformations as in \eqref{rhfthg8g9j40gj433}.}
\end{DEF}
\vspace{\dp0}
\end{minipage}
\end{center}

\subsection{Transformation Laws}
To obtain conditions under which algebraic deformations will be equivalent, in analogy with Theorem 
\ref{fhvurh84j093j3j33333}, we need to see how the coefficient functions of $\vartheta_{\al\be}$ transform under the action of $\lam = (\lam_\al)$. To see this, let $(x_\xi|\q_\xi)$ be coordinates on $\Uc_\al$ with $(x|\q) = (x_\xi|\q_\xi)|_{\xi=0}$. We write,
\[
\lam_\al \equiv \lam_\al(x_\xi|\q_\xi) = (\lam_\al^+(x_\xi|\q_\xi) |\lam_\al^-(x_\xi|\q_\xi))
\]
where $\lam^+_\al$ and $\lam^-_\al$ are the even and odd components of $\lam$. The even and odd components of $\lam_\al$ to second order is written,
\begin{align}
\lam_\al^+(x_\xi|\q_\xi) &= u_\al(x) + \sum_i\xi_i\q~u_\al^i(x) + \sum_{i<j}\xi_i\xi_j~u^{ij}_\al(x);
\label{fnvkjbvkjbvnoieno1}
\\
\lam_\al^-(x_\xi|\q_\xi) &= \q~w_\al(x) + \sum_i\xi_i~w_\al^i(x) + \sum_{i<j} \xi_i\xi_j\q~w^{ij}_\al(x).
\label{fnvkjbvkjbvnoieno12}
\end{align}
Since $\lam_\al$ is superconformal the coefficient functions are subject to the following relations (c.f., Example \ref{fjbkbvhjbrbvenvoe}):
\begin{align}
\frac{\pt u_\al}{\pt x} &= w_\al(x)^2;
\label{fjvbfkbvuirbvuovnoe1}
\\
u^i_\al(x)
&=
w_\al(x)w^i_\al(x);
\label{fjvbfkbvuirbvuovnoe12}
\\
\frac{\pt u_\al^{ij}}{\pt x}
&= 
\frac{1}{2}w_\al(x)w^{ij}_\al(x) - \frac{1}{2} Wr.\left(w^i_\al, w^j_\al\right).
\label{fjvbfkbvuirbvuovnoe13}
\end{align}
Note that $\lam$ as described above is \emph{not necessarily} a superconformal automorphism. We will specialise to this case after the following calculation. 
Assuming that $(\Uc, \vartheta)$ and $(\tilde\Uc,\tilde \vartheta)$ are equivalent under $\lam$ we have by \eqref{rhfthg8g9j40gj433}:
\begin{align}
\lam_\be \circ \tilde\vartheta_{\al\be} = \vartheta_{\al\be}\circ \lam_\al.
\label{vhutuirhg4hgj4gjp3jgpojp}
\end{align}
Explicitly:
\begin{align}
\mathrm{RHS}_{\eqref{vhutuirhg4hgj4gjp3jgpojp}}
&=
\vartheta_{\al\be}\circ \left(\lam_\al^+(x_\xi|\q_\xi)|\lam^-_\al(x_\xi|\q_\xi)\right)
\notag
\\
&=
\left(\vartheta^+_{\al\be}\circ \left(\lam_\al^+(x_\xi|\q_\xi)|\lam^-_\al(x_\xi|\q_\xi)\right)
|
\vartheta^-_{\al\be}\circ \left(\lam_\al^+(x_\xi|\q_\xi)|\lam^-_\al(x_\xi|\q_\xi)\right)\right)
\label{fjvbdvbbvuiebvoinpie}
\end{align}
Note that it will be sufficient to identify the transformation laws for the even component of $\vartheta$ since those of the odd component can then be derived from the superconformal relations \eqref{fjbvrbviubvonoioie1}, \eqref{fjbvrbviubvonoioie12} and \eqref{fjbvrbviubvonoioie13} for $\vartheta$; and \eqref{fjvbfkbvuirbvuovnoe1}, \eqref{fjvbfkbvuirbvuovnoe12} and \eqref{fjvbfkbvuirbvuovnoe13} for $\lam$. Evaluating the even component in \eqref{fjvbdvbbvuiebvoinpie} now, we find:\footnote{The primed symbols indicate differentiation with respect to $x$, e.g., $f^\p = \pt f/\pt x$.}
\begin{align}
\vartheta^+_{\al\be}\circ &\big(\lam_\al^+(x_\xi|\q_\xi)|\lam^-_\al(x_\xi|\q_\xi)\big)
\notag\\
=&~
f_{\al\be}\circ \lam_\al^+(x_\xi|\q_\xi)
+
\sum_i\xi_i~\lam_\al^-(x_\xi|\q_\xi)~f_{\al\be}^i\circ \lam_\al^+(x_\xi|\q_\xi)
+
\sum_{i<j}\xi_i\xi_j~g^{ij}_{\al\be}\circ \lam_\al^+(x_\xi|\q_\xi)
\notag\\
=&~
f_{\al\be}(u_\al(x)) + f^\p_{\al\be}(u_\al(x)) \left(\sum_i\xi_i\q~u_\al^i(x) + \sum_{i<j} \xi_i\xi_j~u^{ij}_\al(x)\right)
\notag\\
&~+ 
\sum_i \xi_i\q~w_\al(x) f^i_{\al\be}(u_\al(x))
-
\sum_{i<j}\xi_i\xi_j \left( w_\al^i(x)f^j_{\al\be}(u_\al(x)) - w_\al^j(x)f^i_{\al\be}(u_\al(x))\right)
\notag\\
&~+
\sum_{i<j}\xi_i\xi_j ~g_{\al\be}^{12}(u_\al(x))
\notag\\
=&~
f_{\al\be}
\notag
\\
&~
+ \sum_i \xi_i\q~\left( f^\p_{\al\be} u^i_\al + w_\al f^i_{\al\be}\right)
+
\sum_{i<j}\xi_i\xi_j\left(f^\p_{\al\be}u_\al^{ij} -  w^i_\al f^j_{\al\be} + w^j_\al f^i_{\al\be} + g_{\al\be}^{ij}\right).
\notag
\end{align}
Similarly to \eqref{fjvbdvbbvuiebvoinpie}, the even component of the left-hand side is as follows:
\begin{align}
\mathrm{LHS}_{\eqref{vhutuirhg4hgj4gjp3jgpojp}}
=&~
\lam^+_\be\circ \left(\tilde\vartheta^{+}_{\al\be}(\tilde x_\xi|\tilde \q_\xi)|\tilde\vartheta_{\al\be}^{-}(\tilde x_\xi|\tilde\q_\xi)\right)
\notag\\
=
&~u_\be (\tilde f_{\al\be}(x)) + u_\be^\p(\tilde f_{\al\be}(x))\left(
\sum_i
\xi_i\q~\tilde f^i_{\al\be}  +\sum_{i<j} \xi_i\xi_j~\tilde g^{ij}_{\al\be}
\right)
\notag\\
&~+
\sum_i \xi_i\q~\tilde\zeta_{\al\be} u^i_\be +\sum_{i<j} \xi_i\xi_j\big(u^i_\be\tilde\psi^j_{\al\be} - u^j_\be \tilde\psi^i_{\al\be}\big)
+
\sum_{i<j}\xi_i\xi_j~ u^{ij}_\be(\tilde f_{\al\be}(x))
\notag\\
=&~
u_\be 
\notag
\\
\notag
&~
+ \sum_i\xi_i\q~\big(u^\p_\be \tilde f^i_{\al\be} + \tilde\zeta_{\al\be}u^i_\be \big)
+
\sum_{i<j}\xi_i\xi_j~\big( u^\p_\be\tilde g^{ij}_{\al\be} + u^i_\be\tilde\psi^j_{\al\be} - u^j_\be\tilde\psi^i_{\al\be} + u^{ij}_\be\big)
\end{align}
Equating $\mathrm{RHS}_{\eqref{vhutuirhg4hgj4gjp3jgpojp}}$ with $\mathrm{LHS}_{\eqref{vhutuirhg4hgj4gjp3jgpojp}}$ leads to the following relations between components of the transition functions $\vartheta_{\al\be}$ and $\tilde\vartheta_{\al\be}$:
\begin{align}
f_{\al\be}(u_\al(x)) &= u_\be(\tilde f_{\al\be}(x));
\label{rfg84gf87gf9h9f8h3}
\\
f^\p_{\al\be}(u_\al(x))u_\al^i(x) + w_\al(x)&f_{\al\be}^i(u_\al(x));
\notag
\\
&=
\notag 
\\
u^\p_\be(\tilde f_{\al\be}&(x))\tilde f^i_{\al\be}(x)
+ \tilde\zeta_{\al\be}(x)u^i_\be(\tilde f_{\al\be});
\notag
\\
f^\p_{\al\be}(u_\al(x)) u^{ij}_\al(x) - w_\al^i(x)f^j_{\al\be}(u_\al(x)) + w_\al^j(x)&f^i_{\al\be}(u_\al(x)) + g^{ij}_{\al\be}(u_\al(x))
\notag
\\
&=
\notag
\\
u_\be^\p(\tilde f_{\al\be}(x)) \tilde g^{ij}_{\al\be}(x) + u^i_\be(\tilde f_{\al\be}(x))& \tilde\psi^j_{\al\be}(x) - u^j_\be(\tilde f_{\al\be}(x))\tilde\psi^i_{\al\be} + u^{ij}_\be(\tilde f_{\al\be}(x)).
\notag
\end{align}

\subsection{Equivalence of Algebraic Deformations}

\subsubsection{Superconformal Automorphisms}
It was mentioned earlier that $\lam = (\lam_\al)$ need not necessarily be an automorphism. It will be a superconformal automorphism precisely when $u_\al(x) = x$ and $w_\al(x) = 1$. Evidently then, by \eqref{fjvbfkbvuirbvuovnoe12} and \eqref{fjvbfkbvuirbvuovnoe13}, for superconformal automorphisms we can equate: 
\begin{align}
u^i_\al(x) = w_\al^i(x)&&\mbox{and}&&
\frac{1}{2}\frac{\pt u^{ij}}{\pt x} &= w^{ij}(x) - \frac{1}{2}Wr.(u^i, u^j).
\label{rgygyeihuhf983h}
\end{align}
Note that when $\lam$ is a $0$-cochain of superconformal automorphisms, the complex structures on the underlying space, being the Riemann surface $C$ with holomorphic glueing data $(f_{\al\be})$ and $(\tilde f_{\al\be})$, must coincide by \eqref{rfg84gf87gf9h9f8h3}. Hence, the complex structure of the underlying space is fixed under an equivalence of algebraic deformations. 

\subsubsection{Equivalences}
We now have the following result detailing when two algebraic deformations will be equivalent in the sense of Definition \ref{787grgirh89hf8033}:
\begin{center}
\begin{minipage}[t]{0.7\linewidth}
\begin{THM}\label{rfhf9h3f8309fj3fjoooo3} 
Let $(\Uc, \vartheta)$ and $(\tilde\Uc, \tilde\vartheta)$ be algebraic deformations of a super Riemann surface $S(C, T^*_{C, -})$ and suppose there exists a $0$-cochain of superconformal automorphisms $\lam = (\lam_\al : \Uc_\al\stackrel{\cong}{\ra}\tilde\Uc_\al)$, given by \eqref{fnvkjbvkjbvnoieno1} and \eqref{fnvkjbvkjbvnoieno12} with $u_\al(x) = x$ and $w_\al(x) = 1$; and sending $(\Uc, \vartheta)\mapsto (\tilde\Uc, \tilde\vartheta)$ as in \eqref{rhfthg8g9j40gj433}. Then we have the following identities relating the coefficient functions of $\vartheta$ and $\tilde\vartheta$ on each intersection $\Uc_\al\cap\Uc_\be$ to second order:
\begin{align}
\psi_{\al\be}^i\frac{\pt}{\pt \eta} - \tilde\psi_{\al\be}^i\frac{\pt}{\pt \eta}
=&~
w^i_\be \frac{\pt}{\pt \eta} - w^i_\al\frac{\pt}{\pt\q};
\label{fgyur8g789h30h03}
\\
\label{fjkfbvrbvuinoiniwnviev}
g^{ij}_{\al\be}\frac{\pt}{\pt y}
-
\tilde g^{ij}_{\al\be}\frac{\pt}{\pt y}
=&~
u^{ij}_\be \frac{\pt}{\pt y}
-
u^{ij}_\al \frac{\pt}{\pt x}
\\
\notag
&
~+
\dt\left(\left\{w^i_\al\frac{\pt}{\pt \q}\right\}\right)_{\al\be}\otimes \psi^j_{\al\be}\frac{\pt}{\pt \eta}
\\
\notag
&
~-
\dt\left(\left\{w_\al^j\frac{\pt}{\pt \q}\right\}\right)_{\al\be}\otimes \psi^i_{\al\be}\frac{\pt}{\pt \eta}
\\
\notag
&~
+
\left(w_\al^iw_\be^j - w_\al^jw_\be^i\right)\frac{\pt}{\pt \q}\otimes \frac{\pt}{\pt\eta}
\end{align}
where, for $k = i, j$, the term $\{w^k_\al\pt/\pt\q\}$ above denotes a $0$-cochain of odd vector fields; and $\dt$ is the coboundary operator on \v Cech cocycles.
\end{THM}
\vspace{\dp0}
\end{minipage}
\end{center}
\emph{Proof of Theorem \ref{rfhf9h3f8309fj3fjoooo3}.}
Inspection of the relations succeeding \eqref{rfg84gf87gf9h9f8h3} gives:
\begin{align}
\q~f^i_{\al\be}\frac{\pt}{\pt y}
-
\q~\tilde f^i_{\al\be}\frac{\pt}{\pt y}
&=
\eta~u_\be^i\frac{\pt}{\pt y} - \q~u_\al^i\frac{\pt}{\pt x}
\label{rgf79g9fh8fh03hf93jf33}
\\
g^{ij}_{\al\be}\frac{\pt}{\pt y}
-
\tilde g^{ij}_{\al\be}\frac{\pt}{\pt y}
&=
u^{ij}_\be \frac{\pt}{\pt y}
-
u^{ij}_\al \frac{\pt}{\pt x}
\notag
\\
&~+
\left(u^i_\be\tilde \psi_{\al\be}^j - u^j_\be\tilde\psi^i_{\al\be}\right)\frac{\pt}{\pt y}
-
\left(
w^i_\al f^j_{\al\be} - w^j_\al f^i_{\al\be}
\right)
\frac{\pt}{\pt y}
\label{rhf849hff093jf0j93}
\end{align}
By superconformality (see \eqref{fjbvrbviubvonoioie12} and \eqref{rgygyeihuhf983h}), \eqref{rgf79g9fh8fh03hf93jf33} is equivalent to:
\begin{align}
\psi_{\al\be}^i\frac{\pt}{\pt \eta} - \tilde\psi_{\al\be}^i\frac{\pt}{\pt \eta}
=
w^i_\be \frac{\pt}{\pt \eta} - w^i_\al\frac{\pt}{\pt\q}
\label{rf84hf84hf0h03f93}
\end{align}
which gives \eqref{fgyur8g789h30h03}. Furthermore, again by superconformality, we can rewrite \eqref{rhf849hff093jf0j93} using only $w$ and $\psi$ giving
\begin{align}
g^{ij}_{\al\be}\frac{\pt}{\pt y}
-
\tilde g^{ij}_{\al\be}\frac{\pt}{\pt y}
&=
u^{ij}_\be \frac{\pt}{\pt y}
-
u^{ij}_\al \frac{\pt}{\pt x}
\notag
\\
&~+
\left(w^i_\be\tilde \psi_{\al\be}^j - w^j_\be\tilde\psi^i_{\al\be}\right)\frac{\pt}{\pt y}
-
\zeta_{\al\be}^{-1}\left(
w^i_\al \psi^j_{\al\be} - w^j_\al \psi^i_{\al\be}
\right)
\frac{\pt}{\pt y}.
\label{rgf478gf87hf893hf03jf93}
\end{align}
To eliminate the factor $\zeta_{\al\be}^{-1}$ we use the isomorphism $T_C^{1/2}\otimes T^{1/2}_C\stackrel{\cong}{\ra} T_C$, allowing for an identification of local sections $\pt/\pt \eta \otimes \pt/\pt \eta$ with $\pt/\pt y$ over $(y|\eta)$. Thus, with this identification, \eqref{rgf478gf87hf893hf03jf93} becomes the following cleaner expression:
\begin{align*}
g^{ij}_{\al\be}\frac{\pt}{\pt y}
-
\tilde g^{ij}_{\al\be}\frac{\pt}{\pt y}
&=
u^{ij}_\be \frac{\pt}{\pt y}
-
u^{ij}_\al \frac{\pt}{\pt x}
\\
&~+
\left(w^i_\be\tilde \psi_{\al\be}^j - w^j_\be\tilde\psi^i_{\al\be}\right)\frac{\pt}{\pt \eta}\otimes\frac{\pt}{\pt \eta}
-
\left(
w^i_\al \psi^j_{\al\be} - w^j_\al \psi^i_{\al\be}
\right)
\frac{\pt}{\pt \q}\otimes \frac{\pt}{\pt \eta}
\end{align*}
Now from \eqref{rf84hf84hf0h03f93} see that we can write $\tilde\psi$ in terms of $\psi$ and thereby further reduce the term in \eqref{rgf478gf87hf893hf03jf93}. And so, using that $\tilde\psi^i_{\al\be} = \psi^i_{\al\be} - w^i_\be + \zeta^{-1}_{\al\be}w^i_\al$, we find
\begin{align*}
\left(w^i_\be\tilde \psi_{\al\be}^j - w^j_\be\tilde\psi^i_{\al\be}\right)\frac{\pt}{\pt y}
=&~
\left(w^i_\be \psi_{\al\be}^j - w^j_\be\psi^i_{\al\be}\right)\frac{\pt}{\pt y}
\\
&~-
\left(w^i_\be w^j_\be - \zeta_{\al\be}^{-1}w^i_\be w^j_\al\right)\frac{\pt}{\pt y}
\\
&~+
\left(
w_\be^j w^i_\be - \zeta_{\al\be}^{-1}w_\be^j w_\al^i
\right)
\frac{\pt}{\pt y}
\\
=&~
\left(
w^i_\be \psi_{\al\be}^j - w^j_\be\psi^i_{\al\be}
+ \zeta^{-1}_{\al\be}\left(w_\be^iw_\al^j - w_\be^j w^i_\al\right)
\right)\frac{\pt}{\pt y}
\end{align*}
Hence
\begin{align*}
\eqref{rgf478gf87hf893hf03jf93}
&=
\left(\left(w_\be^i - \zeta^{-1}_{\al\be}w_\al^i\right)\psi^j_{\al\be}
+
\left(\zeta^{-1}_{\al\be} w_\al^j - w_\be^j\right)\psi_{\al\be}^i
+
 \zeta^{-1}_{\al\be}\left(w_\be^i w_\al^j - w_\be^j w^i_\al\right)\right)\frac{\pt}{\pt y}
 \\
 &=
\left( w_\be^i
 \frac{\pt}{\pt \eta}
 -
 w_\al^i\frac{\pt}{\pt \q}
 \right)
 \otimes \psi_{\al\be}^j\frac{\pt}{\pt \eta}
 +
 \left(
 w^j_\al\frac{\pt}{\pt \q}
 -
 w^j_\be \frac{\pt}{\pt \eta}
 \right)
 \otimes \psi_{\al\be}^i\frac{\pt}{\pt \eta}
 \\
 &~
 -
 \left(w^i_\al w^j_\be - w^j_\al w^i_\be\right)\frac{\pt}{\pt \q}\otimes\frac{\pt}{\pt \eta}.
\end{align*}
We thus obtain the claimed expression for \eqref{fjkfbvrbvuinoiniwnviev}. Theorem \ref{rfhf9h3f8309fj3fjoooo3} now follows.
\qed

\section{Main Result: Statement and Proof Outline}
\label{fknbviubovinpivompomp}

\subsection{Statement of Result}
In Theorem \ref{rfjrnvonoieep} an established correspondence was presented between the algebraic and analytic deformations of a complex manifold. We state now the central result of this article which one might view as an analogue of Theorem \ref{rfjrnvonoieep} for super Riemann surfaces.
\begin{center}
\begin{minipage}[t]{0.7\linewidth}
\begin{THM}\label{rfg748gf7f983h04i4i4i43hf03} 
For any super Riemann surface $S(C, T^*_{C, -})=(S(C, T^*_{C, -})^\8, \widetilde\pt)$ there exists a correspondence, to second order, between: (1) the algebraic deformations of $S(C, T^*_{C, -})$ over $\Abb^{0|n}_\Cbb$ up to equivalence; and (2) the analytic deformations of $(S(C, T^*_{C, -})^\8, \widetilde\pt)$ up to equivalence.
\end{THM}
\vspace{\dp0}
\end{minipage}
\end{center}
The correspondence in Theorem \ref{rfg748gf7f983h04i4i4i43hf03} can be stated more clearly as follows: let $\Xc\ra\Abb^{0|n}_\Cbb$ be a deformation of $S(C, T^*_{C, -})$. By Definition \ref{rf4fg784hf98h0f3}, algebraic deformations are superconformal atlases $(\Uc, \vartheta)$ on $\Xc$. In viewing $\Xc$ as a smooth supermanifold with complex structure $\big(S(C, T^*_{C, -})^\8\times \Abb^{0|n}_\Cbb, \widetilde\pt_\xi\big)$, analytic families of super Riemann surfaces over $\Abb^{0|n}_\Cbb$ with central fiber $S(C, T^*_{C, -})$ are identified with the complex structures $\widetilde\pt_\xi$ by Definition \ref{rf4fg784hf98h0f323445}; and are said to be deformations if $\widetilde\pt_\xi$ satisfies the condition in Definition \ref{rf4fg70f324fh894hg9843445}. Theorem \ref{rfg748gf7f983h04i4i4i43hf03} proposes then, \emph{to second order}: any atlas $(\Uc, \vartheta)$ defines an analytic deformation $\widetilde\pt_\xi$ of $S(C, T^*_{C, -})$ over $\Abb^{0|n}_\Cbb$; and conversely, any analytic deformation $\widetilde\pt_\xi$ defines a superconformal atlas $(\Uc, \vartheta)$. And furthermore, any equivalence of one kind of deformation translates into an equivalence of the other.

\subsection{Proof Outline of Theorem $\ref{rfg748gf7f983h04i4i4i43hf03}$}
Throughout, we fix a super Riemann surface $S(C, T^*_{C, -})\equiv(S(C, T^*_{C, -})^\8, \widetilde\pt)$. We break up the proof of Theorem \ref{rfg748gf7f983h04i4i4i43hf03} into two transitional steps: (1) algebraic $\Rightarrow$ analytic; and (2) analytic $\Rightarrow$ algebraic. The main difficulty lies in constructing one type of deformation, given the other. After these constructions are established, it remains to compare equivalences as described in Theorem \ref{fhvurh84j093j3j33333} (analytic) and Theorem \ref{rfhf9h3f8309fj3fjoooo3} (algebraic).

\subsubsection*{From Algebraic to Analytic}
We are given an algebraic deformation of the super Riemann surface $S(C, T^*_{C, -})$ over $\Abb^{0|n}_\Cbb$ which, recall, is a deformation $\Xc\ra\Abb^{0|n}$ of $S(C, T^*_{C, -})$ along with a superconformal atlas $(\Uc,\vartheta)$. From this data we obtain its primary obstruction $\om_{(\Uc, \vartheta)}$, which defines a class in the primary obstruction space of the model underlying $\Xc$. Upon suitably differentiating its image under the Dolbeault isomorphism with respect to the auxiliary, odd parameters $(\xi_i)$, we construct an analytic deformation $\big(S(C, T^*_{C, -})^\8\times \Abb^{0|n}_\Cbb, \widetilde\pt_\xi\big)$. The main issue lies in showing the difference $\widetilde \pt_\xi - \widetilde\pt$ will be a superconformal vector $(0, 1)$-form. 

\subsubsection*{From Analytic to Algebraic}
Starting with an analytic deformation of the super Riemann surface $(S(C, T^*_{C, -})^\8\times \Abb^{0|n}_\Cbb, \overline \pt_\xi)$, we construct an algebraic deformation. This involves firstly constructing a model $(C, T^*_-)$ on which to build a class of supermanifolds; and secondly, showing this model is superconformal.\footnote{c.f., Definition \ref{rfrfef55w7w7w}.} Subsequently, we construct a superconformal atlas, to second order, by firstly looking the preimage of the superconformal vector $(0, 1)$-form $\widetilde\pt_\xi - \widetilde\pt$ under the Dolbeault isomorphism and suitably differentiating with respect to the auxiliary, odd parameters $(\xi_i)$. That there will exist some supermanifold entertaining the given atlas will follows from our unobstructed thickenings result over Riemann surfaces in Theorem \ref{h8937893h0fj309j03}.

\part{Main Result: Proof}
\label{lfrpovviovuoh}

\section{From Algebraic to Analytic}

\subsection{The Primary Obstruction}
Let $\Xc\ra \Abb^{0|n}_\Cbb$ be a deformation of a super Riemann surface $S(C, T^*_{C, -})$; and let $(\Uc, \vartheta)$ be a superconformal atlas on $\Xc$. Recall that it will define a primary obstruction
 $\om_{(\Uc, \vartheta)} \in  H^1(C, \mathcal Hom(T^*_{\Xc, +}, \wedge^2T^*_{\Xc, -}))$ by Theorem \ref{34fiufguhfojp3}. Viewing $\Xc$ as a complex supermanifold, let $(C, T^*_{\Xc, -})$ be its model. In Theorem \ref{rfrfefefefee} we identified some general properties of $\Xc$, one of which being that $T_{\Xc, -}^*$ is an extension of the holomorphic line bundle  $T_{C, -}^*$ (see \eqref{rfh78gf87gf9h3fj3f3}). This leads to the following short exact sequence on $C$,
\begin{align}
0\lra
\oplus^{\binom{n}{2}}\Oc_C
\lra
\wedge^2T^*_{\Xc, -}\lra \oplus^n T^*_{C, -}\lra0
\label{rfg784gf73h8f09fj9333}
\end{align}
Since the reduced space of $\Xc$ is $C$ we have $T^*_{\Xc, +} = T^*_C$ (c.f., footnote \eqref{ftnoet3646f6r4}). Hence in applying $\mathcal Hom_{\Oc_C}(T^*_{\Xc, +},-) = \mathcal Hom_{\Oc_C}(T^*_{C},-)$ to \eqref{rfg784gf73h8f09fj9333} yields the following short exact sequence:
\begin{align}
0
\lra
\oplus^{\binom{n}{2}}T_C
\lra 
\mathcal Hom(T^*_C, \wedge^2T^*_{\Xc, -})
\lra
\oplus^n T_{C, -}
\lra
0
\label{rfh7f97hf98h340f093jf03}
\end{align}
where, for the right-most term (quotient factor) above, we used that $T_{C, -}^*$ is a spin structure, thereby giving an isomorphism $\mathcal Hom (T^*_C, T^*_{C, -})\cong T_{C, -}$.
Hence on cohomology we obtain a map 
\begin{align}
pr_*: H^1(C, \mathcal Hom(T^*_C, \wedge^2T^*_{\Xc,-})) \lra \oplus^nH^1(C, T_{C, -}),
\label{rfgyeuibiunoifworfeo}
\end{align}
which is necessarily surjective since $H^2(C, -) = (0)$ for dimensional reasons.
\begin{center}
\begin{minipage}[t]{0.7\linewidth}
 \vspace{0pt}
\begin{LEM}\label{hf8h8fh30f390jf93jf39} 
The map $pr_*$ in \eqref{rfgyeuibiunoifworfeo} sends the obstruction class $\om_{(\Uc,\vartheta)}$ of a superconformal atlas for $\Xc$ to the class $\Theta(T^*_{\Xc, -})$ of $T^*_{\Xc, -}$ as an extension of bundles in \eqref{rfh78gf87gf9h3fj3f3}.
\end{LEM}
\vspace{\dp0}
\end{minipage}
\end{center}
\emph{Proof of Lemma \ref{hf8h8fh30f390jf93jf39}}.
A proof of this lemma was given by the author in \cite{BETTSRS} in the case $n = 2$. We present the proof for general $n$ here as it will be illustrative for subsequent considerations. Let $(\Uc, \vartheta)$ be a superconformal atlas for $\Xc$, given by the data in \eqref{rbvyrvyuev7998hf083h0} and \eqref{fgjhbhbhbwuodn}; and subject to the relations in \eqref{fjbvrbviubvonoioie1}, \eqref{fjbvrbviubvonoioie12} and \eqref{fjbvrbviubvonoioie13}. Here, the obstruction class $\om_{\Uc,\vartheta)}$ is represented on intersections $\Uc_\al\cap\Uc_\be$ by the following, $\mathcal Hom(T^*_{\Xc,+}, \wedge^2T^*_{\Xc, -})$-valued $1$-cocycle,
\begin{align}
\om_{(\Uc, \vartheta), \al\be}
=
\left(\sum_i\xi_i\q~f^i_{\al\be} + \sum_{i<j} \xi_i\xi_j~g^{ij}_{\al\be}\right)\otimes \frac{\pt}{\pt y}.
\label{nfrfuih98fh3h30n3}
\end{align}
Similarly, the $1$-cocycle representing the class of $T_{\Xc, -}^*$ as the extension of holomorphic bundles on $C$ in \eqref{rfh78gf87gf9h3fj3f3} is capured in the odd component of the transition functions $\vartheta$. In denoting this class by $\Theta(T^*_{\Xc, -})$, on intersections it is represented by the $\oplus^nT_{C,-}$-valued $1$-cocycle
\[
\Theta(T^*_{\Xc, -})_{\al\be}
=
\sum_i \xi_i~\psi^i_{\al\be}\otimes\frac{\pt}{\pt \eta}.
\]
That the above expression defines a $1$-cocycle with respect to equivalences of algebraic deformations can be seen directly from the rule in \eqref{rf84hf84hf0h03f93}. The superconformal relation \eqref{fjbvrbviubvonoioie12} gives a canonical identification
\[
\sum_i \xi_i\left.\frac{\pt\om_{(\Uc,\vartheta), \al\be}}{\pt \xi_i}\right|_{\xi=0}\otimes \frac{\pt}{\pt y}
=
\Theta(T^*_{\Xc, -})_{\al\be}.
\]
The lemma now follows upon viewing the map $\sum_i\xi_i\frac{\pt}{\pt\xi_i}(-)|_{\xi=0}$ as the projection of cocycles representing the corresponding map $pr_*$ on cohomology in \eqref{rfgyeuibiunoifworfeo}.
\qed
\\\\
{\bf Remark.} In Theorem \ref{rfrfefefefee555w7w7w} and the surrounding discussion there, it was mentioned that the primary obstruction class maps onto the `odd Kodaira-Spencer class' of the deformation $\Xc$. 
This statement is precisely Lemma \ref{hf8h8fh30f390jf93jf39}. 

\subsubsection{Derivations from the Primary Obstruction}
We reiterate a key aspect of the proof of Lemma \ref{hf8h8fh30f390jf93jf39} here. For $\Xc\ra\Abb^{0|n}_\Cbb$ a deformation of a super Riemann surface with total space model $(C, T^*_{\Xc, -})$; and $(\Uc,\vartheta)$ a superconformal atlas for $\Xc$ with primary obstruction $\om_{(\Uc, \vartheta)}$, we found,
\begin{align}
pr_*\om_{(\Uc, \vartheta)}
=
\sum_i\xi_i\left.\frac{\pt\om_{(\Uc, \vartheta)}}{\pt \xi_i}\right|_{\xi=0}
=
\Theta(T^*_{\Xc, -}).
\label{rg84gf873hf98h80j03f}
\end{align}
Note that from the expression of $\om_{(\Uc, \vartheta)}$ on intersections in \eqref{nfrfuih98fh3h30n3} it makes sense to differentiate twice with respect to the parameters $\xi_i$ on the base $\Abb^{0|n}_\Cbb$. This gives a $T_C$-valued, $1$-cochain which, on intersections $\Uc_\al\cap\Uc_\be$, is:
\begin{align}
\left(
\sum_{i<j}
\xi_i\xi_j \left.\frac{\pt^2\om_{(\Uc,\vartheta)}}{\pt \xi_j\pt\xi_i}\right|_{\xi=0}
\right)_{\al\be}
=
\sum_{i<j}\xi_i\xi_j~g^{ij}_{\al\be}\otimes \frac{\pt}{\pt y}.
\label{rg84gf873hf98h80j03f2}
\end{align}
In contrast to \eqref{rg84gf873hf98h80j03f}, the above expression need not be define a cohomology class, as it need not satisfy a cocycle relation on triple intersections. The purpose of presenting \eqref{rg84gf873hf98h80j03f} and \eqref{rg84gf873hf98h80j03f2} is to illustrate how one can obtain derivative constructs from the algebraic deformation $(\Uc, \vartheta)$ directly from its primary obstruction.

\subsection{The Dolbeault isomorphism}
Let $\Xc\ra\Abb^{0|n}_\Cbb$ be a deformation of a super Riemann surface $S(C, T^*_{C, -})$ with total space model $(C, T^*_{\Xc, -})$. Let $C^\8$ denote the smooth Riemann surface underlying $C$. From the general statement of the vector bundle Dolbeault isomorphism in \eqref{rhf74g9h3f830fj03333333} we have:
\begin{align}
Dol:
H^1\big(C, \mathcal Hom(T^*_{\Xc, +}, \wedge^2T^*_{\Xc, -})\big)
\stackrel{\cong}{\lra}
H^{0,1}_{\overline\pt}\big(C^\8, \Hom(T^{*}_+\Xc^\8, \wedge^2T^{*}_-\Xc^\8)\big)
\label{egf783f98hf839fj903}
\end{align}
where $T^{*}_\pm\Xc^\8$ are smooth vector bundles on $C^\8$ with $T_{\Xc, \pm}^*$ their sheaf of holomorphic sections. Hence, associated to any primary obstruction $\om_{(\Uc, \vartheta)}$ to splitting a superconformal atlas $(\Uc, \vartheta)$ on $\Xc$ is a unique, smooth, $\overline\pt$-closed, $(0, 1)$ form valued in $\Hom(T^{*}_+\Xc^\8, \wedge^2T^{*}_-\Xc^\8)$, denoted $Dol(\om_{(\Uc,\vartheta)})$. 
Note that this is also true of the extension class $\Theta(T^*_{\Xc, -})$. Under \eqref{egf783f98hf839fj903} we have $Dol(\Theta(T^*_{\Xc, -}))\in H^{0, 1}_{\overline \pt}(C^\8, T_-C^\8)$, where $T_-C^\8$ is the smooth bundle over $C^\8$ whose sheaf of holomorphic sections coincide with $T_{C, -}$. From \eqref{rg84gf873hf98h80j03f} we see that:
\begin{align}
\sum_i\xi_i\left.\frac{\pt Dol(\om_{(\Uc, \vartheta)})}{\pt \xi_i}\right|_{\xi=0}
=
Dol(\Theta(T^*_{\Xc, -}))
\label{rcbyueuyg78g37gf893h93}
\end{align}
In subsequent sections we will find a similar relation to \eqref{rcbyueuyg78g37gf893h93} for the vector $(0, 1)$-form $\vp$ defining the complex structure $\overline \pt_\xi$ on $\Xc^\8$.

\subsection{An Analytic Deformation}
Let $\Xc\ra \Abb^{0|n}_{\Cbb}$ be a deformation of a super Riemann surface $S(C, T^*_{C, -})$; and let $(S(C, T^*_{C, -})^\8, \widetilde\pt)$ be its smooth model. Fix a supercomformal atlas $(\Uc, \vartheta)$ for $\Xc$ with primary obstruction $\om_{(\Uc, \vartheta)}$. On the product $S(C, T^*_{C, -})^\8\times \Abb^{0|n}_\Cbb$ set:
\begin{align}
\widetilde\pt_\xi
\stackrel{\Delta}{=}
\widetilde\pt 
-
\sum_i \xi_i\left.\frac{\pt Dol~\om_{(\Uc, \vartheta)}}{\pt \xi_i}\right|_{\xi=0}
-
\sum_{i<j} \xi_i\xi_j\left.\frac{\pt^2Dol~\om_{(\Uc, \vartheta)}}{\pt \xi_j\pt\xi_i}\right|_{\xi=0}.
\label{rgf64gf783hf983hf08j3093}
\end{align}
We intend to now prove the following.
\begin{center}
\begin{minipage}[t]{0.7\linewidth}
 \vspace{0pt}
\begin{PROP}\label{hg975hg984hg8003jg03} 
$\widetilde\pt_\xi$ in \eqref{rgf64gf783hf983hf08j3093} defines an analytic deformation of $S(C, T^*_{C, -})$ over $\Abb^{0|n}_\Cbb$.
\end{PROP}
\vspace{\dp0}
\end{minipage}
\end{center}
\emph{Proof of Proposition \ref{hg975hg984hg8003jg03}.}
We begin by showing that the difference $\widetilde\pt_\xi - \widetilde\pt$ in \eqref{rgf64gf783hf983hf08j3093} will define an even, superconformal vector $(0, 1)$-form on the product $S(C, T^*_{C, -})^\8\times\Abb_\Cbb^{0|n}$. This will give an analytic family (c.f., Definition \ref{rf4fg784hf98h0f323445}). Superconformality of the atlas $(\Uc, \vartheta)$ will ensure that this family will be a deformation in the sense of Definition \ref{rf4fg70f324fh894hg9843445}. To begin, recall from \eqref{rfh7rttrgf78hf893h0} and \eqref{rfh7rttrgf78hf893h} that the sheaf of superconformal vector fields on a super Riemann surface $S(C, T^*_{C, -}) = (\Xfr, \Dc)$ can be identified with $T_\Xfr/\Dc$. From \eqref{rfh7rttrgf78hf893h0} it is clear that:
\begin{align*}
\left\{
\begin{array}{l}
\mbox{odd, superconformal vector}
\\
\mbox{fields on $(\Xfr, \Dc) = S(C, T^*_{C, -})$}
\end{array}
\right\}
\cong
T_{\Xfr, \mathrm{odd}}/\Dc
\cong 
T_{C, -}.
\end{align*}
And similarly, any \emph{odd}, smooth, vector $(0, 1)$-form valued in $T_-C^\8$ will be superconformal. From this observation; \eqref{rcbyueuyg78g37gf893h93}; and the Dolbeault isomorphism $H^1(C, T_{C, -})\cong H_{\overline \pt}^{0, 1}(C^\8, T_-C^\8)$, it will follow that $\pt Dol~\om_{(\Uc, \vartheta)}/\pt \xi_i|_{\xi=0}$ is represented by a closed, smooth, $T_-C^\8$-valued $(0, 1)$ form on $C^\8$. Generally, $\pt Dol~\om_{(\Uc, \vartheta)}/\pt \xi_i$ will define an odd, superconformal vector $(0, 1)$-form on $S(C, T^*_{C, -})^\8\times\Abb_\Cbb^{0|n}$. The more subtle fact is superconformality of the second sum in \eqref{rgf64gf783hf983hf08j3093}. To confirm this we will use the following observation for vector fields:
\begin{align}
\begin{array}{l}
\mbox{{\bf Claim.} \emph{Suppose $O, O^\p$ are odd vector fields and $\pt$ is an}}
\\
\mbox{\emph{odd derivation with $\pt O^\p = 0$. Then $\pt[O, O^\p] = [\pt O, O^\p]$.}}
\end{array}
\label{jbfjvbkjnvvevj9}
\end{align}
\emph{Proof of Claim.}
By evaluating we find:
\begin{align}
\pt[O, O^\p] 
&= 
\pt(OO^\p + O^\p O)
\notag
\\
&=
(\pt O)O^\p - O(\pt O^\p) + (\pt O^\p)O - O^\p(\pt O)
\notag
\\
&=
[\pt O, O^\p].
\label{rgf745gf97hf893h0f34f3f}
\end{align}
Our conclusion \eqref{rgf745gf97hf893h0f34f3f} follows since $\pt O$ has even parity.\qed
\\\\
With \eqref{jbfjvbkjnvvevj9} now established, we can proceed thusly. In a coordinate neighbourhood $\Uc_\al$ on $\Xc$ with coordinates $(x_\xi|\q_\xi)$ let $D_{(x|\q)}$ be the generator for the superconformal structure on the central fiber, i.e., a local generator for the superconformal structure on $\Uc_\al\cap S(C, T^*_{C, -})$. As a vector field, it will have odd parity. Note that $\pt/\pt\xi_j$ is an odd derivation with 
\[
\frac{\pt}{\pt\xi_j}D_{(x|\q)} = 0.
\]
Hence,
\begin{align*}
\left[\frac{\pt^2}{\pt\xi_j\pt\xi_i}Dol~\om_{(\Uc, \vartheta)}, D_{(x|\q)}\right]
&=
\left[\frac{\pt}{\pt \xi_j}\left(\frac{\pt Dol~\om_{(\Uc, \vartheta)}}{\pt \xi_i}\right), D_{(x|\q)}\right]
\\
&=
\frac{\pt}{\pt\xi_j}\left[ \frac{\pt Dol~\om_{(\Uc, \vartheta)}}{\pt\xi_i}, D_{(x|\q)}\right]
~\mbox{(by \eqref{jbfjvbkjnvvevj9} since $\frac{\pt Dol~\om_{(\Uc, \vartheta)}}{\pt \xi_i}$ is odd)}
\\
&= \frac{\pt}{\pt\xi_j}\left(hD_{(x|\q)}\right)
~\mbox{(since $\pt Dol~\om_{(\Uc, \vartheta)}/\pt\xi_i$ is superconformal)}
\\
&=
\left(\frac{\pt h}{\pt\xi_j}\right)D_{(x|\q)}.
\end{align*}
And so, $\pt^2Dol~\om_{(\Uc, \vartheta)}/\pt\xi_j\pt\xi_i$ is even and superconformal. We can so far conclude that $(S(C, T^*_{C,-})^\8\times\Abb^{0|n}_\Cbb, \widetilde\pt_\xi)$ defines an analytic family. That it is an analytic deformation is a direct consequence of Proposition \ref{fjbvbrviueovinep}.
\qed 
\\\\
In Proposition \ref{hg975hg984hg8003jg03} we described how one can reconstruct an analytic deformation from an algebraic deformation. In the converse direction, reconstructing an algebraic from an analytic deformation is more involved. It requires the reconstruction of a superconformal model, a complex supermanifold; and a superconformal atlas. We begin therefore with some preliminary remarks and theory pertinent to the reconstruction to follow.

\section{Analytic to Algebraic: Preliminaries}
\label{vbrvuyg9v8h04vjssss}

\subsection{Notation}
Suppose we are given a complex structure 
\begin{align}
\widetilde
\pt_\xi
=
\widetilde\pt
-
\vp(\xi)
\label{rg874gf8798h80f30f3}
\end{align}
on the product $S(C, T^*_{C, -})^\8\times\Abb^{0|n}_\Cbb$. If $\vp(\xi)$ in \eqref{rg874gf8798h80f30f3} is an even, smooth, superconformal vector $(0, 1)$-form, then $(S(C, T^*_{C, -})^\8\times\Abb^{0|n}_\Cbb,\widetilde\pt_\xi)$ will define an analytic family of super Riemann surfaces over $\Abb^{0|n}_\Cbb$. We will not need recourse to this assumption at present however and so will consider $\widetilde\pt_\xi$ in full generality. We only assume $\widetilde\pt_\xi - \widetilde\pt \equiv 0$ modulo $(\xi)$.
Set $\Xc^\8 \stackrel{\Delta}{=} S(C, T^*_{C, -})^\8\times\Abb^{0|n}_\Cbb$. As a locally ringed space we can write 
\[
\Xc^\8 = \big(C^\8, \wedge^\bt (T^*_-C^\8\oplus \underline\Cbb^n)\big)
\]
where $C^\8$ is the smooth Riemann surface underlying $C$; $\underline\Cbb^n\ra C^\8$ is the trivial, complex bundle on $C^\8$ of rank $n$; and $T^*_-C^\8\oplus \underline\Cbb^n$ is their Whitney sum. Phrased alternatively, $\Xc^\8$ is a smooth supermanifold modelled on $(C, T^*_-\Xc^\8)$ where $T^*_-\Xc^\8 = T^*_-C^\8\oplus\underline\Cbb^n$. 

\subsection{The Tangent Bundle $T\Xc^\8$}
Vector fields on $\Xc^\8$ are sections of $T\Xc^\8$. With $T_+\Xc^\8 = TC^\8$ and $T^*_-\Xc^\8 = T^*_-C^\8\oplus \underline\Cbb^n$ the vector fields on $\Xc^\8$ can be described as sections over $C^\8$ as follows:
\begin{align*}
\Gam(C^\8, T\Xc^\8) &= \Gam\big(C^\8, (T_+\Xc^\8\oplus T_-\Xc^\8)\otimes \Oc(\Xc^\8)\big)
\\
&\cong
\oplus_{k\geq0}
\Hom_{C^\8}\big(T^*_+\Xc^\8\oplus T^*_-\Xc^\8, \wedge^kT^*_-\Xc^\8).
\end{align*}
The even and odd components are:
\begin{align}
T_{\mathrm{even}}\Xc^\8 
&\cong 
\oplus_k\Hom_{\Oc_{C^\8}}\big(T^*_{(\pm)^k}\Xc^\8, \wedge^kT^*_-\Xc^\8);
\label{fjbvkfbrbiuvbovnioenvpe}
\\
T_{\mathrm{odd}}\Xc^\8 
&\cong 
\oplus_k\Hom_{\Oc_{C^\8}}\big(T^*_{(\pm)^{k+1}}\Xc^\8, \wedge^kT^*_-\Xc^\8).
\label{fjvkbvurbvueniovneivmep}
\end{align}
These are smooth vector bundles over $C^\8$. The difference $\vp(\xi) = \widetilde\pt_\xi - \widetilde\pt$ is a $T_{\mathrm{even}}\Xc^\8$-valued $(0, 1)$-form on $S(C, T^*_{C, -})^\8$. With respect to the inclusion $C^\8\subset S(C, T^*_{C, -})^\8$ we can restrict to get $\vp(\xi)|_{C^\8}$, now a $T_{\mathrm{even}}\Xc^\8|_{C^\8}$-valued, $(0, 1)$-form on $C^\8$.

\subsection{The Holomorphic Tangent Bundle}
Recall that $\Xc^\8$ is the smooth supermanifold $S(C, T^*_{C, -})^\8\times\Abb^{0|n}_\Cbb$ with a choice of complex structure $\widetilde\pt_\xi$. Holomorphy can then be defined with respect to $\widetilde\pt_\xi$ as follows: a function $f\in \Oc(\Xc^\8)$ is holomorphic if and only if $\widetilde\pt_\xi f=0$. The sheaf of holomorphic functions will be denoted $\Oc_{hol.}(\Xc^\8)$. This description generalises to the tangent and cotangent vector fields, giving thereby holomorphic tangent and cotangent bundles. Denote these by
$T_{hol.}\Xc^\8$ and $T^*_{hol.}\Xc^\8$ respectively and let $\Oc(T_{hol.}\Xc^\8)$ resp. $\Oc(T^*_{hol.}\Xc^\8)$ denote their sheaves of holomorphic sections. 
\begin{center}
\begin{minipage}[t]{0.7\linewidth}
\vspace{0pt}
\begin{LEM}\label{jndrnuy78793f983irreerrr} 
The restriction $T_{hol.}\Xc^\8|_C$ defines a holomorphic vector bundle over $C$.
\end{LEM}
\vspace{\dp0}
\end{minipage}
\end{center}
\emph{Proof of Lemma \ref{jndrnuy78793f983irreerrr}}.
A section $X$ of $T\Xc^\8$ is holomorphic if and only if $\widetilde\pt _\xi X =0$. With $\widetilde\pt_\xi = \widetilde\pt - \vp(\xi)$ we see that holomorphic sections $X$ satisfy
\begin{align}
\widetilde\pt X = \vp(\xi) X. 
\label{rfg784gf79h9fhf80h309}
\end{align}
Locally, $X$ depends on parameters $(x, \widetilde x|\q, \xi)$ and expands in (wedge) powers of $\q$ and $\xi$. Hence if $X$ is holomorphic, the constraint in \eqref{rfg784gf79h9fhf80h309} necessarily requires $\widetilde\pt X|_{\xi = 0} = 0$. Moreover, the embedding $C^\8\subset S(C, T^*_{C, -})^\8$ preserves the respective complex structures and so is holomorphic. As the restriction of holomorphic vector fields on $S(C, T^*_{C, -})$ to $C$ will be holomorphic, it follows that $T_{hol.}\Xc^\8|_{C^\8}\ra C^\8$ will define a holomorphic vector bundle. 
\qed

\subsection{The Holomorphic Cotangent Bundle}
\label{9098g786fvuyviuvno}
Regarding the holomorphic cotangents $T^*_{hol.}\Xc^\8$, this is $\Oc_{hol.}(\Xc^\8)$-dual to the holomorphic tangents $T_{hol.}\Xc^\8$. As such, by Lemma \ref{jndrnuy78793f983irreerrr}, it will also define a holomorphic vector bundle over $C$. Let $\Xc^{\mathrm{an.}}$ denote the complex supermanifold associated to $(S(C, T^*_{C, -})^\8\times \Abb^{0|n}_\Cbb, \widetilde\pt_\xi)$. Recall the even and odd components of $T\Xc^\8$ modulo the fermionic ideal, $T_\pm\Xc^\8$. Their holomorphic sections also restrict to define holomorphic vector bundles over $C$. We set 
\[
T_{\Xc^{\mathrm{an.}}, \pm} \stackrel{\Delta}{=} \Oc_{hol.}(T_\pm\Xc^\8|_{C^\8}).
\]
Note that $T_{\Xc^{\mathrm{an.}}, +} = T_C$. The complex supermanifold $\Xc^{\mathrm{an.}}$ is modelled on $(C, T^*_{\Xc^{\mathrm{an.}}, -})$, where $T^*_{\Xc^{\mathrm{an.}}, -}$ is the $\Oc_C$-dual of $T_{\Xc^{\mathrm{an.}}, -}$. 
\begin{center}
\begin{minipage}[t]{0.7\linewidth}
\begin{LEM}\label{jndrnuirreerrr} 
There exists a short exact sequence of sheaves of holomorphic sections over $C$,
\begin{align}
0\lra \oplus^n\Oc_C \lra T^*_{\Xc^{\mathrm{an.}}, -}\lra T_{C, -}^* \lra 0.
\label{fjvbrvbiurvhoinvpoemve}
\end{align}
\end{LEM}
\vspace{\dp0}
\end{minipage}
\end{center}
\emph{Proof of Lemma \ref{jndrnuirreerrr}.}
In the proof of Lemma \ref{jndrnuy78793f983irreerrr} we showed $T_{hol.}\Xc^\8|_C$ is holomorphic. This required firstly showing $T_{hol.}\Xc^\8|_{S(C, T^*_{C,-})^\8}$ is a holomorphic vector bundles over $S(C, T^*_{C, -})^\8$ with respect to the given complex structure $\widetilde\pt$ (c.f., \eqref{rfg784gf79h9fhf80h309}). Hence, the embedding $S(C, T^*_{C, -})^\8\subset \Xc^\8 = S(C, T^*_{C, -})^\8\times\Abb^{0|n}_\Cbb$ is holomorphic and corresponds therefore to an embedding of complex supermanifolds $S(C, T^*_{C, -})\subset \Xc^{\mathrm{an.}}$. The normal bundle to this embedding will be an $\Abb_\Cbb^{0|n}$-bundle over $S(C, T^*_{C, -})$. In denoting its sheaf of sections by $\nu_{S(C, T^*_{C, -})/\Xc^{\mathrm{an.}}}$ we have the normal bundle sequence:
\begin{align}
0\lra T_{S(C, T^*_{C, -})} \lra \Oc(T_{hol.}\Xc^\8|_{S(C, T^*_{C, -})}) \lra \nu_{S(C, T^*_{C, -})/\Xc^{\mathrm{an.}}}\lra0.
\label{rfg654gf874hf8f09jf33333}
\end{align}
Now $S(C, T^*_{C, -})$ is a $(1|1)$-dimensional, complex supermanifold and is therefore split. Hence, as sheaves of $\Oc_C$-modules, $T_{S(C, T^*_{C, -})}\cong T_C\oplus T_{C, -}$. Similarly, observe that the restriction of $\Oc(T_{hol.}\Xc^\8|_{S(C, T^*_{C, -})})$ to $C$ splits into $T_{\Xc^{\mathrm{an.}}, +}\oplus T_{\Xc^{\mathrm{an.}}, -} = T_C\oplus T_{\Xc^{\mathrm{an.}}, -}$. Since the normal bundle is an $\Abb^{0|n}_\Cbb$-bundle over $S(C, T^*_{C, -})$, it has no even component; and its odd component restricted to $C$ will be trivial, i.e., $\nu_{S(C, T^*_{C, -})/\Xc^{\mathrm{an.}}}|_C\cong (0)\oplus (\oplus^n\Oc_C)$ as $\Oc_C$-modules. Hence, on restricting \eqref{rfg654gf874hf8f09jf33333} to $C$ and using that the maps in \eqref{rfg654gf874hf8f09jf33333} are even (i.e., preserving parity) we obtain: a triviality for the even component; and for the odd component we get the sequence,
\[
0 \lra T_{C, -}\lra T_{\Xc^{\mathrm{an.}}, -} \lra \oplus^n\Oc_C\lra0.
\]
The sequence in \eqref{fjvbrvbiurvhoinvpoemve} is evidently dual to the above sequence. \qed

\section{The Extension Class}

\noindent
From a choice of complex structure $\widetilde\pt_\xi$ on $\Xc^\8$ we obtain a model for a class of complex supermanifolds $(C, T^*_{\Xc^{\mathrm{an.}}, -})$ with $T^*_{\Xc^{\mathrm{an.}}, -}$ an extension of holomorphic sheaves on $C$ in Lemma \ref{jndrnuirreerrr}. Denote by $\Theta(T^*_{\Xc^{\mathrm{an.}}, -})$ its extension class. In the following series of investigations we aim to show how it relates to the vector $(0, 1)$-form $\vp(\xi)$ defining the complex structure $\widetilde\pt_\xi$.

\subsection{A Cocycle Representative for $\Theta(T^*_{\Xc^{\mathrm{an.}},-})$}
Starting from the sequence in \eqref{fjvbrvbiurvhoinvpoemve}, applying $\mathcal Hom_{\Oc_C}(T_{C, -}^*, -)$ and using that all the sheaves in consideration are locally free, we obtain another short exact sequence of sheaves
\begin{align}
0
\lra 
\oplus^nT_{C, -}
\lra
\mathcal Hom_{\Oc_C}(T^*_{C, -}, T^*_{\Xc^{\mathrm{an.}}, -})
\lra
\mathcal Hom_{\Oc_C}(T^*_{C, -}, T^*_{C, -})
\lra
0.
\label{fknvjkfvuirbvuonvoienvioenvoi}
\end{align}
The above sequence induces, on cohomology, the following piece:
\[
\xymatrix{
\cdots\ar[r] & \Hom(T_{C, -}^*, T^*_{C,-})\ar[r]^\dt & H^1\big(C,\oplus^nT_{C, -}\big)\ar[r]  & \cdots
}
\]
where $\Hom(T^*_{C, -}, T^*_{C, -}) = H^0(C, \mathcal Hom_{\Oc_C}(T^*_{C, -}, T^*_{C, -}))$. By definition,
\begin{align}
\Theta(T^*_{\Xc, -})\stackrel{\Delta}{=}\dt({\bf 1}_{T^*_{C, -}})
\label{fjvkfbvyurghvoijop}
\end{align}
where ${\bf 1}_{T^*_{C, -}}$ is the identity map $T^*_{C, -}=T^*_{C, -}$. A representative of $\dt({\bf 1}_{T^*_{C, -}})$ can be described as follows: on $\Xc^\8$ let $(x, \widetilde x|\q, \xi)$ and $(y, \widetilde y|\eta,\xi)$ denote overlapping coordinate systems. They (respectively) parametrise open neighbourhoods $\Uc^\8_\al, \Uc_\be^\8\subset \Xc^\8$. Viewing $\q$ and $\eta$ as local sections of $T^*_{C, -}$ we can write:
\begin{align*}
{\bf 1}_{T^*_{C, -}}|_{\Uc^\8_\al} = \q\frac{\pt}{\pt\q}
&&
\mbox{and}
&&
{\bf 1}_{T^*_{C, -}}|_{\Uc^\8_\be} = \eta\frac{\pt}{\pt\eta}.
\end{align*}
Since \eqref{fknvjkfvuirbvuonvoienvioenvoi} is exact, we can locally lift the above sections over $\Uc^\8_\al\cap C^\8$ and $\Uc^\8_\be\cap C^\8$ to expressions of the form $X_\al = f(x, \widetilde x |\q, \xi)\frac{\pt}{\pt\q}$ and $X_\be = g(y,\widetilde y|\eta, \xi)\frac{\pt}{\pt\eta}$, where $f$ and $g$ are local, holomorphic sections valued in $T^*_-\Xc^\8|_{C^\8}$.\footnote{This means we must write $f = f_0\q + \sum_i\xi_i~f^i_1$ and $g =  g_0\eta +\sum_i\xi_i~ g^i_1\xi_i$ for appropriate $f_0, f^i_1$ and $g_0, g^i_1$.} The surjection in  \eqref{fknvjkfvuirbvuonvoienvioenvoi} is induced by $\xi\mapsto 0$. Hence, in order for $X_\al \mapsto {\bf 1}_{T^*_{C, -}}|_{U_\al}$ and $X_\be\mapsto {\bf 1}_{T^*_{C, -}}|_{U_\be}$, we must have:
\begin{align}
X_\al = \left(\q + \sum_i \xi_i~f^i (x, \widetilde x)\right)\frac{\pt}{\pt\q} 
&&
\mbox{and}
&&
X_\be =  \left(\eta + \sum_i \xi_i~g^i (y, \widetilde y)\right)\frac{\pt}{\pt\eta}
\label{kfnkvbhjyugvi4hroijproee3}
\end{align}
The difference $X_{\al\be} \stackrel{\Delta}{=}X_\be-X_\al$ on intersections $\Uc^\8_\al\cap \Uc^\8_\be$ defines a cocycle representative for $\dt({\bf 1}_{T^*_{X, -}})$ and hence, by \eqref{fjvkfbvyurghvoijop}, of the extension class $\Theta(T^*_{\Xc, -})$.

\subsection{The Dolbeault Image of $\Theta(T^*_{\Xc^{\mathrm{an.}}, -})$}
We have so far (since Section \ref{vbrvuyg9v8h04vjssss}) not required the complex structure $\widetilde\pt_\xi$ on $\Xc^\8$ define an analytic family of super Riemann surfaces. That means $\vp(\xi)$ need not be superconformal. From here onwards however we can no longer dispense with this assumption. Hence, in what follows, $\Xc^\8$ will be an analytic family of super Riemann surfaces over $\Abb^{0|n}_\Cbb$ and so $\vp(\xi)$ will be even and superconformal.
Now by Lemma \ref{jndrnuirreerrr} we know that $\Theta(T^*_{\Xc^{\mathrm{an.}}, -})\in \mathrm{Ext}^1(T^*_{C, -}, \oplus^n\Oc_C)\cong H^1(C, \oplus^nT^*_{C, -})$. Recall from \eqref{rhf74g9h3f830fj03333333} that we presented the Dolbeault isomorphism in considerable generality. For our purposes here we specialise to the incarnation $Dol : H^1(C, \oplus^nT^*_{C, -})\stackrel{\cong}{\ra}H^{0, 1}_{\overline \pt}(C^\8, \oplus^nT^*_-C^\8)$ from whence it is clear that $Dol$ sends cocycle representatives of $\Theta(T^*_{\Xc^{\mathrm{an.}}, -})$ to $\overline \pt$-closed, $T_-C^\8$-valued $(0, 1)$-forms.
\begin{center}
\begin{minipage}[t]{0.7\linewidth}
\begin{LEM}\label{jndrnuirreerrr} 
Let $\Xc^\8$ be an analytic family of super Riemann surfaces over $\Abb^{0|n}_\Cbb$ with complex structure $\widetilde\pt_\xi  = \widetilde\pt - \vp(\xi)$. Then the extension class $\Theta(T^*_{\Xc^{\mathrm{an.}}, -})$ satisfies,
\begin{align}
Dol~\Theta(T^*_{\Xc^{\mathrm{an.}}, -})
\equiv
\sum_i\xi_i\left.\frac{\pt \vp}{\pt \xi_i}\right|_{\xi = 0}
\label{fjvbebvyubviu4r98h30}
\end{align}
where the above equivalence is taken modulo the superconformal structure on $S(C, T^*_{C, -})^\8\subset \Xc^\8$.
\end{LEM}
\vspace{\dp0}
\end{minipage}
\end{center}
\emph{Proof of Lemma \ref{jndrnuirreerrr}}.
Recall that holomorphic objects on $\Xc^\8$ lie in the kernel of $\widetilde\pt_\xi$. The terms $X_\al$ and $X_\be$ in \eqref{kfnkvbhjyugvi4hroijproee3} are holomorphic on $\Xc^\8$ and therefore can  be written as smooth sections in the kernel of $\widetilde\pt_\xi$. With $\widetilde \pt_\xi$ as in \eqref{rg874gf8798h80f30f3} and $\vp(\xi)$ as in \eqref{tg74g94hf8hf033}, see that:
\begin{align*}
\widetilde\pt_\xi X_\al = 0
&\iff
\sum_i \xi_i \left(\widetilde\pt f^i(x, \widetilde x) - \chi^i(x, \widetilde x) \right)\frac{\pt}{\pt \q} = 0
~\mbox{and};
\\
\widetilde\pt_\xi X_\be = 0
&\iff
\sum_i \xi_i \left(\widetilde\pt g^i(y, \widetilde y) - \chi^i(y, \widetilde y) \right)\frac{\pt}{\pt \eta} = 0.
\end{align*}
Hence, when restricted to $C^\8$, we see that $\overline \pt f^i = \chi^i(x, \overline x)\frac{\pt}{\pt \q}$ and $\overline \pt g^i = \chi^i(y,\overline y)\frac{\pt}{\pt \eta}$. Since $\chi^i$ is a globally defined, smooth vector $(0, 1)$-form valued in $T_-C^\8$ we can equate $\chi^i(x, \overline x) = \chi^i(y, \overline y)$ on the intersection $\Uc^\8_\al\cap \Uc^\8_\be\cap C^\8$, giving $\overline \pt X_{\al\be} = 0$. Hence, the restriction of $X_{\al\be}$ to $\Uc^\8_\al\cap \Uc^\8_\be\cap C^\8$ is holomorphic on $C^\8$. The image of this 1-cocycle $(X_{\al\be})_{\al,\be}$ under the Dolbeault isomorphism is the closed, $(0, 1)$-form represented over $\Uc^\8_\al\cap C^\8$ by $\overline \pt X_\al$. We find therefore:
\[
Dol~(X_{\al\be})_{\al,\be} 
= 
\sum_i \xi_i~\chi^i(x, \overline x)\frac{\pt}{\pt \q}
=
\frac{1}{2}\sum_i\xi_i\left.\frac{\pt \vp(\xi)}{\pt \xi_i}\right|_{\xi = 0}
\mod 
D_{(x|\q)}~\mbox{(c.f., \eqref{rfh7rttrgf78hf893h0}).}
\]
This confirms \eqref{fjvbebvyubviu4r98h30} and completes the proof of Lemma \ref{jndrnuirreerrr}.
\qed

\subsection{The Inverse Image of $\vp$}
In Lemma \ref{jndrnuirreerrr} we obtained an expression for the image of the extension class $\Theta(T^*_{\Xc^{\mathrm{an.}}, -})$ under $Dol$. Since $Dol$ is an isomorphism, $\Theta(T^*_{\Xc^{\mathrm{an.}}, -})$ can be written as the inverse image of $\vp$ under $Dol$ in an appropriate sense. Regarding this object $\vp(\xi)$, note that it defines a smooth, vector $(0, 1)$-form on $C^\8$, valued in $T_{\mathrm{even}}\Xc^\8|_{C^\8}$. As such, it will be $\overline\pt$-closed and so we can form the sheaf cohomology class $Dol^{-1}\vp$. This is related to the extension class $\Theta(T^*_{\Xc^{\mathrm{an.}}, -})$ as follows:
\begin{center}
\begin{minipage}[t]{0.7\linewidth}
\vspace{0pt}
\begin{LEM}\label{jndrnuir88y78reerrr} 
\[
\Theta(T^*_{\Xc^{\mathrm{an.}}, -})
=
\sum_i \xi_i\left.\frac{\pt Dol^{-1}\vp}{\pt \xi_i}\right|_{\xi=0}
\]
\end{LEM}
\vspace{\dp0}
\end{minipage}
\end{center}
\emph{Proof of Lemma \ref{jndrnuir88y78reerrr}.}
In the proof of Lemma \ref{hf8h8fh30f390jf93jf39} the map $\sum_i\xi_i\pt(-)/\pt \xi_i|_{\xi=0}$ was interpreted as the projection $pr_*$ in \eqref{rfgyeuibiunoifworfeo}. At the level of $(0, 1)$-forms, it gives a projection $pr_*: A^{0,1}(T_{\mathrm{even}}\Xc^\8) \ra A^{0, 1}(\Hom(T^*C^\8, T^*_-C^\8))$. On inspection of \eqref{fjbvkfbrbiuvbovnioenvpe}, see that $T_{\mathrm{even}}\Xc^\8\cong \Hom(T^*C^\8, \wedge^2T_-^*\Xc^\8)$ modulo the ideal $(\xi^3)$. And modulo the ideal $(\xi^2)$ we recover $A^{0, 1}(\oplus^n\Hom(T^*C^\8, T^*C^\8))$ as can be seen from the second exterior power applied to the exact sequence in \eqref{fjvbrvbiurvhoinvpoemve} at the level of smooth vector bundles (c.f., \eqref{rfh7f97hf98h340f093jf03}). By construction of $pr_*$, we have a commutative diagram:
\begin{align}
\xymatrix{
A^{0, 1}(T_{\mathrm{even}}\Xc^\8)
\ar[drr]^{pr_*}
\ar[d]_{\mod (\xi^3)}
\\
A^{0, 1}\big(\Hom(T^*C^\8, \wedge^2T^*_-\Xc^\8)\big)\ar[rr]_{\mod(\xi^2)} & &  A^{0, 1}(\oplus^n\Hom(T^*C^\8, T^*_-C^\8))
}
\label{fknvkjfbvurborvi3ipmopmve}
\end{align}
That is, $pr_*$ of a $T_{\mathrm{even}}\Xc^\8$-valued, $(0, 1)$-form is equivalent to reducing that form by the ideal $(\xi^2)$; and this is equivalent to applying $\sum_i\xi_i\pt(-)/\pt\xi_i|_{\xi=0}$. Now since any $(0, 1)$-form on a Riemann surface will be $\overline \pt$-closed, the Dolbeault isomorphism can be applied to relate vector $(0, 1)$ forms with sheaf cohomology classes. We obtain therefore a commutative diagram:
\begin{align}
\xymatrix{
\ar[d]_{Dol^{-1}}
A^{0, 1}\big(\Hom(T^*C^\8, \wedge^2T^*_-\Xc^\8)\big)\ar[rr]^{pr_*} & &  A^{0, 1}(\oplus^n\Hom(T^*C^\8, T^*_-C^\8))\ar[d]^{Dol^{-1}}
\\
H^1(\mathcal Hom(T^*_C, \wedge^2T^*_{\Xc^{\mathrm{an.}}, -}))
\ar[rr]_{pr_*} & & 
H^1(\oplus^nT_{C, -})
}
\label{klfnvkbvevbiuoeeee}
\end{align}
In recalling that $\Theta(T^*_{\Xc^{\mathrm{an.}}, -})\in H^1(C, \oplus^nT_{C, -})$ we now have:
\begin{align}
\Theta(T^*_{\Xc^{\mathrm{an.}}, -})
&=
Dol^{-1}\left(\sum_i\xi_i\left.\frac{\pt \vp}{\pt \xi_i}\right|_{\xi = 0}\right)
&&\mbox{(from Lemma \ref{jndrnuirreerrr})}
\notag
\\
&=
Dol^{-1}pr_*\vp
\notag
\\
&=
pr_*Dol^{-1}\vp&&\mbox{(from \eqref{klfnvkbvevbiuoeeee})}
\label{kfnvfjvkjbvuoivn3ovmpo3m}
\\
\notag
&=
\sum_i \xi_i\left.\frac{\pt Dol^{-1}\vp}{\pt \xi_i}\right|_{\xi=0},
\end{align}
which is precisely the statement of Lemma \ref{jndrnuir88y78reerrr}.
\qed

\section{The Superconformal Model}

\noindent
Any model $(X, T^*_{X, -})$ for a class of complex supermanifolds, $X$ is referred to as the \emph{reduced space} and $T^*_{X, -}$ is referred to as the \emph{odd cotangent bundle}. In this section make use Lemma \ref{jndrnuirreerrr} and Lemma \ref{jndrnuir88y78reerrr} to comment the odd cotangent bundle $T^*_{\Xc^{\mathrm{an.}}, -}$ up to isomorphism.

\subsection{The Odd Cotangent Bundle $T^*_{\Xc^{\mathrm{an.}}, -}$}
Starting from the smooth supermanifold $\Xc^\8 = S(C, T^*_{C, -})^\8\times\Abb_\Cbb^{0|n}$ with complex structure $\widetilde\pt_\xi$, we have the difference $\vp(\xi) = \widetilde\pt_\xi - \widetilde\pt$. Note that $\pt\vp/\pt\xi_i|_{\xi = 0}$ will be a closed, vector $(0, 1)$-form valued in $T_-C^\8$. Its pre-image then under the Dolbeault isomorphism gives a class $Dol^{-1}\pt\vp/\pt\xi_i|_{\xi=0} \in H^1(C, T_{C, -})$. Now for any locally free sheaf $\Fc$ on a complex space with structure sheaf $\Oc$ there is an isomorphism $H^1(\Fc^*)\cong \mathrm{Ext}^1(\Fc, \Oc)$, where $\Fc^* =\mathcal Hom_\Oc(\Fc, \Oc)$. Hence, associated to any class in $H^1(\Fc^*)$ will be an extension of holomorphic bundles $0\ra \Oc\ra\Fc^{\p*}\ra\Fc^*\ra0$. Using that $H^1(C, T_{C, -})\cong \mathrm{Ext}^1(T^*_{C, -}, \Oc_C)$, we obtain from $Dol^{-1}\pt\vp/\pt\xi_i|_{\xi=0}$ an extension $0\ra \Oc_C\ra T^*_{i,-}\ra T^*_{C, -}\ra0$ for each $i$. The Dolbeault pre-image $Dol^{-1}(\sum_i\xi_i\pt\vp/\pt\xi_i|_{\xi=0})$ corresponds then to an extension $0\ra \oplus^n\Oc_C\ra T^*_-\ra T^*_{C, -}\ra0$. In this way, from the data of a complex structure $\widetilde \pt_\xi$ on $\Xc^\8 = S(C, T^*_{C, -})^\8\times\Abb^{0|n}_\Cbb$, we can form a model $(C, T^*_-)$ on which to build a class of supermanifolds. Regarding the odd cotangents $T^*_-$ we have:
\begin{center}
\begin{minipage}[t]{0.7\linewidth}
\vspace{0pt}
\begin{LEM}\label{jndrnuied3f3f3rreerrr} 
There exists an isomorphism between $T^*_-$ and $T^*_{\Xc^{\mathrm{an.}}, -}$.
\end{LEM}
\vspace{\dp0}
\end{minipage}
\end{center}
\emph{Proof of Proposition \ref{jndrnuied3f3f3rreerrr}}.
Recall from Lemma \ref{jndrnuirreerrr} that $T^*_{\Xc^{\mathrm{an.}}, -}$ will define an extension class $\Theta(T^*_{\Xc^{\mathrm{an.}}, -})\in \mathrm{Ext}^1(T^*_{C, -},\oplus^n\Oc_C)$, similarly to $T^*_-$.  To show $T^*_{\Xc^{\mathrm{an.}}, -}\cong T^*_-$ then, it suffices to show that these extension classes coincide. Regarding the extension class of $T^*_-$, we have by construction,
\begin{align}
Dol^{-1}\sum_i\xi_i\left.\frac{\pt\vp}{\pt \xi_i}\right|_{\xi=0}
=
\Theta(T^*_-).
\label{rnvjhrvugf938hion3v378v83i3}
\end{align}
We now have the following deductions:
\begin{align*}
\Theta(T^*_{\Xc, -}) 
&=
Dol^{-1}Dol~\Theta(T^*_{\Xc, -})
\\
&=
Dol^{-1}\sum_i \xi_i\left.\frac{\pt \vp}{\pt\xi_i}\right|_{\xi=0}
&&\mbox{(by Lemma \ref{jndrnuirreerrr})}
\\
&=
\Theta(T^*_-)&&\mbox{(from \eqref{rnvjhrvugf938hion3v378v83i3})}
\end{align*}
This completes the proof.
\qed

\subsection{A Decomposition and Superconformality}
We now consider the consequences of Lemma \ref{jndrnuir88y78reerrr}. Recall that we have isomorphisms 
\[
\mathrm{Ext}^1(T^*_{C, -}, \oplus^n\Oc_C)\cong H^1(\oplus^nT_{C, -})\cong \oplus^nH^1(T_{C, -})\cong \oplus^n\mathrm{Ext}^1(T^*_{C, -}, \Oc_C).
\] 
Therefore, associated to $T^*_-$ is a tuple of extensions $(T^*_{i, -})_{i = 1,\ldots, n}$. 
\begin{center}
\begin{minipage}[t]{0.7\linewidth}
\vspace{0pt}
\begin{LEM}\label{jndrnuied3f222rreerrr} 
There exists an isomorphism between $T^*_-$ and $\oplus_iT^*_{i, -}$ where $T^*_{i, -}$ is an extension of $T^*_{C, -}$ by $\Oc_C$ with extension class
\[
\Theta(T^*_{i, -}) = \xi_i~Dol^{-1}\left.\frac{\pt \vp}{\pt \xi_i}\right|_{\xi=0}
\]
\end{LEM}
\vspace{\dp0}
\end{minipage}
\end{center}
\emph{Proof of Lemma \ref{jndrnuied3f222rreerrr}}.
As with $T^*_-$, we define the extension $T^*_{i, -}$ up to isomorphism, and so by its extension class 
\begin{align}
\xi_i~Dol^{-1}\pt\vp/\pt\xi_i|_{\xi=0}
=
\Theta(T^*_{i, -}) 
\label{rfjh4uthg974hg804j333322}
\end{align}
Recall now the extension $T^*_{\Xc^{\mathrm{an.}}, -}$. By Lemma \ref{jndrnuir88y78reerrr} (and commutativity of the diagram \eqref{klfnvkbvevbiuoeeee} in particular), we can move $Dol^{-1}$ past $\xi_i$ and $\pt/\pt\xi_i$. We therefore have
\begin{align*}
\Theta(T^*_-) &= \Theta(T^*_{\Xc^{\mathrm{an.}}, -})&&\mbox{(by Lemma \ref{jndrnuied3f3f3rreerrr})}
\\
&=
\sum_i\xi_i\left.\frac{\pt Dol^{-1}\vp}{\pt \xi_i}\right|_{\xi = 0} 
&&\mbox{(by Lemma \ref{jndrnuir88y78reerrr})}
\\
&=
\sum_i\xi_i~Dol^{-1}\left.\frac{\pt \vp}{\pt \xi_i}\right|_{\xi = 0} 
&&\mbox{(by \eqref{klfnvkbvevbiuoeeee}, c.f., \eqref{kfnvfjvkjbvuoivn3ovmpo3m})}
\\
&=
\sum_i \Theta(T^*_{i, -})
&&\mbox{(by \eqref{rfjh4uthg974hg804j333322})}.
\end{align*}
This completes the proof.\qed
\\\\
We can now conclude the following.
\begin{center}
\begin{minipage}[t]{0.7\linewidth}
\vspace{0pt}
\begin{PROP}\label{jndrnuidddf222rreerrr} 
Let $\Xc^\8 = (S(C, T^*_{C, -})^\8\times\Abb^{0|n}_\Cbb, \widetilde\pt_\xi)$ be an analytic deformation. Then $(C, T^*_{\Xc^{\mathrm{an.}}, -})$ will be a superconformal model.
\end{PROP}
\vspace{\dp0}
\end{minipage}
\end{center}
\emph{Proof of Lemma \ref{jndrnuidddf222rreerrr}.} 
This is a consequence of the relation between $\Theta(T^*_{i, -})$ and $\vp(\xi) = \widetilde\pt_\xi - \widetilde\pt$ in Lemma \ref{jndrnuied3f222rreerrr}.
\qed

\section{A Superconformal Atlas}
\label{f9488h4h04j08408jv0}

\subsection{To Second Order (Ansatz)}
Recall that an algebraic deformation is, by Definition \ref{rf4fg784hf98h0f3}, a superconformal atlas. Associated to any atlas $(\Uc,\vartheta)$ is an abelian cohomology class $\om_{(\Uc, \vartheta)}$ which, in Theorem \ref{34fiufguhfojp3}, was referred to as an obstruction class. Schematically, we write
\begin{align}
\vartheta = {\bf 1} + \om_{(\Uc, \vartheta)} +\cdots
\label{rbugg4gh983hf03jf093j}
\end{align}
to indicate an atlas $(\Uc, \vartheta)$ with obstruction $\om_{(\Uc, \vartheta)}$. With $(\Uc, \vartheta)$ an atlas for a complex supermanifold $\Xc$ modelled on $(X, T^*_{X, -})$, the term ${\bf 1}$ in \eqref{rbugg4gh983hf03jf093j} represents the data coming from the model $(X, T^*_{X, -})$, i.e., is represented by the glueing data of the underlying space $X$ and the sheaf $T^*_{X, -}$ (c.f., \eqref{rbvyrvyuev7998hf083h0} and \eqref{fgjhbhbhbwuodn}). Note that the expression in \eqref{rbugg4gh983hf03jf093j} makes sense more generally for any appropriate, abelian $1$-cochain as long as it is such that $\vartheta$ will satisfy the cocycle conditions.

With this preamble for the notation we have:
\begin{center}
\begin{minipage}[t]{0.7\linewidth}
\vspace{0pt}
\begin{LEM}\label{jndrf4f4f4g4rnuir88y78reerrr} 
Let $\widetilde\pt_\xi = \widetilde\pt - \vp(\xi)$ be an analytic deformation of a super Riemann surface $\Scl$ over $\Abb^{0|n}_\Cbb$. Then any atlas $(\Uc, \vartheta)$ for the total space of a deformation $\Xc\ra \Abb^{0|n}_\Cbb$ of $\Scl$ given by the ansatz
\begin{align*}
\vartheta
=
{\bf 1} + 
\sum_i \xi_i\left.\frac{\pt Dol^{-1}\vp}{\pt \xi_i}\right|_{\xi=0}
+
\frac{1}{2}
\sum_{i<j}\xi_i\xi_j
\left.\frac{\pt^2 Dol^{-1}\vp}{\pt \xi_j\pt\xi_i}\right|_{\xi=0}
+\cdots
\end{align*}
will be superconformal to second order.
\end{LEM}
\vspace{\dp0}
\end{minipage}
\end{center}
\emph{Proof of Lemma \ref{jndrf4f4f4g4rnuir88y78reerrr}}.
Since $\widetilde\pt_\xi$ is an analytic deformation, $\vp$ will be an even, superconformal vector $(0, 1)$-form. Suppose it is given by \eqref{tg74g94hf8hf033} and let $\om = Dol^{-1}\vp$ modulo $(\xi^3)$. Then on intersections $\om$ will be represented by a difference of smooth, superconformal vector fields. We will firstly show that this difference will itself be a $1$-cocycle valued in holomorphic, superconformal vector fields. To begin, let $(\Uc^\8_\al)$ be a covering for $\Xc^\8 = S(C, T^*_{C, -})^\8\times \Abb^{0|n}_\Cbb$.
We can represent $\om$ on $\Uc^\8_\al\cap \Uc^\8_\be$ by,
\begin{align}
\om_{\al\be}
&=
\left(\nu_\be - \nu_\al\right)
+
\left(X_\be - X_\al\right)
\label{vruibrutih093jvj30}
\end{align}
where $\nu_\al, \nu_\be$ are odd; and $X_\al, X_\be$ are even, smooth superconformal vector fields. By construction we have $\overline \pt\nu_\al = \overline\pt\nu_\be = \sum_i\xi_i~\pt Dol^{-1}\vp/\pt\xi_i|_{\xi=0}$; and similarly $\overline \pt X_\be = \overline \pt X_\al = \sum_{i<j}\xi_i\xi_j~\pt Dol^{-1}\vp/\pt \xi_j\pt\xi_i|_{\xi=0}$. Now with the super Riemann surface $\Scl = S(C, T^*_{C, -})$ given, let $(\Ufr_\al, \rho)$ be a covering. If we denote by $(x, \widetilde x|\q, \xi)$ and $(y, \widetilde y|\eta, \xi)$ coordinates on $\Uc^\8_\al$ and $\Uc^\8_\be$ respectively, then $(x|\q)$ and $(y|\eta)$ will be local, holomorphic coordinates on $\Scl$ and so coordinates on $\Ufr_\al$, resp., $\Ufr_\be$. On intersections,
\begin{align*}
y = \rho^+_{\al\be}(x|\q) = f_{\al\be}(x)
&&\mbox{and}&&
\eta = \rho_{\al\be}^-(x|\q) = \zeta_{\al\be}(x)\q.
\end{align*}
Here, $(f_{\al\be})$ are the (holomorphic) transition data for the underlying space $C$; while $(\zeta_{\al\be})$ are the transition data for the spin structure $T^*_{C, -}$. Since $T^*_{C, -}$ is a spin structure, we have the relation $\zeta^2_{\al\be} = \pt f_{\al\be}/\pt x$ (c.f., Theorem \ref{fjbvbviueovinep} and \eqref{fjbvrbviubvonoioie1}). With this relation we will now show that $\om_{\al\be}$ in \eqref{vruibrutih093jvj30} will define a sum of even and odd superconformal vector fields. Starting with the odd superconformal vector fields, write:
\begin{align*}
\nu_\al = \sum_i\xi_i~v_\al^i(x, \widetilde x)\left(\frac{\pt}{\pt \q} - \q\frac{\pt}{\pt x}\right)
&&
\mbox{and}
&&
\nu_\be = \sum_i \xi_i~v_\be^i(y, \widetilde y)\left(\frac{\pt}{\pt \eta} - \eta\frac{\pt}{\pt y}\right).
\end{align*}
Observe that over intersections $\Ufr_\al\cap \Ufr_\be$ on the super Riemann surface $S(C, T^*_{C, -})$,
\begin{align*}
\frac{\pt}{\pt \q} = \zeta_{\al\be}(x)\frac{\pt}{\pt\eta}
&&
\mbox{and}
&&
\frac{\pt}{\pt x}
=
\frac{\pt f_{\al\be}}{\pt x}\frac{\pt}{\pt y} + \frac{\pt \zeta_{\al\be}}{\pt x}\q\frac{\pt}{\pt \eta}
=
\zeta_{\al\be}(x)^2\frac{\pt}{\pt y} + \frac{\pt \zeta_{\al\be}}{\pt x}\q\frac{\pt}{\pt \eta}.
\end{align*}
Hence we find
\begin{align}
\nu_\be - \nu_\al
&=
\sum_i \xi_i
\Big(
v_\be(y,\widetilde y)
-
\zeta_{\al\be}(x)v_\al(x,\widetilde x)
\Big)
\left(\frac{\pt}{\pt \eta} - \eta\frac{\pt}{\pt y}\right)
\label{fjbvrbiuvoihip3j9j39p}
\end{align}
and so the difference $\nu_\be-\nu_\al$ will define a $1$-cochain valued in the sheaf of superconformal vector fields. It is holomorphic since, by assumption, $\overline \pt\om_{\al\be} = 0$. Regarding the difference $X_\be - X_\al$ of even superconformal vector fields, write:
\begin{align*}
X_\al &= \sum_{i<j}\xi_i\xi_j\left(X_\al^{ij}(x,\widetilde x)\frac{\pt}{\pt x} + \frac{1}{2}\frac{\pt X_\al^{ij}}{\pt x}\q\frac{\pt}{\pt \q}\right)
\\
X_\be &= \sum_{i<j}\xi_i\xi_j\left(X_\be^{ij}(y,\widetilde y)\frac{\pt}{\pt y} + \frac{1}{2}\frac{\pt X_\be^{ij}}{\pt y}\eta\frac{\pt}{\pt \eta}\right)
\end{align*}
As with the even case we have,
\begin{align}
X_\al 
&= \sum_{i< j}\xi_i\xi_j
\left(
X^{ij}_\al\left(\zeta_{\al\be}^2\frac{\pt}{\pt y} + \frac{\pt\zeta_{\al\be}}{\pt x}\zeta_{\al\be}^{-1}\eta\frac{\pt}{\pt \eta}\right)
+
\frac{1}{2}\frac{\pt X^{ij}_\al}{\pt x} \eta\frac{\pt}{\pt \eta}
\right)
\notag
\\
&=
\sum_{i<j}\xi_i\xi_j
\left(
(X^{ij}_\al\zeta^2_{\al\be})\frac{\pt}{\pt y}
+
\frac{1}{2}\left(X^{ij}_\al\zeta_{\al\be}^{-2}\frac{\pt \zeta^2_{\al\be}}{\pt x}
+
\frac{\pt X_\al^{ij}}{\pt x}
\right)
\eta\frac{\pt}{\pt \eta}
\right)
\label{fkvnkrbviubvounoincieonveou}
\\
&=
\sum_{i<j}\xi_i\xi_j
\left(
(X^{ij}_\al\zeta^2_{\al\be})\frac{\pt}{\pt y}
+
\frac{1}{2}\zeta_{\al\be}^{-2}\frac{\pt \big(X^{ij}_\al\zeta_{\al\be}^2\big)}{\pt x}
\eta\frac{\pt}{\pt \eta}
\right)
\notag
\\
&=
\sum_{i<j}\xi_i\xi_j
\left(
(X^{ij}_\al\zeta^2_{\al\be})\frac{\pt}{\pt y}
+
\frac{1}{2}\frac{\pt x}{\pt y} \frac{\pt \big(X^{ij}_\al\zeta_{\al\be}^2\big)}{\pt x}
\eta\frac{\pt}{\pt \eta}
\right)
\label{kfndbvbvuibevoen}
\\
\notag
&=
\sum_{i<j}\xi_i\xi_j
\left(
(X^{ij}_\al\zeta^2_{\al\be})\frac{\pt}{\pt y}
+
\frac{1}{2} \frac{\pt \big(X^{ij}_\al\zeta_{\al\be}^2\big)}{\pt y}
\eta\frac{\pt}{\pt \eta}
\right)
\end{align}
where in \eqref{fkvnkrbviubvounoincieonveou} we used $\frac{\pt\zeta_{\al\be}}{\pt x}\frac{1}{\zeta_{\al\be}} = \frac{1}{2}\frac{\pt \zeta^2_{\al\be}}{\pt x}\frac{1}{\zeta^2_{\al\be}}$; and in \eqref{kfndbvbvuibevoen} we used $\zeta^{-2}_{\al\be} = (\pt f_{\al\be}/\pt x)^{-1} = \pt x/\pt y$. Hence as in \eqref{fjbvrbiuvoihip3j9j39p} we find, 
\begin{align}
X_\be
-
X_\al
&=
\sum_{i<j}
\xi_i\xi_j
\left(
\big(X^{ij}_\be - X^{ij}_\al\zeta_{\al\be}^2\big)
+
\frac{1}{2}\frac{\pt \big(X^{ij}_\be - X^{ij}_\al\zeta_{\al\be}^2\big)}{\pt y}\eta\frac{\pt}{\pt \eta}
\right)
\label{fjvbrbvybvuvoi3nvpinuib}
\end{align}
and so the difference $X_\be - X_\al$ will be a $1$-cochain valued in even, superconformal vector fields. From the revelations in \eqref{fjbvrbiuvoihip3j9j39p} and \eqref{fjvbrbvybvuvoi3nvpinuib}, we can write,
\begin{align}
\om_{\al\be}
&\stackrel{\Delta}{=}
\left\{\sum_i \xi_i\left.\frac{\pt Dol^{-1}\vp}{\pt \xi_i}\right|_{\xi=0}
+
\sum_{i<j}\xi_i\xi_j
\left.\frac{\pt^2 Dol^{-1}\vp}{\pt \xi_j\pt\xi_i}\right|_{\xi=0}\right\}_{\al\be}
\notag
\\
&=
\sum_i\xi_i~\psi^i_{\al\be}(x)
\left(\frac{\pt}{\pt \eta} - \eta\frac{\pt}{\pt y}\right)
+
\sum_{i<j} \xi_i\xi_j
\left(g_{\al\be}^{ij}\frac{\pt}{\pt y}
+
\frac{1}{2}\frac{\pt g^{ij}_{\al\be}}{\pt y}\eta\frac{\pt}{\pt \eta}\right)
\label{fjvbkrbuibvui4oi4nvpi3nb}
\end{align}
for some holomorphic functions $\psi^i_{\al\be}$ and $g_{\al\be}^{ij}$. Note that since $\widetilde\pt_\xi$ was assumed to define an analytic deformation, Proposition \ref{jndrnuidddf222rreerrr} will imply that $\{\pt Dol^{-1}\vp/\pt \xi_i|_{\xi=0}\}_i$ are pairwise linearly dependent (c.f., Proposition \ref{fjbvbrviueovinep}). In comparing now with the relations in \eqref{fjbvrbviubvonoioie1}---\eqref{fjbvrbviubvonoioie13} it becomes evident that any atlas $(\Uc, \vartheta)$ with $\vartheta = {\bf 1} + \om_{(\Uc, \vartheta)} + \cdots$, where $(\om_{(\Uc, \vartheta)})_{\al\be} = \om_{\al\be}$ for all $\al, \be$, will be superconformal modulo the ideal $(\xi^3)$, i.e., to second order.
\qed

\subsection{Existence}
In Lemma \ref{jndrf4f4f4g4rnuir88y78reerrr} we assumed the existence of an atlas $(\Uc,\vartheta)$ for the total space of a deformation $\Xc\ra\Abb^{0|n}_\Cbb$ with $\vartheta$ given by the ansatz in Lemma \ref{jndrf4f4f4g4rnuir88y78reerrr}. Superconformality of $(\Uc, \vartheta)$ then followed. Presently we argue that there will always exist a supermanifold with such an atlas associated to any analytic deformation.
\begin{center}
\begin{minipage}[t]{0.7\linewidth}
\vspace{0pt}
\begin{PROP}\label{jndrf4g4rnuir88y78reerrr} 
Let $\widetilde\pt_\xi = \widetilde\pt - \vp(\xi)$ be an analytic deformation of a super Riemann surface $\Scl$ over $\Abb^{0|n}_\Cbb$. Then there will exist some atlas $(\Uc, \vartheta^\p)$ for the total space of a deformation $\Xc^\p\ra \Abb^{0|n}_\Cbb$ of $\Scl$ such that
\begin{align*}
\vartheta^\p
=
{\bf 1} + 
\sum_i \xi_i\left.\frac{\pt Dol^{-1}\vp}{\pt \xi_i}\right|_{\xi=0}
+
\frac{1}{2}
\sum_{i<j}\xi_i\xi_j
\left.\frac{\pt^2 Dol^{-1}\vp}{\pt \xi_j\pt\xi_i}\right|_{\xi=0}
+\cdots
\end{align*}
\end{PROP}
\vspace{\dp0}
\end{minipage}
\end{center}
\emph{Proof of Proposition \ref{jndrf4g4rnuir88y78reerrr}.}
Starting from the expression of $\om$ on intersections $U_\al\cap U_\be$ in \eqref{fjvbkrbuibvui4oi4nvpi3nb}, note that the following part
\[
\om^\p_{\al\be} = \sum_i\eta\xi_i~\psi^i_{\al\be}\frac{\pt}{\pt y} + \sum_{i<j}\xi_i\xi_j~g^{ij}_{\al\be}\frac{\pt}{\pt y}
\]
will define a $1$-cocycle valued in the second obstruction space of the model $(C, T^*_{-})$; where $T_{-}^*$ is an extension of holomorphic bundles $T^*_{C, -}$ and $\oplus^n\Oc_C$, with extension class $\Theta(T^*_{-})$ represented by the $1$-cocycle $\Theta(T^*_{-})_{\al\be} = \sum_i\xi_i~\psi^i_{\al\be}\pt/\pt \eta$. By Theorem \ref{h8937893h0fj309j03} there will exist a supermanifold $\Xc$ modelled on $(C, T^*_{-})$ with atlas $(\Uc, \vartheta)$ such that $\om_{(\Uc, \vartheta)} = \om^\p$. Denote by $\Xc^{(2)}\subset \Xc$ the second infinitesimal neighbourhood of $C\subset \Xc$. As a locally ringed space $\Xc^{(2)} = (C, \Oc_\Xc/J^3)$ where $J\subset \Oc_\Xc$ is the fermionic ideal.\footnote{This is the ideal generated by the odd, nilpotent elements.} Now set 
\begin{align*}
\vartheta^\p 
&\stackrel{\Delta}{=} \left(\vartheta\mod (\xi^2)\right) + \sum_{i<j}\xi_i\xi_j~\left.\frac{\pt^2Dol^{-1}\vp}{\pt\xi_j\pt \xi_i}\right|_{\xi=0}
\\
&=
\left(\vartheta_{\al\be}\mod (\xi^2)\right)
+
\sum_{i<j} \xi_i\xi_j\left(g_{\al\be}^{ij}\frac{\pt}{\pt y}  + \frac{1}{2}\frac{\pt g_{\al\be}^{ij}}{\pt x}\eta\frac{\pt}{\pt \eta}\right)
&&\mbox{(on intersections $U_\al\cap U_\be$)}.
\end{align*}
Any obstructions to enforcing the cocycle condition for $\vartheta^\p$ are detected by second-degree cohomology. As there exist no such cohomology on Riemann surfaces, there are therefore no obstructions to enforcing the cocycle condition for $\vartheta^\p$. As such $(\Uc, \vartheta^\p)$ defines an atlas for an infinitesimal thickening $\Xc^{\p(3)}$ of $\Xc^{(2)}\subset \Xc$. As any thickening over a Riemann surface is unobstructed, it will embedd in some supermanifold $\Xc^\p$. Hence there will exist some supermanifold $\Xc^\p$ with atlas $(\Uc, \vartheta^\p)$ such that $\vartheta^\p$ is given by the expression in the statement of this proposition. This completes the proof.
\qed

\section{Correspondence of Equivalence Classes}

\noindent
We have so far obtained a correspondence of representative deformations. That is, given an algebraic deformation $(\Uc, \vartheta)$, Proposition \ref{hg975hg984hg8003jg03} shows how we can construct an analytic deformation $\widetilde\pt_\xi$. Conversely, given an analytic deformation $\widetilde\pt_\xi$, we see how to reconstruct an algebraic deformation $(\Uc, \vartheta)$ in Lemma \ref{jndrf4f4f4g4rnuir88y78reerrr} and Proposition \ref{jndrf4g4rnuir88y78reerrr}. In order to finish the proof of Theorem \ref{rfg748gf7f983h04i4i4i43hf03} now, it remains to show that this correspondence will preserve the respective equivalences, detailed in Theorem \ref{fhvurh84j093j3j33333} and Theorem \ref{rfhf9h3f8309fj3fjoooo3} for the analytic and algebraic deformations respectively.

\subsection{The Linear Term}
We firstly consider the terms in the respective deformations which are linear with respect to the auxiliary parameters $(\xi_i)$. For an algebraic deformation $(\Uc, \vartheta)$ this term is the $1$-cocycle $\{\psi^i_{\al\be}\frac{\pt}{\pt\eta_i}\}_{\al,\be}$; and for an analytic deformation it is the odd, superconformal vector $(0, 1)$-form $\chi^i(x, \widetilde x)\left(\frac{\pt}{\pt\q} - \q\frac{\pt}{\pt x}\right)$. Recall now the following correspondences of deformations: by Proposition \ref{hg975hg984hg8003jg03} (c.f., \eqref{rgf64gf783hf983hf08j3093})  we have
\begin{align}
\left\{\psi^i_{\al\be}\frac{\pt}{\pt\eta}\right\}_{\al, \be}
&\equiv 
\frac{1}{2}\left\{\psi_{\al\be}^i\left(\frac{\pt}{\pt \eta} - \eta\frac{\pt}{\pt y}\right)\right\}
\notag
\\
&\longmapsto 
\frac{1}{2}
\left.\frac{\pt Dol~\om_{(\Uc, \vartheta)}}{\pt \xi_i}\right|_{\xi=0}
=
\frac{1}{2}\chi^i\left(\frac{\pt}{\pt \q }- \q\frac{\pt}{\pt x}\right)
\equiv 
\chi^i\frac{\pt}{\pt\q}
\label{fjbvjrvuyviubvouninoded}
\end{align}
and conversely, from the statements leading up to Proposition \ref{jndrf4g4rnuir88y78reerrr} we have
\begin{align}
\chi^i\frac{\pt}{\pt\q}
&\equiv \frac{1}{2}\chi^i\left(\frac{\pt}{\pt\q} - \q\frac{\pt}{\pt x}\right)
\notag
\\
&\longmapsto 
\frac{1}{2}
\left\{\left.\frac{\pt Dol^{-1}\vp}{\pt \xi_i}\right|_{\xi=0}\right\}_{\al, \be}
=
\frac{1}{2}\left\{\psi^i_{\al\be}\left(\frac{\pt}{\pt\eta} - \eta\frac{\pt}{\pt y}\right)\right\}_{\al, \be}
\equiv 
\left\{\psi^i_{\al\be}\frac{\pt}{\pt\eta}\right\}_{\al,\be}
\label{fjkvbrhvbirbvuenvoinvpon}
\end{align}
where the above equivalences `$\equiv$' are taken modulo the superconformal structure (c.f., \eqref{rfh7rttrgf78hf893h0}). The equivalences in \eqref{fjbvjrvuyviubvouninoded} and \eqref{fjkvbrhvbirbvuenvoinvpon} are made with respect to the superconformal structure on the super Riemann surface $S(C, T^*_{C, -})$. The Dolbeault isomorphism sends coboundary terms in \v Cech cohomology to $\overline \pt$-exact forms. On inspection of \eqref{Pfnvfbvjkvnvoeefe} in Theorem \ref{fhvurh84j093j3j33333} and \eqref{fgyur8g789h30h03} in Theorem \ref{rfhf9h3f8309fj3fjoooo3} we have the following correspondence:
\begin{align}
\left\{w^i_{\be}\otimes \frac{\pt}{\pt \eta}
-
w^i_{\al}\otimes\frac{\pt}{\pt \q}\right\}_{\al, \be} \longleftrightarrow 
\widetilde\pt \nu^i
\label{rg78g78gf97389hf039j093}
\end{align}
Hence, any linear equivalence of an algebraic deformation (i.e., an equivalence modulo $(\xi^2)$) will translate to a linear equivalence of the corresponding analytic deformation and vice-versa.

\subsection{The Quadratic Term}
The quadratic terms in the respective deformations are the terms proportional to $(\xi_i\xi_j)_{i<j}$. For the algebraic deformation $(\Uc,\vartheta)$ it is the $1$-cochain $\{g^{ij}_{\al\be}\frac{\pt}{\pt y}\}_{\al,\be}$. And for the analytic deformation $\widetilde\pt_\xi$ it is the even, superconformal vector $(0, 1)$-form $h^{ij}\frac{\pt}{\pt x} + \frac{1}{2}\frac{\pt h^{ij}}{\pt x}\q\frac{\pt}{\pt\q}$. As with the correspondence of the linear terms in \eqref{fjbvjrvuyviubvouninoded} and \eqref{fjkvbrhvbirbvuenvoinvpon} we have:
\begin{align}
\left\{
g^{ij}_{\al\be}\frac{\pt}{\pt y}
\right\}_{\al, \be}
\longmapsto 
\left.\frac{\pt^2 Dol~\om_{(\Uc, \vartheta)}}{\pt \xi_j\pt\xi_i}\right|_{\xi=0}
=
h^{ij}\frac{\pt}{\pt x} + \frac{1}{2}\frac{\pt h^{ij}}{\pt x}\q\frac{\pt}{\pt\q}
\equiv 
h^{ij}\frac{\pt}{\pt x}
\label{rjbkvy4ivu4ovoicioeneoe}
\end{align}
and conversely (c.f., \eqref{fjvbkrbuibvui4oi4nvpi3nb}),
\begin{align}
h^{ij}&\frac{\pt}{\pt x}
\longmapsto 
\left\{
\left.\frac{\pt^2 Dol^{-1}\vp}{\pt\xi_j\pt\xi_i}\right|_{\xi=0}
\right\}_{\al,\be}
=
\left\{
g^{ij}_{\al\be}\frac{\pt}{\pt y}
+
\frac{1}{2}\frac{\pt g^{ij}_{\al\be}}{\pt y}\eta\frac{\pt}{\pt\eta}
\right\}_{\al,\be}
\equiv
\left\{
g^{ij}_{\al\be}\frac{\pt}{\pt y}
\right\}_{\al,\be}
\label{fjbcevyu4iucneiocnjkebciue}
\end{align}
where the above equivalences `$\equiv$' are taken modulo the superconformal structure on $S(C, T^*_{C, -})$ (c.f., \eqref{rfh7rttrgf78hf893h}). Recall that the transformation of the quadratic terms under equivalence were documented in \eqref{fnvknvkjnivnoienvo} in Theorem \ref{fhvurh84j093j3j33333}; and \eqref{fjkfbvrbvuinoiniwnviev} in Theorem \ref{rfhf9h3f8309fj3fjoooo3}. We begin by recalling the transformation \eqref{fjkfbvrbvuinoiniwnviev} below:
\begin{align}
g^{ij}_{\al\be}\frac{\pt}{\pt y} - \tilde g^{ij}_{\al\be}\frac{\pt}{\pt y}
=&~
u^{ij}_\be\frac{\pt}{\pt y} - u^{ij}_\al\frac{\pt}{\pt x}
\label{fjbvhjrvyu494vh04444}
\\
&~+
\dt\left(\left\{ w^i_\al\frac{\pt}{\pt\q}\right\}\right)_{\al\be}\otimes \psi^j_{\al\be}\frac{\pt}{\pt\eta}
\label{fjnkjbvybivunio3}
\\
&~-
\dt\left(\left\{ w^j_\al\frac{\pt}{\pt\q}\right\}\right)_{\al\be}\otimes \psi^i_{\al\be}\frac{\pt}{\pt\eta}
\label{rbuyrgf4gf7h893}
\\
&~+
\left(w^i_\al w^j_\be - w^j_\al w^i_\be\right)\frac{\pt}{\pt\q}\otimes\frac{\pt}{\pt\eta}.
\label{vuhruh894h0f9j903}
\end{align}
We will now use that $\pt/\pt\eta$ has \emph{odd} parity as a vector field, meaning it will anti-commute with itself. Hence, regarding \eqref{fjnkjbvybivunio3}, we find:
\begin{align}
\dt\left(\left\{ w^i_\al\frac{\pt}{\pt\q}\right\}\right)_{\al\be}\otimes &~\psi^j_{\al\be}\frac{\pt}{\pt\eta}
\notag
\\
&=
\left( w^i_\be \frac{\pt}{\pt \eta} - w^i_\al\frac{\pt}{\pt \q}\right)\otimes \psi^j_{\al\be}\frac{\pt}{\pt \eta}
\notag
\\
&=
\left(w^i_\be \frac{\pt}{\pt \eta}\right)\otimes \left(\psi^j_{\al\be}\frac{\pt}{\pt \eta}\right)
-
\left(w^i_\al \frac{\pt}{\pt \q}\right)\otimes \left(\psi^j_{\al\be}\frac{\pt}{\pt \eta}\right)
\notag
\\
&=
\left[
\left(
w^i_\be\frac{\pt}{\pt \eta}
\right)
\otimes
\left(
\psi^j_{\al\be}\frac{\pt}{\pt \eta}
\right)
+
\left(
\psi^j_{\al\be}\frac{\pt}{\pt \eta}
\right)
\otimes
\left(
w^i_\al\frac{\pt}{\pt \q}
\right)
\right]
\label{fjvbrbgig978984h0j0v9jp4v43v}
\end{align}
Similarly, regarding \eqref{rbuyrgf4gf7h893} we have:
\begin{align}
\dt\left(\left\{ w^j_\al\frac{\pt}{\pt\q}\right\}\right)_{\al\be}\otimes&~ \psi^i_{\al\be}\frac{\pt}{\pt\eta}
\notag\\
&=
\left(\psi^i_{\al\be}\frac{\pt}{\pt \eta}\right)
\otimes 
\left(w^j_\be \frac{\pt}{\pt \eta}\right)
+
\left(w^j_\al \frac{\pt}{\pt \q}\right)\otimes \left(\psi^i_{\al\be}\frac{\pt}{\pt \eta}\right)
\label{fbjhvbhjrvuybiuvb3vh893h}
\end{align}
Finally, where the expression in \eqref{vuhruh894h0f9j903} is concerned, note that it can be rewritten in the following way:
\begin{align}
\eqref{vuhruh894h0f9j903}
&=
\left(w_\al^i \frac{\pt}{\pt\q}\right)\otimes\dt\left(\left\{ w^j_\al\frac{\pt}{\pt\q}\right\}\right)_{\al\be}
-
\left(w^j_\al\frac{\pt}{\pt\q}\right)\otimes\dt\left(\left\{ w^i_\al\frac{\pt}{\pt\q}\right\}\right)_{\al\be}
\label{fknkdbvviubuoeino333ep}
\end{align}
To recover global, smooth vector $(0, 1)$-forms on $C^\8$ we use the correspondences established in \eqref{fjbvjrvuyviubvouninoded}, \eqref{fjkvbrhvbirbvuenvoinvpon} and \eqref{rg78g78gf97389hf039j093}. From this latter-most correspondence \eqref{rg78g78gf97389hf039j093} we can recover a smooth vector field on $C^\8$ from the $0$-cochain $w^k = \{w^k_\al\frac{\pt}{\pt \q}\}_{\al\in I}, k = i, j$ as follows: On each intersection $\Uc^\8_\al\cap \Uc^\8_\be\cap C^\8$ we can write $(\dt w^k)_{\al\be} = \mu^k_\be\frac{\pt}{\pt \eta} - \mu^k_\al\frac{\pt}{\pt \q}$, where $\mu^k = \{\mu^k_\al\frac{\pt}{\pt\q}\}_{\al\in I}$ is a $0$-cochain of smooth vector fields. Rearranging this equation gives the relation: $(w^k_\be - \mu^k_\be)\frac{\pt}{\pt \eta} = (w^k_\al - \mu^k_\al)\frac{\pt}{\pt \q}$ on all intersections $\Uc^\8_\al\cap \Uc^\8_\be\cap C^\8$. Hence there exists a smooth globally defined, odd vector field\footnote{i.e., a smooth section of the bundle $T_-C^\8\ra C^\8$} $\nu^k$ on $C^\8$ such that $\nu^k|_{\Uc^\8_\al\cap C^\8} = (w^k_\al - \mu^k_\al)\frac{\pt}{\pt \q}$. We therefore have: $\nu^k = w^k - \mu^k$ and so a bijective correspondence: $w^k\mapsto w^k - \mu^k$. From this and \eqref{fjkvbrhvbirbvuenvoinvpon} we recover the following smooth vector $(0, 1)$-forms from holomorphic data on intersections:
\begin{align}
\eqref{fjvbrbgig978984h0j0v9jp4v43v}
\longleftrightarrow&~
\nu^i\frac{\pt}{\pt \eta}\otimes \chi^j\frac{\pt}{\pt \eta} + \chi^j\frac{\pt}{\pt \eta}\otimes \nu^i\frac{\pt}{\pt \eta}
\label{fhvurvyu4vbui3ouio3n3ecece}
\\
\equiv&~
\frac{1}{2^2} 
\left[
\nu^i\left(\frac{\pt}{\pt \eta} - \eta\frac{\pt}{\pt y}\right)
,
\chi^j\left(\frac{\pt}{\pt \eta} - \eta\frac{\pt}{\pt y}\right)
\right]
\notag
\\
&~+
\frac{1}{2^2}
\left[
\chi^j\left(\frac{\pt}{\pt \eta} - \eta\frac{\pt}{\pt y}\right)
,
\nu^i\left(\frac{\pt}{\pt \eta} - \eta\frac{\pt}{\pt y}\right)
\right]
\notag
\\
\equiv&~
\frac{1}{2} \left[
\nu^i\left(\frac{\pt}{\pt \eta} - \eta\frac{\pt}{\pt y}\right)
,
\chi^j\left(\frac{\pt}{\pt \eta} - \eta\frac{\pt}{\pt y}\right)
\right]
\label{7g487g94h8093f4444}
\\
=&~
\frac{1}{2}[\nu^i,\chi^j]
\label{rhf784gf78f98h38fj39jf33344}
\end{align} 
where in the lines succeeding \eqref{rhf784gf78f98h38fj39jf33344} we use the correspondence in \eqref{rfh7rttrgf78hf893h0} in addition to the mapping sending tensor products of vector fields to their (super)Lie bracket; in \eqref{7g487g94h8093f4444} that the (super)Lie bracket of \emph{odd} vector fields is symmetric; and in \eqref{rhf784gf78f98h38fj39jf33344} we identify $\nu^i$ and $\chi^j$ with the superconformal, vector $(0, 1)$-forms in the previous line. Similarly, we have:
\begin{align}
\eqref{fbjhvbhjrvuybiuvb3vh893h}
\longleftrightarrow
\frac{1}{2}
\left[\nu^j, \chi^i\right].
\label{rf74gf874hojf908j0fj039jf}
\end{align}
Finally, concerning \eqref{fknkdbvviubuoeino333ep} we have:
\begin{align}
\eqref{fknkdbvviubuoeino333ep}
\longleftrightarrow
\nu^i
\otimes \widetilde\pt \nu^j
-
\nu^j\otimes \widetilde\pt \nu^i
\equiv
\frac{1}{2^2}
\left(
[\nu^i, \widetilde\pt \nu^j]
-
[\nu^j, \widetilde\pt \nu^i]
\right)
\label{rhf94h89h40j390fj3fk3444}
\end{align}
And so from \eqref{rhf784gf78f98h38fj39jf33344}, \eqref{rf74gf874hojf908j0fj039jf} and \eqref{rhf94h89h40j390fj3fk3444} we can conclude:
\begin{align}
\Big\{g^{ij}_{\al\be}\otimes\frac{\pt}{\pt y}
-
\tilde g^{ij}_{\al\be}&\otimes\frac{\pt}{\pt y}
\Big\}_{\al, \be}
\notag
\\
&~\longleftrightarrow
\notag
\\
\widetilde\pt \nu^{ij}
&~+
\frac{1}{2}
\Big(
[\nu^i,\chi^j] - [\nu^j,\chi^i]
\Big)
+
\frac{1}{4}\left(
[\nu^i, \widetilde\pt \nu^j]
-
[\nu^j, \widetilde\pt \nu^i]
\right)
\notag
\\
&= h^{ij\p}\frac{\pt}{\pt y} - h^{ij}\frac{\pt}{\pt y}
\label{rgf764gf784gf98h80f30}
\end{align}
where $\widetilde\pt\nu^{ij}$ is the Dolbeault form of the coboundary $\dt\{u_\al^{ij}\pt/\pt x\}$ in \eqref{fjbvhjrvyu494vh04444}; and where \eqref{rgf764gf784gf98h80f30} is deduced on comparing with the transformation law for $h^{ij}$ in \eqref{fbvrbvyurbvuir98h38h3}. Hence equivalences of algebraic deformations translate to equivalences of analytic deformations, to second order. This concludes the proof of Theorem \ref{rfg748gf7f983h04i4i4i43hf03}.
\qed

\part{Applications}
\label{ppjpvburbviyvbiuvr}

\section{Analytic Deformations and Non-Splitting}

\noindent
For algebraic deformations, recall that the total space $\Xc$ of the deformation will be a complex supermanifold. Hence it makes sense to ask whether these deformations \emph{split} by asking this question of the total space $\Xc$. In this section we consider this problem for analytic deformations.

\subsection{A Non-Split Supermanifold}
Recall from Proposition \ref{jndrnuidddf222rreerrr} that associated to any analytic deformation $\Xc^\8 = (S(C, T^*_{C, -})^\8\times\Abb^{0|n}_\Cbb, \widetilde\pt_\xi)$ is a superconformal model $(C, T^*_{\Xc^{\mathrm{an.}}, -})$. In Proposition \ref{jndrf4g4rnuir88y78reerrr} then, we deduced the existence of a complex supermanifold modelled on  $(C, T^*_{\Xc^{\mathrm{an.}}, -})$ with an atlas of a  certain kind. This deduction forms the basis for the following.
\begin{center}
\begin{minipage}[t]{0.7\linewidth}
\begin{THM}\label{jf0jg43r3jgp3345topgj4gj4} 
Let $\Xc^\8$ be an analytic deformation of a super Riemann surface $(S(C, T^*_{C, -})^\8,\widetilde\pt)$ and suppose it is non-trivial. Then there will exist a non-split, complex supermanifold modelled on $(C, T^*_{\Xc^{\mathrm{an.}}, -})$, where $T^*_{\Xc^{\mathrm{an.}}, -}$ is as in (sub)Section \ref{9098g786fvuyviuvno}.
\end{THM}
\vspace{\dp0}
\end{minipage}
\end{center}
\emph{Proof of Theorem \ref{jf0jg43r3jgp3345topgj4gj4}.} 
Let $\Xc^\8 = (S(C, T^*_{C, -})^\8\times \Abb_\Cbb^{0|n}, \widetilde\pt_\xi)$ be an analytic deformation with $\widetilde\pt_\xi = \widetilde\pt - \vp(\xi)$, with $\vp(\xi)$ as in \eqref{tg74g94hf8hf033}. To second order, analytic deformation parameters describing $\Xc^\8$ can be represented by the $3$-tuple $(\chi^i, \chi^j, h^{ij})_{i<j;i, j = 1, \ldots, n}$. With $\nu = (\nu^i, \nu^j, \nu^{ij})_{i<j;i, j = 1, \ldots, n}$ representing a superconformal vector field to second order (see \eqref{fvnibviuh89rhf3j9jeoeio}), Theorem \ref{fhvurh84j093j3j33333} and \eqref{fbvrbvyurbvuir98h38h3} give the transformation rule:
\begin{align}
\left(
\begin{array}{c}
\chi^i
\\
\chi^j
\\
h^{ij}
\end{array}
\right)
\stackrel{\nu}{\longmapsto}
\left(
\begin{array}{l}
\chi^i + \widetilde\pt \nu^i
\\
\chi^j + \widetilde\pt \nu^j
\\
h^{ij} + \widetilde\pt \nu^{ij}
+\nu^i\chi^j - \nu^j\chi^i + \frac{1}{2}\big(\nu^i\widetilde\pt \nu^j - \nu^j\widetilde\pt\nu^i\big)
\end{array}
\right).
\label{rggh66u40u0940jf94f4}
\end{align}
An analytic deformation is \emph{trivial} if there exists some superconformal vector field $\nu = (\nu^i, \nu^j, \nu^{ij})_{i<j;i, j = 1, \ldots, n}$ sending  $(\chi^i, \chi^j, h^{ij})_{i<j;i, j = 1, \ldots, n}\mapsto (0, 0, 0)$. Now in order to prove Theorem \ref{jf0jg43r3jgp3345topgj4gj4} we will show that the atlas in Lemma \ref{jndrf4f4f4g4rnuir88y78reerrr} will define a non-vanishing, primary obstruction which, by Theorem \ref{34fiufguhfojp3}, defines a non-split supermanifold. Suppose firstly that $\Xc^\8$ is non-trivial to linear order, so that $(\chi^i)_{i =1, \ldots, n}$ transforms non-trivially. By \eqref{rggh66u40u0940jf94f4} this means $\chi^i$ is not exact. Recall that $\chi^i = \pt \vp/\pt\xi_i|_{\xi=0}$. Its preimage under the Dolbeault isomorphism will be non-trivial. By Lemma \ref{jndrnuir88y78reerrr}, its pre-image coincides with $\Theta(T^*_{\Xc^{\mathrm{an.}}, -})$, so therefore this extension class is non-trivial, i.e., $\Theta(T^*_{\Xc^{\mathrm{an.}}, -})\neq0$. Let $(\Uc, \vartheta)$ be the atlas from Lemma \ref{jndrf4f4f4g4rnuir88y78reerrr}. As it is superconformal to second order, its primary obstruction will map onto $\Theta(T^*_{\Xc^{\mathrm{an.}}, -})$ by Lemma \ref{hf8h8fh30f390jf93jf39}. As $\Theta(T^*_{\Xc^{\mathrm{an.}}, -})\neq0$ then $\om_{(\Uc, \vartheta)}\neq0$ so therefore $(\Uc, \vartheta)$ will define an atlas for some non-split supermanifold. Now suppose $\Xc^\8$ is trivial to linear order, but non-trivial to second order. This means $(\chi^i)_{i = 1,\ldots, n}$ is exact. From \eqref{rggh66u40u0940jf94f4} then we can deduce there will \emph{not} exist any $\nu$ sending $h^{ij}\mapsto 0$. From the correspondence between analytic and algebraic deformations in \eqref{fjbcevyu4iucneiocnjkebciue} we have:
\begin{align}
h^{ij} \longmapsto \left\{\left.\frac{\pt^2Dol^{-1}\vp}{\pt\xi_j\pt\xi_i}\right|_{\xi=0}\right\}_{\al,\be}
&\equiv
\left\{g^{ij}_{\al\be}\frac{\pt}{\pt y}\right\}_{\al, \be}
\label{fvkfvjkbjbruvbnvoirnvpe}
\\
&=
\left\{\left.\frac{\pt^2\om_{(\Uc, \vartheta)}}{\pt\xi_j\pt\xi_i}\right|_{\xi=0}\right\}_{\al, \be}
\label{rivruihurhj90j89fh78fghfo4344}
\end{align}
where \eqref{rivruihurhj90j89fh78fghfo4344} follows from the general expression for the primary obstruction to splitting in \eqref{nfrfuih98fh3h30n3} (see also \eqref{rg84gf873hf98h80j03f2}). Since the mapping in \eqref{fvkfvjkbjbruvbnvoirnvpe} preserves equivalences by Theorem \ref{rfg748gf7f983h04i4i4i43hf03}, there cannot exist any transformation sending $(g^{ij}_{\al\be})_{\al, \be}\mapsto 0$. Hence by \eqref{rivruihurhj90j89fh78fghfo4344} there cannot exist any transformation sending $\om_{(\Uc,\vartheta)}\mapsto0$ and so $\om_{(\Uc, \vartheta)}\neq 0$. Theorem \ref{jf0jg43r3jgp3345topgj4gj4} now follows. 
\qed
\\\\
We say an algebraic deformation $(\Uc, \vartheta)$ is \emph{non-split} if its primary obstruction to splitting is non-vanishing. Arguing similarly to Theorem \ref{jf0jg43r3jgp3345topgj4gj4} now, we have the converse:
\begin{center}
\begin{minipage}[t]{0.7\linewidth}
\begin{THM}\label{jf0jg43r3jgp33eefefefe4gj4} 
Let $(\Uc, \vartheta)$ be a non-split, algebraic deformation of a super Riemann surface $S(C, T^*_{C, -}) = (S(C, T^*_{C, -})^\8, \widetilde\pt)$. Then it defines a non-trivial, analytic deformation of $(S(C, T^*_{C, -})^\8, \widetilde\pt)$.
\qed
\end{THM}
\vspace{\dp0}
\end{minipage}
\end{center}

\subsection{The Associated Complex Superspace}
Recall that to any smooth supermanifold $\Xfr^\8$ with complex structure $\widetilde\pt$, there will be associated some complex supermanifold $\Xfr$. The structure sheaf $\Oc_\Xfr$ comprises (germs of) those functions on $\Xfr^\8$ in the kernel of $\widetilde\pt$. And so, to any analytic deformation $\Xc^\8 = \big(S(C, T^*_{C, -})^\8\times\Abb^{0|n}_\Cbb,\widetilde\pt_\xi\big)$, there is associated a complex supermanifold which we had referred to in (sub)Section \ref{9098g786fvuyviuvno} as $\Xc^{\mathrm{an.}}$. It is a complex supermanifold modelled on $(C, T^*_{\Xc^{\mathrm{an.}}, -})$ which, by Proposition \ref{jndrnuidddf222rreerrr}, is a superconformal model. As a locally ringed space,
\[
\Xc^{\mathrm{an.}} = \big(C, \Oc_{hol.}(\Xc^\8)\big)
\]
where the structure sheaf $\Oc_{hol.}(\Xc^\8)$ is the sheaf of holomorphic functions on $\Xc^\8$.\footnote{Recall, holomorphy is defined by reference to $\widetilde\pt_\xi$.} Now in Proposition \ref{jndrf4g4rnuir88y78reerrr} we deduced the existence of \emph{some} supermanifold $\Xc^\p$ modelled on $(C, T^*_{\Xc^{\mathrm{an.}}, -})$ which will admit a superconformal atlas, to second order, from the given, analytic deformation $\Xc^\8$. It is natural to then ask:
\begin{center}
\begin{minipage}[t]{0.7\linewidth}
~\\
{\bf Question.}
\emph{does the complex supermanifold $\Xc^{\mathrm{an.}}$ associated to the analytic deformation $\Xc^\8$ admit a superconformal atlas?}
\\
\vspace{\dp0}
\end{minipage}
\end{center}
We suspect the answer to the above question is `obviously true', however we will not venture to prove it here and leave it therefore as an open question. Note that as a corollary to an affirmative answer to the above question: $\Xc^\8$ will be non-trivial to second order if and only if $\Xc^{\mathrm{an.}}$ is non-split as a complex supermanifold by Theorem \ref{jf0jg43r3jgp3345topgj4gj4} and Theorem \ref{jf0jg43r3jgp33eefefefe4gj4}.

\section{Algebraic Deformations and Supermoduli}

\subsection{Statement of Result}
With the correspondence established in Theorem \ref{rfg748gf7f983h04i4i4i43hf03} we can readily address a question raised by the author in \cite{BETTSRS} on the relation between deformations of super Riemann surfaces and supermoduli. Denote by $\Mfr_g$ the supermoduli space of curves. This is the universal parameter space for super Riemann surfaces of genus $g$. It is a superspace in its own right; and as such Donagi and Witten in \cite{DW1, DW2} argued that its primary obstruction to splitting, denoted $\om_{\Mfr_g}$, is non-vanishing. In \cite{DW2} a dual description of $\om_{\Mfr_g}$ was offered by reference to analytic deformation parameters. That is, a  pairing was obtained between the primary obstruction $\om_{\Mfr_g}$ and analytic deformation parameters over $\Abb^{0|2}_\Cbb$, $(\chi^1, \chi^2, h^{12})$. It is invariant under equivalences of these parameters. This pairing can be generalised to analytic deformation parameters over $\Abb^{0|n}_\Cbb$, $(\chi^i, \chi^{j}, h^{ij})_{i<j = 0,\ldots, n}$ (see Lemma \ref{jfdd4sssso4jopgj4gj4} below). From the correspondence presented in Theorem \ref{rfg748gf7f983h04i4i4i43hf03} between (equivalence classes of) analytic and algebraic deformations, we can then conclude:
\begin{center}
\begin{minipage}[t]{0.7\linewidth}
\begin{THM}\label{jf04jg0jg4jgpo4jopgj4gj4} 
For any algebraic deformation $(\Uc, \vartheta)$ of a super Riemann surface $S(C, T^*_{C, -})$ there exists an invariant pairing between $\om_{\Mfr_{g}}$ and $\om_{(\Uc, \vartheta)}$.
\end{THM}
\vspace{\dp0}
\end{minipage}
\end{center}
Invariance of the pairing in Theorem \ref{jf04jg0jg4jgpo4jopgj4gj4} means $\langle \om_{\Mfr_g}, \om_{(\Uc, \vartheta)}\rangle = \langle \om_{\Mfr_g}, \om_{(\tilde\Uc, \tilde\vartheta)}\rangle$ for any two equivalent algebraic deformations $(\Uc, \vartheta)\equiv (\tilde\Uc, \tilde\vartheta)$. Before we give a proof of Theorem \ref{jf04jg0jg4jgpo4jopgj4gj4} it will be necessary to digress briefly on the pairing between analytic deformations and supermoduli obstructions.

\subsection{Donagi and Witten's Invariant Pairing}
In \cite{DW2}, en route to their deduction of non-projectedness of the supermoduli space, the following result is proved:
\begin{center}
\begin{minipage}[t]{0.7\linewidth}
\begin{THM}\label{jf04sssso4jopgj4gj4} 
Let $\om_{\Mfr_g}$ denote the primary obstruction to splitting the supermoduli space of genus $g$ curves $\Mfr_g$; and let $(\Xc^\8, \widetilde \pt_\xi)$ be an analytic family of super Riemann surfaces over $\Abb^{0|2}_\Cbb$. Then there exists an invariant pairing between $\om_{\Mfr_g}$ and $\widetilde\pt_\xi$.
\qed
\end{THM}
\vspace{\dp0}
\end{minipage}
\end{center}
By invariance of the pairing $\langle\om_{\Mfr_g},\widetilde\pt_\xi\rangle$ stated in Theorem \ref{jf04sssso4jopgj4gj4} it is meant that: $\langle\om_{\Mfr_g},\widetilde\pt_\xi\rangle = \langle\om_{\Mfr_g},\widetilde\pt^\p\rangle$ for any two, linearly equivalent, analytic  families.\footnote{recall the terminology from Theorem \ref{fhvurh84j093j3j33333}.} With $(S(C, T^*_{C, -})^\8, \widetilde\pt)$ the central fiber of $(\Xc^\8, \widetilde\pt_\xi)$; $\widetilde\pt_\xi = \widetilde\pt -\vp(\xi)$ and $\vp(\xi)$ as in \eqref{tg74g94hf8hf033}, the following local formula\footnote{up to a sign factor, depending on the conventions adopted here} was provided in \cite[p. 32]{DW2},
\begin{align}
\Big\langle\om_{\Mfr_g},&\left.\widetilde\pt_\xi\Big\rangle\right|_{\left(S(C, T^*_{C, -})^\8, \widetilde\pt\right)}
\notag
\\
\stackrel{\Delta}{=}&~
\int_{C\times C}
\left\langle
\om_{\Mfr_g},\chi^1\boxtimes\chi^2
\right\rangle
-
2\pi \sqrt{-1}
\int_C
\left\langle \mathrm{Res}~\om_{\Mfr_g}
,
h^{12}
\right\rangle,
\label{rfg84gf7hf83f903h98f38fh03}
\end{align}
where the expression in \eqref{rfg84gf7hf83f903h98f38fh03} is defined via Serre duality. We refer to \cite{DW2} for further elaboration on this point. 

\subsection{Proof of Theorem $\ref{jf04jg0jg4jgpo4jopgj4gj4}$}
We begin with the following generalisation of Donagi and Witten's invariant pairing from Theorem \ref{jf04sssso4jopgj4gj4}:
\begin{center}
\begin{minipage}[t]{0.7\linewidth}
\begin{LEM}\label{jfdd4sssso4jopgj4gj4} 
The invariant pairing in Theorem \ref{jf04sssso4jopgj4gj4} holds more generally for any analytic family over $\Abb^{0|n}_\Cbb$, where $n\geq2$. It is given by the local expression:
\begin{align*}
\Big\langle\om_{\Mfr_g},&\left.\widetilde\pt_\xi\Big\rangle\right|_{\left(S(C, T^*_{C, -})^\8, \widetilde\pt\right)}
\\
\stackrel{\Delta}{=}&~
\sum_{i<j}
\left(\int_{C\times C}
\left\langle
\om_{\Mfr_g},\chi^i\boxtimes\chi^j
\right\rangle
-
2\pi \sqrt{-1}
\int_C
\left\langle \mathrm{Res}~\om_{\Mfr_g}
,
h^{ij}
\right\rangle
\right).
\end{align*}
\end{LEM}
\vspace{\dp0}
\end{minipage}
\end{center}
\emph{Proof of Lemma \ref{jfdd4sssso4jopgj4gj4}.}
Each summand on the right-hand side of the local expression for $\langle\om_{\Mfr_g},\widetilde\pt_\xi\rangle$ in the statement of this lemma will be invariant under linear equivalence in Theorem \ref{fhvurh84j093j3j33333}. This can be deduced by the same arguments given by Donagi and Witten in \cite{DW2} regarding the expression for pairings of deformations over $\Abb^{0|2}_\Cbb$ in \eqref{rfg84gf7hf83f903h98f38fh03}. Hence the entire expression, i.e., the sum over $\{i, j: i<j\}$, will be invariant.
\qed
\\\\
Now recall the relation between analytic and algebraic deformations in \eqref{rgf64gf783hf983hf08j3093}. For any algebraic deformation $(\Uc, \vartheta)$ of a super Riemann surface $S(C, T^*_{C, -})\equiv (S(C, T^*_{C, -})^\8, \widetilde\pt)$ over $\Abb^{0|n}_\Cbb$ with primary obstruction $\om_{(\Uc, \vartheta)}$, write:
\[
\widetilde\pt_{(\Uc, \vartheta)}
\stackrel{\Delta}{=}
\widetilde\pt 
-
\sum_i \xi_i\left.\frac{\pt Dol~\om_{(\Uc, \vartheta)}}{\pt \xi_i}\right|_{\xi=0}
-
\sum_{i<j} \xi_i\xi_j\left.\frac{\pt^2Dol~\om_{(\Uc, \vartheta)}}{\pt \xi_j\pt\xi_i}\right|_{\xi=0}.
\]
By Proposition \ref{hg975hg984hg8003jg03}, $\widetilde\pt_{(\Uc, \vartheta)}$ will define an analytic deformation and we can set:
\[
\left.\left\langle \om_{\Mfr_g}
, 
\om_{(\Uc,\vartheta)}\right\rangle\right|_{S(C, T^*_{C, -})}
\stackrel{\Delta}{=}
\left.\left\langle 
\om_{\Mfr_g}
,
Dol~\om_{(\Uc, \vartheta)}
\right\rangle
\right|_{(S(C, T^*_{C,-})^\8, \widetilde\pt)}.
\]
By Theorem \ref{rfg748gf7f983h04i4i4i43hf03}, equivalences of algebraic deformations translate to equivalences of analytic deformations and so to linear equivalences. As such, the expression $\left\langle \om_{\Mfr_g}, \om_{(\Uc,\vartheta)}\right\rangle$ will be invariant under equivalences of algebraic deformations. Theorem \ref{jf04jg0jg4jgpo4jopgj4gj4} now follows.\qed

\section{The No-Exotic-Obstructions Conjecture}
\label{rhf47g794hg893093j9444}

\noindent
In (sub)Section \ref{rhg74hg89h03j09g408} we presented a conjecture from a previous article by the author. It concerned the higher obstructions to splitting the total space of a super Riemann surface deformation. We recall this conjecture below for convenience:
\begin{center}
\begin{minipage}[t]{0.7\linewidth}
\begin{CONJ1}
Let $(\Uc, \vartheta)$ be an algebraic deformation of a super Riemann surface $S(C, T^*_{C, -})$ over $\Abb^{0|n}_\Cbb$. If its primary obstruction $\om_{(\Uc, \vartheta)}$ vanishes, then $(\Uc, \vartheta)$ is split.
\end{CONJ1}
\vspace{\dp0}
\end{minipage}
\end{center}
In what follows we sketch a potential route to proving the above conjecture. Following \cite{DW1, BETTHIGHOBS}, an atlas for a supermanifold $\Xc$ is said to be \emph{exotic} if:
\begin{enumerate}[$\bt$]
	\item $\Xc$ is split as a supermanifold and;
	\item the atlas defines a non-vanishing obstruction to splitting.
\end{enumerate}
From Theorem \ref{34fiufguhfojp3}, any atlas with non-vanishing primary obstruction will be non-split and so cannot be exotic. Therefore, any exotic atlas will define a `higher' obstruction to splitting. The above conjecture  asserts: there do not exist any exotic, algebraic deformations of a super Riemann surface. Hence, any obstruction to splitting the total space of an algebraic deformation must be primary or must vanish. This statement was confirmed explicitly for deformations over $\Abb^{0|2}_\Cbb$. In the general setting, equipped now with Theorem \ref{jf04jg0jg4jgpo4jopgj4gj4}, a proof of the above conjecture might read along the following lines: 
\begin{center}
\begin{minipage}[t]{0.7\linewidth}
\emph{Proposed Proof Sketch of Conjecture.}
Suppose $\langle \om_{\Mfr_g}, -\rangle$ is non-vanishing on a subset $\Uscr\subseteq \Mfr_g$. This means at any super Riemann surface $[S(C, T^*_{C, -})]\in \Uscr$ and for any algebraic deformation $(\Uc, \vartheta)$ that $\langle \om_{\Mfr_g}, \om_{(\Uc, \vartheta)}\rangle = 0$ if and only if $\om_{(\Uc, \vartheta)} =0$. We suppose moreover that $\om_{\Mfr_g}$ will invariantly pair with \emph{any} obstruction to splitting an algebraic deformation $(\Uc, \vartheta)$. Assuming now that $(\Uc, \vartheta)$ is exotic, we have $\om_{(\Uc, \vartheta)}^{\mathrm{primary}} = 0$. Hence $\langle \om_{\Mfr_g}, \om_{(\Uc, \vartheta)}^{\mathrm{primary}}\rangle = 0$. Since this pairing is invariant under equivalences by the second of our suppositions, we must also have $\langle \om_{\Mfr_g}, \om_{(\Uc, \vartheta)}\rangle =0$. Therefore $\om_{(\Uc, \vartheta)} = 0$ by the first supposition. Hence $(\Uc, \vartheta)$ must be split. This argument reveals, over regions $\Uscr\subseteq \Mfr_g$ where $\langle\om_{\Mfr_g}, -\rangle$ is non-vanishing, that there cannot exist any exotic, algebraic deformations. The conjecture would therefore hold over these regions in supermoduli space. It remains then to argue $\Uscr = \Mfr_g$.
\qed
\vspace{\dp0}
\end{minipage}
\end{center}
~\\



\addtocontents{toc}{\vspace{\normalbaselineskip}}
\bibliographystyle{alpha}
\bibliography{Bibliography}

\hfill
\\
\noindent
\small
\textsc{
Kowshik Bettadapura 
\\
\emph{Yau Mathematical Sciences Center} 
\\
Tsinghua University
\\
Beijing, 100084, China}
\\
\emph{E-mail address:} \href{mailto:kowshik@mail.tsinghua.edu.cn}{kowshik@mail.tsinghua.edu.cn}

\end{document}